\documentclass[mnsc,nonblindrev]{informs3_nojournal}

\usepackage{natbib}
 \bibpunct[, ]{(}{)}{,}{a}{}{,}%
 \def\bibfont{\small}%
\TheoremsNumberedThrough     
\ECRepeatTheorems

\EquationsNumberedThrough    

\MANUSCRIPTNO{}

\usepackage[T1]{fontenc}
\usepackage[utf8]{inputenc}
\usepackage{xcolor}
\usepackage{soul}
\usepackage{subfigure}
\usepackage{placeins}
\usepackage[colorlinks,citecolor=black,linkcolor=black,bookmarks=false,hypertexnames=true]{hyperref} 
\usepackage{algorithm}
\usepackage{algorithmic}

\usepackage{thmtools,mathtools}

\usepackage[margin=1in]{geometry}
\usepackage{graphicx}

\usepackage{soul, multirow, bm, amssymb, mathtools, array, enumitem}

\newcommand{\PP}{\mathbb{P}}

\newcommand{\ones}{\bm{1}}

\newcommand{\prodset}{\mathcal{N}}

\newcommand{\ropq}{\rho}

\newcommand{\trad}[2]{\textsc{Trad}(#1, #2)}
\newcommand{\rs}[1]{\bm{r}_{#1}}

\newcommand{\ctrad}[2]{\pi \rpr{#1 \mid \trad{#2}{\rs{#2}}}}
\newcommand{\ctradequal}[3]{\pi \rpr{#1 \mid \trad{#2}{#3}}}

\newcommand{\Snp}{S \cup \cpr{0}}
\newcommand{\rev}[1]{\textsc{Rev}\rpr{#1}}

\let\originalleft\left
\let\originalright\right
\renewcommand{\left}{\mathopen{}\mathclose\bgroup\originalleft}
\renewcommand{\right}{\aftergroup\egroup\originalright}

\newcommand{\cpr}[1]{\left\{#1\right\}}
\newcommand{\spr}[1]{\left[#1\right]}
\newcommand{\rpr}[1]{\left(#1\right)}

\usepackage{xcolor}

\newif\ifrevision
\revisiontrue   %
\revisionfalse %

\ifrevision
  \newcommand{\cor}[1]{\textcolor{blue}{#1}}
\else
  \newcommand{\cor}[1]{#1}
\fi

\newcounter{omarcounter}

\begin{document}
    
    \RUNAUTHOR{El Housni, Elmachtoub, Sheth, and Shi}
     
    \RUNTITLE{Price and Assortment Optimization with Opaque Products}
    
    \TITLE{Price and Assortment Optimization under the Multinomial Logit Model with Opaque Products}
    
    \ARTICLEAUTHORS{%
    \AUTHOR{Omar El Housni}
    \AFF{School of Operations Research and Information Engineering, Cornell Tech, Cornell University, New York, NY 10044, \EMAIL{ oe46@cornell.edu}}
    \AUTHOR{Adam N. Elmachtoub}
    \AFF{Department of Industrial Engineering and Operations Research and Data Science Institute, Columbia University, New York, NY 10027, \EMAIL{adam@ieor.columbia.edu}}
    \AUTHOR{Harsh Sheth}
    \AFF{Department of Industrial Engineering and Operations Research and Data Science Institute, Columbia University,  New York, NY 10027, \EMAIL{hts2112@columbia.edu} }
    \AUTHOR{Jiaqi Shi}
    \AFF{Department of Industrial Engineering and Operations Research and Data Science Institute, Columbia University,  New York, NY 10027, \EMAIL{js5778@columbia.edu} }
    } 
    
    \ABSTRACT{%
    \noindent An opaque product is a product for which only partial information is disclosed to the buyer at the time of purchase. Opaque products are common in sectors such as travel and online retail, where the brand or product color is hidden in the opaque product. Opaque products enable sellers to target customers who prefer a price discount in exchange for being flexible about the product they receive. In this paper, we integrate opaque products and traditional products together into the multinomial logit (MNL) choice model and study the associated price and assortment optimization problems. For the price optimization problem, we show a surprising result that uniform pricing is optimal which implies it has the same optimal pricing solution
    and value as the traditional MNL model. Although adding an opaque product may enhance revenue given arbitrary traditional product prices, our result shows that this advantage disappears when all prices are optimized jointly.  \cor{For the assortment problem, we first study the setting where the seller chooses both the offered assortment and the support of the opaque product, and prove that there exists an optimal assortment that is nested by revenue.} We then study the case where the opaque support is required to coincide with the offered assortment which is common in practice. Here, we show that the optimal assortment is nested-by-valuation for uniformly priced products. For arbitrary prices, we propose a nested-by-revenue-and-valuation heuristic with a theoretical approximation guarantee that performs extremely well in our numerical study. 
    }
    \maketitle

    \section{Introduction}
\begin{figure}[h!] 
\caption{Examples of opaque products.}
\label{fig:examples}	
  \centering
\subfigure[Opaque product -- ``Managers Special'' on Hertz.]{\includegraphics[height=0.3\textwidth]{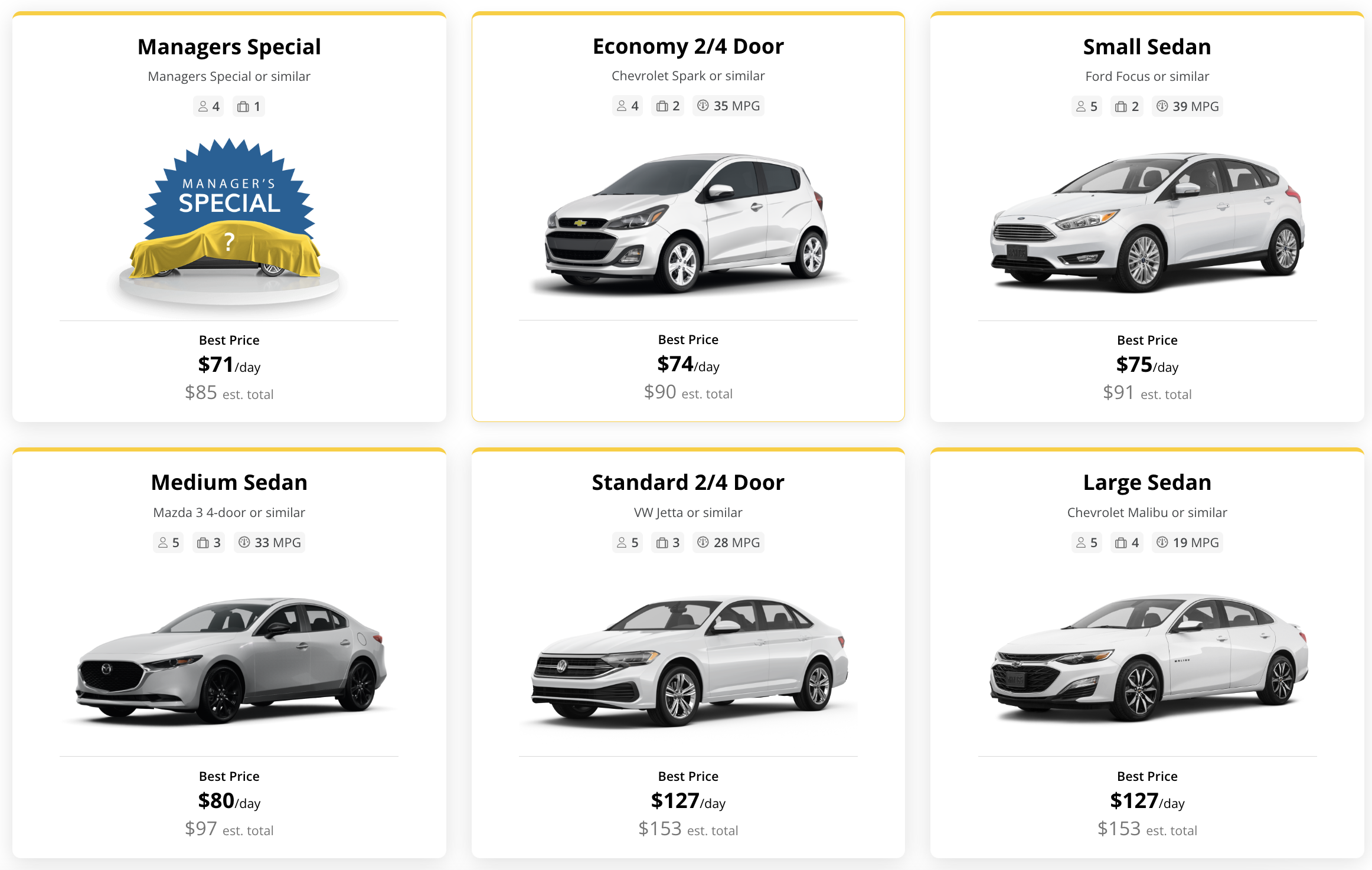}  \label{fig:car}   }  
\hfill
\subfigure[Opaque product -- ``Color Will Vary'' on Amazon. ]{\includegraphics[height=0.3\textwidth]{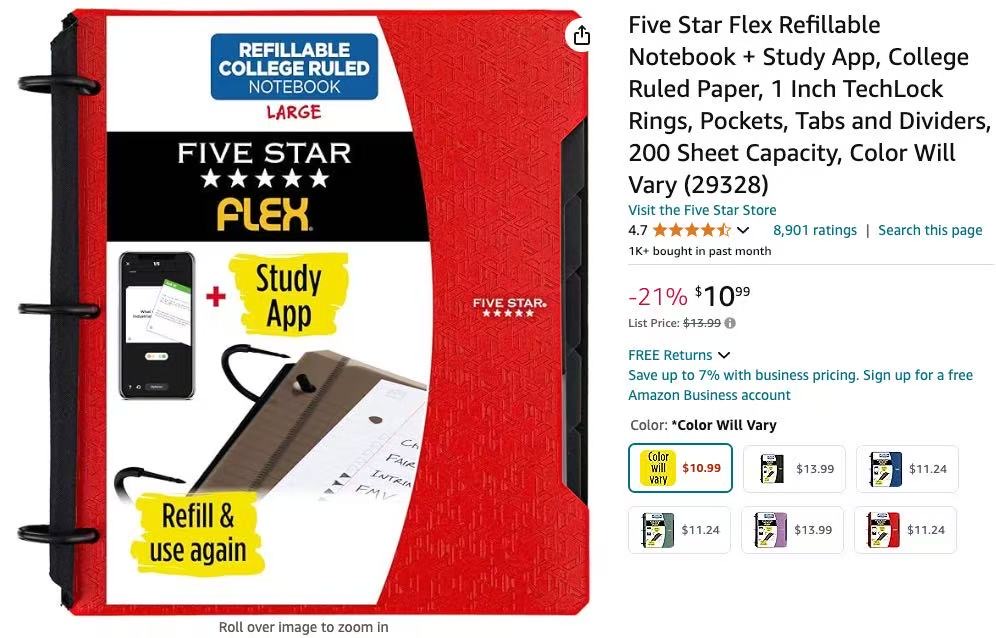} \label{fig:amazon}  } 

\vspace{0.2cm}
\begin{minipage}{\textwidth}
    \footnotesize
    \textit{Note.} (a) A search for rental cars at JFK on Hertz shows the following rates: Economy at \$74, Small Sedan at \$75, Medium Sedan at \$80, Standard at \$127, and Large Sedan at \$127. In contrast, an opaque car—where the company selects the vehicle—is available for \$71. (b) The original price on Amazon is \$13.99 for black and purple. However, other colors like blue, red, and green are offered at a discount of \$2.75, while the opaque product—where the specific color is not disclosed to the customer—receives a higher discount of \$3.
    \end{minipage}
\end{figure}

    An opaque product is a product for which only partial information is disclosed to the buyer at the time of purchase. The practice of selling opaque products has gained traction in industries such as travel (including airlines, hospitality, and car rentals) and online retail.
    For instance, the car rental company Hertz offers a ``Managers Special'' option, where the car type is hidden until after the booking is confirmed (see Figure~\ref{fig:examples}\subref{fig:car}). In another example, Amazon.com offers notebooks in various colors alongside a ``Color Will Vary'' option at a discounted price (see Figure~\ref{fig:examples}\subref{fig:amazon}).  In both examples, customers are faced with the choice between selecting the traditional products or opting for the opaque product at a lower price.

    Offering opaque products provides two primary advantages for sellers. First, it allows sellers to potentially increase revenue by selling to flexible customers at a discount without overly cannibalizing demand for traditional products \citep{mike2021power}. Second, it provides the seller the flexibility to manage their product inventories more efficiently \citep{yehua2015retailing}. To fully leverage these benefits, it is also important to (when possible) (i) optimize the prices of offered products (e.g., determining the best prices for the ``Color Will Vary'' option and other specific colors as in the Amazon example in Figure~\ref{fig:examples}\subref{fig:amazon}), and (ii) select the appropriate product assortment, for instance, deciding whether to offer all available colors or only a subset. Choice models, which effectively capture customer choice behavior between substitutable products, have been identified as powerful tools to enhance revenue by solving associated price and assortment optimization problems \citep{talluri2004revenue,vulcano2010om}. In this paper, we integrate opaque products alongside traditional products into the classical multinomial logit (MNL) choice model. We study two important optimization problems under this model: price optimization and assortment optimization.
    
    In our framework, we consider a seller managing a set of substitutable products that are similar with the same cost but differentiated by some feature such as color or brand. We refer to these products as \emph{traditional products}. The seller may choose an assortment of these products to offer to customers, with prices that may be either exogenously set or endogenously determined. In addition, the seller offers an \emph{opaque product} at a seller-chosen (discounted) price, where the customer purchases without knowing the exact traditional product they will receive at the time of purchase. \cor{Associated with the opaque product is a known \emph{opaque support}, defined as a subset of the offered assortment from which the realized product is selected. Hence, the opaque product provides the seller with flexibility to satisfy demand using any traditional product in the opaque support.} Customers draw valuations for each of the traditional products from independent Gumbel distributions (following the MNL framework), and select the product with the highest utility (valuation minus price). Given these valuation distributions, we refer to the model where only traditional products are offered as the \emph{traditional MNL model}, and the model where both traditional products and the opaque product are offered as the \textit{Opaque-MNL model}. 

    A key issue to address in the Opaque-MNL model is the customer's valuation for the opaque product.  Early works on ``probabilistic selling'', such as \cite{fay2008probabilistic}, assume that the customer is aware of the probability of receiving every possible traditional product when purchasing the opaque product, leading to a valuation of the opaque product based on the expected value over the possible alternatives. However, this approach does not apply to the scenarios illustrated in Figures \ref{fig:examples}\subref{fig:car} and \ref{fig:examples}\subref{fig:amazon}, where the allocation probabilities are unknown to customers.  \cor{In this work, we consider $\bm{k}$\textbf{-risk-averse} customers who value the opaque product as the $k$-th \textit{smallest} order statistic of their valuations for the traditional products in the opaque support. Here $k$ is a parameter that can take any value between 1 and the size of the assortment. For example, when $k =1$, this captures the original \textit{risk-averse} assumption proposed in \cite{mike2021power}, wherein customers value the opaque product as the minimum valuation among all traditional products that can fulfill it.  Such customers are considered risk-averse because they will not regret their purchase decision even if the opaque product is fulfilled with the least preferred product. This assumption is particularly valid when customers inherently distrust the seller, lack prior information about which product will fulfill their purchase, or recognize that certain products provide no value. When $k$ represents the median, it signifies the customer believes they will receive a product somewhere near the middle of their preference list. For larger $k$, the opaque product is valued more highly, which captures customers with a more optimistic or gambling-oriented attitude toward opaque purchases.  We next summarize our contributions of this paper.}

    \textbf{Price Optimization.} We study the problem of determining the revenue-maximizing prices for traditional products and the opaque product under the Opaque-MNL model. Recall that for the traditional MNL model, it is well established that setting all products at the same price is optimal (see, e.g., \citealt{hopp2005product}).
    Surprisingly, \cor{for any $k$-risk-averse customer}, we discover that uniform pricing remains optimal for the Opaque-MNL model (Theorem \ref{thm-pricing}), yielding the same optimal pricing solution and value as the traditional MNL model. While adding an opaque product priced at the revenue-maximizing opaque price can enhance the revenue under the traditional MNL model in general, this is no longer the case when traditional prices are jointly optimized rather than provided exogenously. In such scenarios, introducing an opaque product does not enhance revenue. 
    \cor{To prove Theorem \ref{thm-pricing}, we overcome the analytical challenges posed by the exponential terms and the non-linearity of the $\min$ operator in the revenue expression. Our approach relies on two key technical steps. First, we simplify the opaque revenue expression by showing that the inclusion-exclusion summation can be restricted solely to the subset of products priced higher than the opaque product, effectively resolving the complexity of the minimization term. Second, using a change of variables, we decompose the opaque revenue into a fixed traditional revenue component and a remainder term governed by a concave auxiliary function. By exploiting this concavity, we establish that the revenue under the Opaque-MNL model is always bounded above by the revenue of the traditional MNL model at the optimal uniform price, thereby proving that uniform pricing remains optimal.}


    \textbf{Assortment Optimization.} We examine the problem of identifying the revenue-maximizing assortment under the Opaque-MNL model, which is important when the seller cannot change the prices of the traditional products. \cor{This is common for third-party travel websites that sell car rentals and airfare, and common when selling retail products where traditional prices are set by the manufacturer. We first study the most flexible design setting, in which the seller jointly chooses the offered assortment and the opaque support, i.e., the set of traditional products that may be used to fulfill the opaque product. The joint problem is challenging because the seller must simultaneously optimize the assortment, the opaque support, and the opaque price, which results in the revenue function containing exponentially many terms. Despite these difficulties, we prove that the optimal assortment is nested-by-revenue (Theorem \ref{Thm: revenue-ordered}).  
    }
    
    We then study the assortment problem in the setting where the opaque support is required to coincide with the offered assortment. This setting is important when the seller cannot choose the opaque support separately, which is common in practice since it is easier to implement.
    Once this flexibility is removed, the problem becomes substantially harder. In particular, the nested-by-revenue structure of the optimal assortment fails, and several key properties of the traditional MNL model, such as substitutability and inclusion of the highest-revenue product in the optimal assortment, do not extend to the Opaque-MNL model. 
    Nevertheless, we prove that the revenue-maximizing assortment is nested-by-valuation when all traditional products share the same price (Theorem \ref{thm:nested-val}). This case, where products are uniformly priced, is particularly relevant as it reflects a natural setting for selling opaque products when traditional products differ by only one feature (such as color). For the more general case of non-uniformly priced products, we propose a natural nested-by-revenue-and-valuation heuristic, drawing on properties from the traditional MNL model and Opaque-MNL model under uniform prices. We prove that this heuristic offers a 1/2-approximation guarantee (Theorem \ref{prop: NRV}), and demonstrate that it performs extremely well in numerical experiments.

    \subsection{Literature Review}
    \textbf{Price optimization under the MNL model.} 
     Under the traditional multinomial logit (MNL) model, it is well-known that the profit function of multiple substitutable products is concave in the choice probabilities \citep{song2021demand,dong2009dynamic}. \cite{hopp2005product} and \cite{wang2012capacitated} show that the markup or adjusted markup is product-invariant at optimality under the MNL model. Our review is limited to studies based on the MNL model, but we refer readers to other examples that involve the constant markup property of the optimal solution including \cite{li2011pricing} and \cite{gallego2014multiproduct} under the nested logit model, \cite{li2015d} under the $d$-level nested logit model, \cite{li2017optimal} under the paired combinatorial logit (PCL) model, and \cite{zhang2018multiproduct} under the generalized extreme value (GEV) models.

    \textbf{Assortment optimization under the MNL model.} 
      Under the traditional MNL model, \cite{talluri2004revenue} show that the revenue-maximizing assortment is nested by revenue. \cite{rusmevichientong2010dynamic},  \cite{sumida2021revenue}, \cite{elmachtoub2022revenue} study  MNL assortment problems under cardinality constraints, totally unimodular constraints, and product retirement constraints,  respectively. \cite{ma2023assortment}  studies when deterministic assortments are optimal among broader randomized mechanisms and shows that this holds under several common choice models, including MNL and Markov Chain models.     
      \cite{bai2023assortment} and \cite{bai2024assortment} consider variants of the MNL model, focusing on the endogenous context-dependent MNL and the utility-based rank cutoffs MNL, respectively. 
      \cite{bront2009column,rusmevichientong2014assortment} and \cite{el2023joint} study the assortment problem under a mixed of multinomial logit models. For other representative work on assortment optimization under other choice models, we refer to \cite{davis2014assortment} for the nested logit model, \cite{blanchet2016markov} for the markov chain model, and \cite{aouad2018approximability} for the preference list-based model.

    \textbf{Opaque selling strategy.}
    The work by \citet{mike2021power} is closely related to ours in the existing literature. Their work focuses on identifying conditions under which the opaque selling strategy dominates discriminatory pricing and single pricing strategies for exchangeable valuation distributions. In contrast, our work focuses on price and assortment optimization problems for $k$-risk-averse customers under the multinomial logit choice model with opaque products. In addition, their research examines both risk-averse and risk-neutral customers. \cor{Subsequently to our paper, \cite{chen2025approximation} study a multi-period, multi-product dynamic pricing problem for risk-averse customers with opaque products. Their model adopts an opaque-specific independent-demand framework, in which each customer arrives with a focal traditional product in mind and compares it only with the opaque product when that product is included in the opaque support. In contrast, we study a static MNL-based model for $k$-risk-averse customers that consider purchasing every traditional product offered alongside the opaque product. 
    }

    Opaque selling for risk-neutral customers has also been extensively studied under the name of probabilistic selling in the literature \citep{jiang2007price, fay2008probabilistic,  jerath2010revenue, zhang2015verticaldiff}. These papers typically focus on settings involving two or three products, with customer valuations being drawn from the Hotelling or Salop's circle choice models, demonstrating that there is value in offering opaque products when  different customer segments value products differently. \cite{anderson2012choice} consider fitting a nested logit choice model using star level and neighborhood as features, based on a dataset from Hotwire. They apply dynamic programming to characterize the optimal pricing policy for the opaque hotel bookings. \citet{strauss2021dynamic} also consider the pricing problem with opaque products, particularly in the setting of flexible time slots. They assume a nested logit model for customer valuations, with traditional and opaque products belonging in separate nests. Other works in probabilistic selling have explored diverse scenarios, such as bounded rational customers \citep{huang2021dynamic,huang2014sell}, customers' salient thinking behavior \citep{zheng2019probabilistic}, and the impact of inventory information disclosure \citep{li2020mitigating}.  Other works quantify the value of opaque products in inventory management and supply chain \citep{fay2015timing, yehua2015retailing, ren2022opaque,  zhang2024less, fay2024effect, freund2025power,elmachtoub2026value}.

    Opaque products have also been widely studied in the revenue management literature under the term ``flexible products". \citet{gallegofirst2004revenue} study a model with two traditional products and one flexible product, where the allocation of the flexible product occurs at the end of the selling horizon. \cite{gallegosecond2004managing} extend the concept of flexible products to network revenue management. 
    They also consider choice models in their work and provide algorithms for assortment and allocation planning with asymptotic guarantees when the customer valuations are drawn from general ``attraction'' models (which encapsulate the MNL choice model). Their work differs from ours as they do not assume the opaque or flexible product's valuation is derived from those of the traditional products, and instead allow for it to have an arbitrary random valuation. Similarly, \citet{liu2022opaque} also provide an algorithm for assortment and allocation planning with finite inventories based on the fluid relaxation of the problem. The fluid problem requires known choice probabilities to satisfy the regularity assumptions, which are violated for our model.  \cite{piri2023firm} study firm-centric and customer-centric flexibility products and show that integrating the two can create a synergistic profit effect. \cite{zhu2024performance} use dynamic programming to propose a policy with a constant factor approximation ratio for the network revenue management problems with flexible products.


    \section{Model Formulation}   We consider a seller with a set $\prodset := \cpr{1,\ldots, n}$ of substitutable products, where $\bm r = \cpr{r_i \mid i \in \prodset}$ denotes the vector of product prices. The seller may choose an assortment $S \subseteq \cal N$ to offer, along with an opaque product, denoted by $q$, at price $\rho$. We refer to these $n$ products as traditional products to distinguish them from the opaque product $q$. In addition, the no-purchase option is represented as product 0. The customer has a valuation $V_i$ for each product $i \in \cal N \cup \{0\}$ and a valuation $V_q$ for the opaque product. Following the classical MNL framework, we assume that $V_i$ has the Gumbel distribution with intrinsic valuation $v_i$, given by:
     \begin{align}\label{eq: gumbel}
         V_i := v_i + \epsilon_i,  \ i \in \cal N \quad \text{and} \quad V_0 = \epsilon_0,
     \end{align}
    where $(\epsilon_i)_{i=0}^n$ are i.i.d. random variables with distribution $Gumbel(0,1)$. The exact distribution of $V_q$ will be discussed in detail later. For a traditional product $i \in \cal N$, the utility of $i$ is defined as the valuation minus price, i.e.,  $U_i := V_i -r_i $. Similarly,  $U_q := V_q - \rho$ for the opaque product, and $U_0 :=V_0  $ for the no-purchase option.

    In this work, we consider two selling models: (i) traditional MNL model - involving only the sale of traditional products, denoted by $\textsc{Trad}$ and (ii) Opaque-MNL - involving the sale of traditional products alongside an opaque product, denoted by $\textsc{Opq}$. For every selling model, the customer maximizes their own utility among all available options (including the no-purchase option) and selects the one with the highest utility, consistent with the Random Utility Maximization (RUM) framework. We next describe the choice probabilities and expected revenues for each selling model.

    \subsection{Traditional MNL Model}
    We first consider the model where only traditional products are offered. 
    Fix any assortment $S \subseteq \cal N $. Let $\bm{r}_S := \{r_i \mid i \in S\}$ denote the corresponding price vector. To simplify notation, we will refer to this selling model as $\trad{S}{\bm r_S}$ hereafter. According to the Gumbel distribution assumption in \eqref{eq: gumbel}, the classical MNL model yields that the purchase probability of any product $i$ under $\trad{S}{\bm r_S}$ is given by
    \begin{equation} \label{eq:mnl}
		\ctrad{i}{S} := \PP\spr{U_i > U_j,  \forall j \in \Snp \setminus \cpr{i} } = \frac{e^{v_i - r_i}}{1 + \sum_{j \in S} e^{v_j - r_j}}.
    \end{equation}
    Naturally, we have $\ctrad{i}{S} = 0$ for $ i \notin S $. Further, we let $\rev{\trad{S}{\rs{S}}} $ denote
    the expected revenue under selling model $\trad{S}{\bm r_S}$, which is expressed as:
    \begin{equation} \label{eq:mnl-revenue}
	\rev{\trad{S}{\rs{S}}} = \sum_{i \in S} r_i \ctrad{i}{S} =  \frac{\sum_{i \in S} r_i e^{v_i - r_i}}{1 + \sum_{i \in S} e^{v_i - r_i}}.
    \end{equation}

\newcommand{\Low}{L}
\newcommand{\coeff}{(-1)^{|I|-k}\binom{|I|-1}{k-1}}
\newcommand{\Trad}{\textsc{Trad}}
\newcommand{\Opq}{\textsc{Opq}}
\newcommand{\Rev}{\textsc{Rev}}

\subsection{Opaque-MNL Model}
    The seller may also choose to offer an opaque product in addition to the traditional products. Specifically, the seller provides an assortment $S \subseteq \cal N$ of traditional products and an opaque product $q$. \cor{The seller also needs to choose the opaque support $T \subseteq S $, meaning that whenever the opaque product is purchased,  the seller fulfills it using one of the traditional products in $T$. The key to deriving the choice probabilities under the Opaque-MNL model lies in the definition of the opaque product's valuation $V_q$. In this paper, we consider \textit{$k$-risk-averse} customers who are characterized by a fixed parameter $k\in\{1,2,\ldots,n\}$, where $k$ reflects the customer’s attitude toward the uncertainty induced by opacity and is independent of the offered assortment and opaque support. Given an opaque support $T$, such customers value the opaque product as the $k$-th \textit{smallest} valuation among the traditional products in opaque support $T$. Formally, for any $T \subseteq\mathcal N$,
    \begin{align*}
        V_q^T := V_{(k)}^T,
    \end{align*}
    where $V_{(k)}^T$ denotes the $k$-th order statistic of $\{V_j:j\in T\}$.
    If $k>|T|$, we adopt the convention that
    \begin{align*}
        V_{(k)}^T := \max_{i\in T} V_i .
    \end{align*} 
    In the following, we refer to this selling model with price vector $\bm r_S$ of traditional products, opaque support $T \subseteq S $, and opaque price $\rho$ as $\Opq_{(k)}(T, S,\bm r_S,\rho)$.
    In particular, when $k=1$, we recover the risk-averse customer setting studied in \citealt{mike2021power}, who value the opaque product as the minimum valuation among all traditional products in the opaque support, i.e.,  $ V^T_q = \min_{i \in T} V_i = V_{(1)}^T$. Those customers are considered risk-averse in the sense that they will not regret their purchase decision even if the opaque product is fulfilled with the least preferred product. At the other extreme, when $k \ge |T|$, customers are optimistic and value the opaque product as their most preferred product in $T$. Our subsequent analysis  works for any $ k \in \{1,\ldots,n \} $. Values of  $k$ close to the median, i.e., $\lceil \frac{n}{2} \rceil$, would represent customers who believe that the seller might offer them something reasonable in the middle of their preference list, while larger values of $k$ capture customers’ optimistic or gambling behavior.}
    

    We now derive the purchase probabilities and expected revenue under $\Opq_{(k)}(T, S,\bm r_S,\rho)$. \cor{Let $r^T_{(k)}$ denote the $k$-th smallest price among the traditional products in $T$. We also use the convention that $r^T_{(k)} = \max_{i \in T}r_i $ for $k > |T|$. } It is important to note that, without loss of generality, we can assume $\ropq \leq r^T_{(k)}$. This assumption is justified since the valuation of the opaque product is defined as the $k$-th order statistic among all traditional products in $T$, implying that the opaque product maximizes a customer's utility only if its price is lower than that of the $k$-th cheapest traditional price in $T$. Under the condition $\ropq \leq r^T_{(k)}$,  we show that the purchase probabilities, $  \pi \rpr{ \cdot \mid \Opq_{(k)}(T, S, \bm r_S,\rho)}$, admit closed-form solutions under the $k$-risk-averse customer assumption by applying the inclusion-exclusion principle. For convenience, we denote $\bm r_S^{I,\rho \wedge r }$ as the price vector where the prices of products in $I$ are replaced by $ \min \{ \rho, r_i \}$, i.e., $r^{I, \rho \wedge r }_{S,i} =\min \{ \rho, r_i \}$ for $ i \in I $ and $r^{I, \rho \wedge r }_{S,i} =r_i $ for $i \in S \setminus I $.

\cor{
\begin{lemma} \label{lem:choice-prob}
Fix $k \in \{ 1, 2,\ldots, n \}$. Define $\Low : = \{i \in T: \rho \leq r_i\}$. Under $\Opq_{(k)}(T, S, \bm r_S,\rho) :$ 
\begin{itemize}
 \item[(a)]  For $k \leq |T| $ and $\rho\leq r^T_{(k)} $, the purchase probabilities are given by:
 \begin{align*}
\pi \rpr{i \mid \Opq_{(k)}(T, S, \bm r_S,\rho)} =& \sum_{\substack{I\subseteq T\\ |I|\ge k}} (-1)^{|I|-k}\binom{|I|-1}{k-1} \pi \left(i \middle| \Trad \left(S, \bm{r}_S^{I, \rho \wedge r }\right)\right),  \quad i \in S \cup \{0\} \setminus L,   \\
 \pi \rpr{i \mid \Opq_{(k)}(T, S, \bm r_S,\rho)} =& \sum_{\substack{I\subseteq T \\ |I|\ge k}} (-1)^{|I|-k}\binom{|I|-1}{k-1} \mathbf{1}\{i \notin I\} \pi \left(i \middle| \Trad \left(S, \bm{r}_S^{I, \rho \wedge r }\right)\right),   \quad i \in  L.  \\ 
\pi \rpr{q \mid \Opq_{(k)}(T, S,\bm r_S,\rho)} =& \sum_{\substack{I\subseteq T \\ |I|\geq k}} (-1)^{|I|-k}\binom{|I|-1}{k-1} \sum_{j\in I \cap \Low }\pi \left(j \middle| \Trad \left(S, \bm{r}_S^{I,\rho \wedge r}\right)\right). 
 \end{align*}
\item[(b)] For $k \leq |T| $ and $\rho\leq r^T_{(k)} $, the expected revenue is
 \begin{align}\label{eq:rev-k-order}
     \Rev \left( \Opq_{(k)}(T, S,\bm r_S,\rho) \right)  =& \sum_{\substack{I\subseteq T\\ |I|\ge k}}  (-1)^{|I|-k}\binom{|I|-1}{k-1} \Rev \left( \Trad\left(S, \bm r_S^{I,\rho \wedge r} \right) \right) . 
 \end{align}
 For $\rho >  r^T_{(k)}  $, the opaque product will never be selected, and it naturally follows that 
      \begin{align}\label{eq: rho>rmin}
    \Rev \left( \Opq_{(k)}(T, S,\bm r_S,\rho) \right)  = \rev{\trad{S}{\bm r_S}}.
    \end{align}
 For $ k \geq |T| $ and $\rho \geq 0$, the opaque revenue collapses to a single MNL revenue:
 \begin{align}\label{eq:rev-k>|S|}
     \Rev \left( \Opq_{(k)}(T, S,\bm r_S,\rho) \right)  =& \Rev \left( \Trad\left(S, \bm r_S^{T,\rho \wedge r} \right) \right) . 
 \end{align}
\end{itemize}
\end{lemma}
}

\cor{
\textit{Proof.  } (a) For $k\leq |T|$ and $\rho \leq r^T_{(k)} $, the choice probability of choosing traditional product $i \in S$ is
\begin{align*}
&  \pi \rpr{i \mid \Opq_{(k)}(T, S, \bm r_S,\rho)} \\
 =& \mathbb{P}  \left[ U_i > U_j, \forall j \in S \cup \{0\} \setminus \{i\} \;\&\; U_i > V^T_{(k)} - \rho \right] \\
=& \sum_{\substack{I\subseteq T\\ |I|\ge k}}
(-1)^{|I|-k}\binom{|I|-1}{k-1} \mathbb{P}  \left[   U_i > U_j, \forall j \in S \cup \{0\} \setminus   \{i\} \;\&\;  \bigcap_{j \in I}\{U_i > V_j - \rho\} \right] \\
 =& \sum_{\substack{I\subseteq T\setminus\{i\}\\ |I|\ge k}} (-1)^{|I|-k}\binom{|I|-1}{k-1} \mathbb{P}  \left[   U_i > U_j, \forall j \in S \cup \{0\} \setminus (I \cup \{i\} ) \;\&\; \bigcap_{j \in I} \{ U_i > V_j - \min\{ \rho, r_j\} \}  \right] \\
&+ \sum_{\substack{I\subseteq T, i \in I\\ |I|\geq k}}
(-1)^{|I|-k}\binom{|I|-1}{k-1} \mathbb{P}  \left[   U_i > U_j, \forall j \in S \cup \{0\} \setminus  I \;\& \bigcap_{j \in I \setminus \{i\} } \{ U_i > V_j - \min\{ \rho, r_j\} \} \;\&\; U_i > V_i - \rho  \right]. 
\end{align*}
	The first equality follows from the definition of $V^T_q$ for the $ k-$th order statistic. The second equality holds by the general inclusion-exclusion principle: $\ones \{U_i>V^T_{(k)} -\rho \}= \sum_{\substack{I\subseteq T, |I|\ge k}} (-1)^{|I|-k}\binom{|I|-1}{k-1}\ones \{ U_i>V_j-\rho,\ \forall j\in I \}.$ The third equality follows from dividing the summation into two parts according to whether $i \in I$. We next refine the above probability based on whether $i \in S \setminus T$, $i \in \Low$, or $i \in T \setminus \Low  $.  For $i \in S \setminus T  $, the second summation is simply 0 since $i \notin T$. Because $ T \setminus\{i\}  = T $ and by the definitions of choice probabilities \eqref{eq:mnl} and $\bm r_S^{I, \rho \wedge r}$, we thus have that 
    \begin{align}\label{eq:opq-prob-i notin T} 
   \pi \rpr{i \mid \Opq_{(k)}(T, S, \bm r_S,\rho)} =\sum_{\substack{I\subseteq T \\ |I|\ge k}} (-1)^{|I|-k}\binom{|I|-1}{k-1} 
    \pi \left(i \middle| \Trad \left(S, \bm{r}_S^{I, \rho \wedge r }\right)\right) , \quad i \in S \setminus T .  
\end{align}
For $i \in \Low = \{j\in T: \rho \leq r_j \} $, since $V_i-\rho\ge U_i = V_i - r_i$, the second term is also 0 and thus
\begin{align}
   \pi \rpr{i \mid \Opq_{(k)}(T, S, \bm r_S,\rho)} =\sum_{\substack{I\subseteq T \setminus\{i\}\\ |I|\ge k}} (-1)^{|I|-k}\binom{|I|-1}{k-1} 
\pi \left(i \middle| \Trad \left(S, \bm{r}_S^{I, \rho \wedge r }\right)\right) , \quad i \in \Low .  \label{eq:opq-prob-i} 
\end{align}
For $i \in T \setminus  \Low$, since $U_i = V_i - r_i > V_i - \rho $, we have that
\begin{align}
& \pi \rpr{i \mid \Opq_{(k)}(T, S, \bm r_S,\rho)}  \notag \\
=& \sum_{\substack{I\subseteq T\setminus\{i\}\\ |I|\ge k}} (-1)^{|I|-k}\binom{|I|-1}{k-1} 
\pi \left(i \middle| \Trad \left(S,  \bm{r}_S^{I, \rho \wedge r }\right)\right) +  \sum_{\substack{I\subseteq T, i \in I\\ |I|\geq k}} (-1)^{|I|-k}\binom{|I|-1}{k-1} \pi \left(i  \middle|  \Trad \left(S,  \bm{r}_S^{I \setminus \{i\} , \rho \wedge r }\right)\right) \notag \\
=& \sum_{\substack{I\subseteq T\setminus\{i\}\\ |I|\ge k}} (-1)^{|I|-k}\binom{|I|-1}{k-1} 
\pi \left(i  \middle|  \Trad \left(S,  \bm{r}_S^{I, \rho \wedge r }\right)\right) +  \sum_{\substack{I\subseteq T, i \in I\\ |I|\geq k}} (-1)^{|I|-k}\binom{|I|-1}{k-1} \pi \left(i  \middle|  \Trad \left(S, \bm{r}_S^{I , \rho \wedge r }\right)\right) \notag \\
=&  \sum_{\substack{I\subseteq T\\ |I|\ge k}} (-1)^{|I|-k}\binom{|I|-1}{k-1} 
\pi \left(i \middle| \Trad \left(S,  \bm{r}_S^{I, \rho \wedge r }\right)\right) , \quad i \in T \setminus \Low.  \label{eq:opq-prob-i (2)} 
\end{align}
Again, the first equality holds by the definitions of choice probabilities \eqref{eq:mnl} and $\bm r_S^{I, \rho \wedge r}$. The second equality follows from $\pi \left(i | \Trad \left(S, \bm{r}_S^{I \setminus \{i\} , \rho \wedge r }\right)\right)  = \pi \left(i |  \Trad \left(S, \bm{r}_S^{I , \rho \wedge r }\right)\right)  $ for $i \in T \setminus \Low$ since $ \min\{r_i, \rho\} = r_i $. Similarly, the no-purchase probability is computed as
\begin{align}
    \pi \rpr{0 \mid \Opq_{(k)}(T, S, \bm r_S,\rho)} =& \sum_{\substack{I\subseteq T\\ |I|\ge k}} (-1)^{|I|-k}\binom{|I|-1}{k-1} \pi \left(0  \middle|  \Trad \left(S,  \bm{r}_S^{I, \rho \wedge r }\right)\right). \label{eq:opq-prob-0}
\end{align}
For the opaque product, the purchase probability is given by
\begin{align*}
\pi \rpr{q \mid \Opq_{(k)}(T, S, \bm r_S,\rho)}     =& 1 -  \pi \rpr{0 \mid \Opq_{(k)}(T, S, \bm r_S,\rho)}  -  \sum_{i\in S} \pi \rpr{i \mid \Opq_{(k)}(T, S, \bm r_S,\rho)}  \\
=& 1 -  \sum_{\substack{I\subseteq T\\ |I|\ge k}} (-1)^{|I|-k}\binom{|I|-1}{k-1} \left[ \sum_{i \in \{0\} \cup (S \setminus L) \cup (L \setminus I) } \pi \left(i  \middle|  \Trad \left(S,  \bm{r}_S^{I, \rho \wedge r }\right)\right) \right] \\
=& \sum_{\substack{I\subseteq T\\ |I|\ge k}} 
(-1)^{|I|-k}\binom{|I|-1}{k-1} \sum_{i \in I \cap L } \pi \left(i  \middle|  \Trad \left(S,  \bm{r}_S^{I, \rho \wedge r }\right)\right),
\end{align*}
where the second equality follows from \eqref{eq:opq-prob-i notin T}, \eqref{eq:opq-prob-i}, \eqref{eq:opq-prob-i (2)}, \eqref{eq:opq-prob-0}, and the observation that $ i \in L \cap I $ is excluded from the summation. The third equality follows from $ \sum_{\substack{I\subseteq T, |I|\ge k}} 
(-1)^{|I|-k}\binom{|I|-1}{k-1} = 1$. }

\cor{ (b) Fix $k \leq |T| $ and $\rho \leq r^T_{(k)}$. Given the purchase probabilities in (a), the expected revenue is
\begin{align*}
    &  \Rev \left( \Opq_{(k)}(T, S,\bm r_S,\rho) \right)  \\
    =& \sum_{i \in S} r_i \cdot \pi \rpr{i \mid \Opq_{(k)}(T, S, \bm r_S,\rho)}  + \rho \cdot \pi \rpr{q \mid \Opq_{(k)}(T, S, \bm r_S,\rho)}  \\
   =& \sum_{\substack{I\subseteq T\\ |I|\ge k}}  (-1)^{|I|-k}\binom{|I|-1}{k-1}  
   \left[ \sum_{i \in S \setminus(I \cap L)    } r_i  \pi \left(i \middle| \Trad \left(S, \bm{r}_S^{I, \rho \wedge r }\right)\right)  + \rho \sum_{i \in I \cap L } \pi \left(i \middle| \Trad \left(S, \bm{r}_S^{I, \rho \wedge r }\right)\right)  \right] \\
    =& \sum_{\substack{I\subseteq T\\ |I|\ge k}} (-1)^{|I|-k}\binom{|I|-1}{k-1} 
    \Rev \left( \Trad\left(S, \bm{r}_S^{I \cap L, \rho} \right) \right) 
    = \sum_{\substack{I\subseteq T\\ |I|\ge k}} 
(-1)^{|I|-k}\binom{|I|-1}{k-1} \Rev \left( \Trad\left(S, \bm{r}_S^{I,\rho \wedge r} \right) \right) .
\end{align*} }

\cor{
 When  $\rho>r^T_{(k)}$, let $ \mathcal{H} :=\{j\in T: r_j\le r^T_{(k)}\}$. Then $|\mathcal{H}|\geq k$ and for each $j \in \mathcal{H}$ we have $r_j\le r^T_{(k)}<\rho$,
hence $ U_j = V_j - r_j \ge V_j - r^T_{(k)}$. For the opaque utility, it holds that $U_q = V^T_{(k)} - \rho < V^T_{(k)} - r^T_{(k)} $. Now, among the $|\mathcal{H}|\geq k$ valuations $\{V_j:j\in \mathcal{H}\}$, the maximum is at least $V^T_{(k)}$, because in any subset of size at least $k$, the largest element is at least the $k$-th order statistic of the full set, i.e., $\max_{j\in \mathcal{H}} V_j \geq V^T_{(k)} $. Combining the above results, we have that 
\begin{align*}
    \max_{j\in \mathcal{H}} U_j \geq  \max_{j\in \mathcal{H}} \left(V_j - r^T_{(k)} \right) = \left(\max_{j\in \mathcal{H}} V_j\right) - r^T_{(k)} \geq V^T_{(k)} - r^T_{(k)} > U_q.
\end{align*} 
This implies that the utility of the opaque product will never be the highest one. Therefore, the selling model $  \Opq_{(k)}(T, S,\bm r_S,\rho)$ reduces to the traditional model $\trad{S}{\bm r_S}$, implying that $  \Rev \left( \Opq_{(k)}(T, S,\bm r_S,\rho) \right) = \rev{\trad{S}{\bm r_S}}.  $ }

\cor{When $k \geq |T|$, since $  V_{(k)}^T := \max_{i\in T} V_i = V_{(|T|)}^T$, the opaque revenue follows
\begin{align*}
  \Rev \left( \Opq_{(k)}(T, S,\bm r_S,\rho) \right)  =   \Rev \left( \Opq_{(|T|)}(T, S,\bm r_S,\rho) \right)  =  \Rev \left( \Trad\left(S, \bm{r}_S^{T,\rho \wedge r} \right) \right) .
\end{align*} 
Note that this holds for any $\rho \geq 0$ as $\bm{r}_S^{T,\rho \wedge r} = \bm{r}_S$ when $\rho > r_{(k)}^T = r_{(|T|)}^T  $. \hfill \Halmos }

We emphasize that the derivation of the choice probabilities for the $k$-risk-averse customers is novel and does not appear in \cite{mike2021power}, even for the 1-risk-averse case considered therein.  For $\rho \leq r^T_{(k)}$, Eq. \eqref{eq:rev-k-order} shows that the opaque revenue is an inclusion-exclusion sum of the revenues earned from the traditional MNL model where a subset of products are offered at the minimum of opaque price and the corresponding traditional price. While this closed-form expression provides valuable insights, its exponential number of terms poses significant challenges for solving price and assortment problems.
  


Building on these results, we next address these challenges in the context of price and assortment optimization problems under the opaque selling model. In Sections \ref{sec:price-support} and \ref{sec:joint asst support}, we study the flexible setting where the seller may choose any opaque support within the offered assortment. In Section \ref{sec:asst = support}, we then turn to the case where the opaque support is required to equal the offered assortment.

\newcommand{\opqast}[2]{\textsc{Opq}(#1, #2)}
\newcommand{\opqsast}[1]{\textsc{Opq}^*\rpr{#1}}
\newcommand{\opqsastequal}[2]{\textsc{Opq}^*\rpr{#1,#2} }
\newcommand{\tradast}[1]{\textsc{Trad}\rpr{#1}}

\section{Price Optimization} \label{sec:price-support}

We first consider price optimization under the Opaque-MNL model. Fix an assortment $S \subseteq \cal N$, \cor{an opaque support $T \subseteq S $, and $ k \in \{1,2,\ldots, n \} $.  We are interested in the following price optimization problem:
\begin{align}\label{opt: opq-price}
    \max_{{\bm r_S} \geq \bm{0},  \rho \geq 0  }  \Rev \left( \Opq_{(k)}(T, S,\bm r_S,\rho) \right)  . \tag{\sf Opq-Pricing}
\end{align}
}

Compared to the price optimization problem in the traditional MNL model, the primary challenge of solving \eqref{opt: opq-price} stems from the complex expression of $ \Rev \left( \Opq_{(k)}(T, S,\bm r_S,\rho) \right)$, which encompasses $2^{|T|} - \sum_{m=0}^{k-1} \binom{|T|}{m}  $ traditional revenue terms, as shown in Eq. \eqref{eq:rev-k-order}. Since the opaque revenue contains exponentially many terms, it is challenging to solve the first-order conditions directly. Recall that for the traditional MNL model, it is well established that offering every product at the same price is optimal (see, e.g., \citealt{hopp2005product}). Surprisingly, we discover that uniform pricing remains optimal for the Opaque-MNL model. Specifically, for any assortment $S \subseteq \cal N$, let $r_S^*$ represent the optimal uniform price under the traditional MNL model. Theorem \ref{thm-pricing} (proved in Section \ref{sec:pricing-proof}) shows that the optimal prices of \eqref{opt: opq-price} are $ \rho = r_i = r_S^* $, for all $ i \in S$. For $r \geq 0$, let $r \ones_S \in \mathbb{R} $ denote the price vector for the assortment $S$, where all prices are set to $r$.

\begin{theorem}\label{thm-pricing}
Uniform pricing is optimal for the price optimization problem \eqref{opt: opq-price}. 
For any assortment $S \subseteq \cal N$, \cor{opaque support $T \subseteq S $, $k \in \{ 1, \ldots, n \} $, } ${\bm r_S} \geq \bm{0}$, and $\rho \geq 0   $,   
\begin{equation}\label{eq:thm1-price}
  \Rev \left( \Opq_{(k)}(T, S,\bm r_S,\rho) \right)  \leq 
   \Rev \left( \Opq_{(k)}(T, S, r_S^*\ones_S,r_S^*) \right) = \rev{\tradast{S,r_S^*\ones_S}}.
\end{equation}
\end{theorem}

We point out that there are two potential ways to enhance revenue under the traditional MNL model $\trad{S}{\bm r_S}$: (i) adding an opaque product priced at the revenue-maximizing opaque price, and (ii) repricing every product to the optimal uniform price $r_S^*$. Theorem \ref{thm-pricing} establishes that repricing is a more powerful tool for revenue maximization than adding an opaque product. 
Note that \cite{mike2021power} show that the opaque selling strategy always (weakly) dominates the traditional selling model for 1-risk-averse customers ($k=1$) when valuation distributions are exchangeable.\footnote{Exchangeability generalizes the i.i.d. assumption: the joint distribution of customer valuations is invariant under any relabeling of the products.}  In contrast, Theorem \ref{thm-pricing} shows that the two strategies are equivalent under the MNL model, even if we violate the exchangeability assumption (meaning that valuations $v_i$ can be different). A possible explanation is as follows. In a segmented market where discriminatory pricing captures more revenue than single pricing, \cite{mike2021power} show that the opaque selling model can further enhance revenue by targeting more price-sensitive customers without lowering the prices of traditional products. However, in a market where single pricing is optimal (as in the MNL model), the market is not effectively price segmented.  Consequently, opaque selling models may not be advantageous in generating additional revenue, as the market lacks the necessary segmentation for such a strategy to be effective.

\cor{ 
Proving Theorem \ref{thm-pricing} is a considerable challenge due to the exponentially many terms involved in the revenue expression. The non-linearity introduced by the $\min \{\rho, r_i\}$ (denoted as $\rho \wedge r$) term in the opaque revenue expressions further complicates the analysis. Nevertheless, we resolve these difficulties by first simplifying the revenue expression (Lemma \ref{lemma: k - k-j}), showing that the inclusion-exclusion summation can be restricted solely to the subset of products $L = \{i \in T : r_i \ge \rho\} $. This observation effectively handles the complexity of the minimization term. Subsequently, using a change of variables, we decompose the opaque revenue into a fixed traditional revenue component and a remainder term governed by a concave auxiliary function. By exploiting this concavity, we establish that the revenue under the Opaque-MNL model is always bounded above by the revenue of the traditional MNL model at the optimal uniform price $r_S^*$ (Lemma \ref{lemma: main result}), thereby proving that uniform pricing remains optimal. In Section \ref{sec:pricing-proof}, we prove Theorem \ref{thm-pricing}.
}

\subsection{Proof of Theorem~\ref{thm-pricing} }\label{sec:pricing-proof}

Fix an assortment $S \subseteq \mathcal{N} $, an opaque support $T \subseteq S $, and $k \in \{1,\ldots,n \} $. If $k \geq |T| $, the opaque revenue reduces to a MNL revenue according to Eq. \eqref{eq:rev-k>|S|} in Lemma \ref{lem:choice-prob}, and thus Theorem \ref{thm-pricing} follows immediately. We also note that when $\rho > r^T_{(k)} $, \eqref{eq:thm1-price} holds naturally due to Eq. \eqref{eq: rho>rmin} and the optimality of $r_S^*.$ Hence, we restrict our attention to the case where $\rho \leq r^T_{(k)} $ and $ k \leq |T| -1 $ for the remainder of the proof. 

\cor{
Our first step is to simplify the opaque revenue expression in \eqref{eq:rev-k-order} by replacing the term $\min\{r_j, \rho\}$. Note that for any $j \in \Low = \{i \in T: r_i \geq \rho \}   $, $\min\{r_j, \rho \}  = \rho $ and for any $j \notin \Low$,  $ \min\{r_j, \rho \}  = r_j  $. Based on this observation, the price vector satisfies $\bm{r}_S^{I,\rho \wedge r} = \bm{r}_S^{I \cap \Low ,\rho }  $. This allows us to restrict the inclusion-exclusion summation in \eqref{eq:rev-k-order} to the subset $\Low$ rather than the entire set $T$, as established in Lemma \ref{lemma: k - k-j} below. Let $r^T_{(0)}:=-\infty$.
\begin{lemma}\label{lemma: k - k-j}
Fix $  T \subseteq S \subseteq \mathcal{N} $ and $k \in \{1,\ldots, |T|\}$. Suppose $ r^T_{(j)} < \rho \leq r^T_{(j+1)} $ for some $ j \in \{0, 1,\ldots, k-1\}$  and recall that $ \Low = \{i \in T: r_i \geq \rho \} $. The opaque revenue can be simplified as
\begin{align*}
     \Rev \left( \Opq_{(k)}(T, S,r_S,\rho) \right) =& \sum_{\substack{I\subseteq L\\ |I|\ge k -j}}  (-1)^{|I|-(k-j)}\binom{|I|-1}{k-j-1} \Rev \left( \Trad\left(S, r_S^{I  ,\rho } \right) \right) .
\end{align*}
\end{lemma} 
}

\cor{
Next, we present our main technical tool in Lemma \ref{lemma: main result}. 
This lemma establishes that the opaque revenue $   \Rev \left( \Opq_{(k)}(T, S,\bm{r}_S,\rho) \right)$  is upper bounded by the traditional MNL revenue at the optimal uniform price when the opaque price $\rho \leq \min_{i \in T} r_i$.
\begin{lemma}\label{lemma: main result}
    For any $T \subseteq S $, $k = 1,\ldots, |T|$, and $\rho \leq \min_{i \in T} r_i$, we have that
    \begin{align*}
        \Rev \left( \Opq_{(k)}(T , S,\bm{r}_S,\rho) \right) = \sum_{\substack{I\subseteq T \\ |I|\ge k}}  (-1)^{|I|-k}\binom{|I|-1}{k-1} \Rev \left( \Trad\left(S, \bm{r}_S^{I  ,\rho } \right) \right) \leq  \Rev \left( \Trad\left(S, r_S^* \mathbf{1}_S \right) \right). 
    \end{align*}
\end{lemma}
}

\cor{
We are now equipped to prove the main theorem. For any $\rho \in [0, r^T_{(k)}]$, there exists some index $j \in \{0,1,\ldots,k-1\}$ such that $r^T_{(j)} < \rho \leq r^T_{(j+1)}$. Consequently, the size of the subset $\Low$ is $|\Low| = |T| - j$. By combining the lemmas, we obtain that
    \begin{align*}
         \Rev \left( \Opq_{(k)}(T, S,\bm{r}_S,\rho) \right) =& \sum_{\substack{I\subseteq L\\ |I|\ge k -j}}  (-1)^{|I|-(k-j)}\binom{|I|-1}{k-j-1} \Rev \left( \Trad\left(S, \bm{r}_S^{I  ,\rho } \right) \right)  \\
         \leq & \Rev \left( \Trad\left(S, r_S^* \mathbf{1}_S \right) \right),
    \end{align*}
where the equality follows from Lemma \ref{lemma: k - k-j} and the inequality follows from applying  Lemma \ref{lemma: main result} with $ T = \Low $ and order statistic $k' = k -j  \leq |T| -j = |\Low| $. Finally, the equality in Theorem~\ref{thm-pricing} follows by setting $\bm r_S= r_S^* \mathbf{1}_S  $ and $\rho=r_S^*$ in Eq.~\eqref{eq:rev-k-order}, which gives $ \Rev \left( \Opq_{(k)}(T, S, r_S^* \mathbf{1}_S,r_S^*) \right)  = \sum_{\substack{I\subseteq T, |I|\ge k}}  (-1)^{|I|-k}\binom{|I|-1}{k-1} \Rev \left( \Trad\left(S, r_S^* \mathbf{1}_S \right) \right) = \Rev \left( \Trad\left(S, r_S^* \mathbf{1}_S \right) \right) $.  \hfill \Halmos
}

We prove Lemma \ref{lemma: main result} here and defer the proof of Lemma \ref{lemma: k - k-j} to Appendix \ref{app:pricing-general}.

\cor{
\underline{ \textit{Proof of Lemma \ref{lemma: main result}. }}
 We use change of variables to streamline the analysis. For $i \in T$, let 
 \begin{align*}
     w_i = \frac{e^{v_i - \rho}}{1 + \sum_{j \in T} e^{v_j-\rho} + \sum_{j \in S \setminus T} e^{v_j -r_j }} \in (0,1) \text{ and } \beta_i = e^{\rho - r_i} \in (0,1].
 \end{align*}
 We first decompose $\Rev \left( \Trad\left(S, \bm{r}_S^{I ,\rho } \right) \right)$ into $\Rev \left( \Trad\left(S, \bm{r}_S^{T  ,\rho } \right) \right)$ and a residual. Through algebraic manipulation, we observe that for any $I \subseteq T$: 
    \begin{align}
        \Rev \left( \Trad\left(S, \bm{r}_S^{I  ,\rho } \right) \right)  =& \frac{\sum_{i \in I}\rho e^{v_i-\rho} + \sum_{i \in S \setminus I} r_i e^{v_i -r_i }}{ 1 + \sum_{i \in I} e^{v_i-\rho} + \sum_{i \in S \setminus I} e^{v_i -r_i }   } \notag \\
      =& \frac{ \frac{\sum_{i \in T}\rho e^{v_i-\rho} + \sum_{i \in S \setminus T} r_i e^{v_i -r_i }}{ 1 + \sum_{i \in T} e^{v_i-\rho} + \sum_{i \in S \setminus T} e^{v_i -r_i }  } + \frac{\sum_{i \in T \setminus I} \left( r_i e^{v_i -r_i} -\rho e^{v_i - \rho} \right)  }{ 1 + \sum_{i \in T} e^{v_i-\rho} + \sum_{i \in S \setminus T} e^{v_i -r_i } } }{ 1 + \frac{\sum_{i \in T \setminus I} e^{v_i -r_i} - e^{v_i -\rho} }{1 + \sum_{i \in T} e^{v_i-\rho} + \sum_{i \in S \setminus T} e^{v_i -r_i }}   } \notag \\   
     =& \frac{  \Rev \left( \Trad\left(S, \bm{r}_S^{T  ,\rho } \right) \right) + \frac{\sum_{i \in T \setminus I} \left( r_i e^{v_i -r_i} -\rho e^{v_i - \rho} \right) }{ 1 + \sum_{i \in T} e^{v_i-\rho} + \sum_{i \in S \setminus T} e^{v_i -r_i } }}{ 1 + \frac{\sum_{i \in T \setminus I} e^{v_i -r_i} - e^{v_i -\rho} }{1 + \sum_{i \in T} e^{v_i-\rho} + \sum_{i \in S \setminus T} e^{v_i -r_i }}   } \notag \\
     =&   \frac{  \Rev \left( \Trad\left(S, \bm{r}_S^{T  ,\rho } \right) \right) + \sum_{i \in T \setminus I} w_i \left( \left(\rho -\ln \beta_i \right)\beta_i - \rho \right) }{ 1 + \sum_{i \in T \setminus I} w_i (\beta_i  - 1)  } \notag \\
     =& \Rev \left( \Trad\left(S, \bm{r}_S^{T  ,\rho } \right) \right) + \sum_{i \in T \setminus I} w_i \frac{-\ln \beta_i \cdot \beta_i + \left[\rho -\Rev \left( \Trad\left(S, \bm{r}_S^{T  ,\rho } \right) \right)   \right] (\beta_i -1) }{1 - \sum_{j \in T \setminus I} w_j(1-\beta_j) }  \notag \\
     =& \Rev \left( \Trad\left(S, \bm{r}_S^{T  ,\rho } \right) \right) + \sum_{i \in T \setminus I}\frac{w_i f(\beta_i, \rho, \bm{r}_{S \setminus T}) }{1 - \sum_{j \in T \setminus I} w_j(1-\beta_j) }, \label{eq:decomposition}
    \end{align}
    where we define the auxiliary function
    \begin{align*}
        f(\beta_i, \rho, \bm{r}_{S \setminus T}) := -\ln \beta_i \cdot \beta_i + \left[\rho -\Rev \left( \Trad\left(S, \bm{r}_S^{T  ,\rho } \right) \right)   \right] \cdot (\beta_i -1).
    \end{align*}
    Crucially, $f(\beta_i, \rho, \bm{r}_{S \setminus T})$ is concave in $\beta_i$ as $f''(\beta) = -1/\beta <0$. Then, the opaque revenue can be re-expressed as:
    \begin{align}
          & \Rev \left( \Opq_{(k)}(T, S,\bm{r}_S,\rho) \right) \notag\\
          =& \sum_{\substack{I\subseteq T \\ |I|\ge k}}  (-1)^{|I|-k}\binom{|I|-1}{k-1} \left[ \Rev \left( \Trad\left(S, \bm{r}_S^{T  ,\rho } \right) \right) +\frac{  \sum_{i \in T \setminus I} w_i  f(\beta_i, \rho, \bm{r}_{S \setminus T}) }{1 - \sum_{j \in T \setminus I} w_j(1-\beta_j) }\right]   \notag \\
        =& \Rev \left( \Trad\left(S, \bm{r}_S^{T  ,\rho } \right) \right)  + \sum_{i \in T}w_i  f(\beta_i, \rho, \bm{r}_{S \setminus T})  \sum_{\substack{I\subseteq T \setminus \{i\}  \\ |I|\ge k}}\frac{  (-1)^{|I|-k}\binom{|I|-1}{k-1} }{1 - \sum_{j \in T \setminus I} w_j(1-\beta_j)} . \label{eq:opq-rev-simplifies}
    \end{align}
    Plugging \eqref{eq:decomposition} into opaque revenue \eqref{eq:rev-k-order} gives the first equality. The second equality follows from the identity $\sum_{\substack{I\subseteq T , |I|\ge k}}  (-1)^{|I|-k}\binom{|I|-1}{k-1} = 1 $ and rearranging terms. 
    }

  \cor{   For any $r \geq \rho$, let $\bar{w} := \sum_{i \in T} w_i \in (0,1)$ and $\beta := e^{\rho - r}$. A similar computation shows that
\begin{align*}
     \Rev \left( \Trad\left(S, \bm{r}_S^{T  ,r } \right) \right) = \Rev \left( \Trad\left(S, \bm{r}_S^{T  ,\rho } \right) \right) + \frac{\bar{w} f(\beta, \rho, \bm{r}_{S \setminus T }) }{1 - \bar{w} (1-\beta)}.
\end{align*}
Furthermore, because $r$ only affects $\beta$, we have that
\begin{align}\label{eq: max-r-beta}
    \max_{r \geq \rho }  \;  \Rev \left( \Trad\left(S, \bm{r}_S^{T  ,r } \right) \right) = \Rev \left( \Trad\left(S, \bm{r}_S^{T  ,\rho } \right) \right)  + \sup_{\beta \in (0,1]} \frac{\bar{w}  f(\beta, \rho, \bm{r}_{S \setminus T}) }{1 - \bar{w} (1-\beta)}.
\end{align}
 Before proceeding with the proof, we first show that the coefficient of $w_i f(\beta_i, \rho, \bm{r}_{S \setminus T }) $ in \eqref{eq:opq-rev-simplifies} is non-negative and provide an upper bound.  Note that $\sum_{i \in T} w_i (1- \beta_i) = \frac{\sum_{i \in T } (e^{v_i -\rho} - e^{v_i -r_i} ) }{1 + \sum_{i \in T} e^{v_i-\rho} + \sum_{i \in S \setminus T} e^{v_i -r_i }}   < 1 $. The proof of Claim \ref{claim: Binomial} is provided in Appendix \ref{app:pricing-general}.
\begin{claim}\label{claim: Binomial}
      For any $T \subseteq S$, $k \leq |T| -1 $, $w_i \in  (0,1)$, $ \beta_i \in (0,1] $, and $\sum_{i \in T} w_i(1-\beta_i) <1 $, it holds that
    \begin{align*}
      0 \leq  \sum_{\substack{I\subseteq T \setminus \{i\}  \\ |I|\ge k}}\frac{  (-1)^{|I|-k}\binom{|I|-1}{k-1} }{1 - \sum_{j \in T \setminus I} w_j(1-\beta_j)} \leq \frac{1}{1 - \sum_{j \in T} w_j (1- \beta_j) }, \quad i \in T.
    \end{align*}
\end{claim}
}

\cor{ Next, our goal is to prove that \eqref{eq:opq-rev-simplifies} is less than or equal to \eqref{eq: max-r-beta}.
Let $T^+ := \{ i \in T, f(\beta_i) \geq 0\}$. If $T^+ = \emptyset$, since the coefficient is non-negative shown in Claim \ref{claim: Binomial} above, we immediately obtain that
\begin{align*}
      \Rev \left( \Opq_{(k)}(T , S,\bm{r}_S,\rho) \right) \leq \Rev \left( \Trad\left(S, \bm{r}_S^{T  ,\rho } \right) \right) \leq  \Rev \left( \Trad\left(S, r_S^* \mathbf{1}_S \right) \right).
\end{align*}
For $T^+  \neq \emptyset $, we denote $\hat{w} = \sum_{i\in T^+} w_i $ and $\hat{\beta} = \frac{\sum_{i \in T^+} w_i \beta_i }{\sum_{i \in T^+} w_i} \in (0,1] $. We bound the second term in \eqref{eq:opq-rev-simplifies} as follows:
\begin{align*}
&\sum_{i \in T} w_i  f(\beta_i, \rho, \bm{r}_{S \setminus T}) \cdot \sum_{\substack{I\subseteq T \setminus \{i\}  \\ |I|\ge k}}\frac{  (-1)^{|I|-k}\binom{|I|-1}{k-1} }{1 - \sum_{j \in T \setminus I} w_j(1-\beta_j)} \\
 \leq &  \frac{  \sum_{i \in T^+}   w_i  f(\beta_i, \rho, \bm{r}_{S \setminus T})}{1 - \sum_{j \in T} w_j (1- \beta_j) } 
    \leq  \frac{\hat{w}  f(\hat{\beta}, \rho, \bm{r}_{S \setminus T}) }{1 - \sum_{j \in T} w_j(1-\beta_j) } 
    \leq \frac{\hat{w}  f(\hat{\beta}, \rho, \bm{r}_{S \setminus T}) }{1 -  \bar{w} + \hat{w}\hat{\beta} } 
  \leq \frac{\bar{w}  f(\hat{\beta}, \rho, \bm{r}_{S \setminus T}) }{1 -  \bar{w} + \bar{w} \hat{\beta}}.
\end{align*}
Here, the first inequality follows from Claim~\ref{claim: Binomial} and dropping out negative terms. The second inequality uses Jensen's inequality so that $\frac{ \sum_{i \in T^+} w_i  f(\beta_i, \rho, \bm{r}_{S \setminus T}) }{\sum_{i\in T^+} w_i} \leq   f(\frac{\sum_{i \in T^+} w_i \beta_i }{\sum_{i \in T^+} w_i}, \rho, \bm{r}_{S\setminus T}) =  f(\hat{\beta}, \rho, \bm{r}_{S \setminus T})$. The third inequality follows from $\bar{w} = \sum_{i \in T }w_i $ and $\sum_{i \in T } w_i \beta_i \geq \sum_{i \in T^+ } w_i \beta_i := \hat{w}  \hat{\beta} $. Finally, the last inequality is true due to the monotonicity in $\hat{w}$ (as $  f(\hat{\beta}, \rho, \bm{r}_{S \setminus T}) \geq 0 $) and $\hat{w} \leq \bar{w}$. Combining these results, we conclude that
\begin{align*}
  \Rev \left( \Opq_{(k)}(T, S,\bm{r}_S,\rho) \right) \leq& \Rev \left( \Trad\left(S, \bm{r}_S^{T  ,\rho } \right) \right) + \frac{\bar{w}  f(\hat{\beta}, \rho, \bm{r}_{S \setminus T}) }{1 -  \bar{w} + \bar{w} \hat{\beta}}  \\
  \leq &\Rev \left( \Trad\left(S, \bm{r}_S^{T  ,\rho } \right) \right) + \sup_{\beta \in (0,1] }\frac{\bar{w}  f(\beta, \rho, \bm{r}_{S \setminus T}) }{1 - \bar{w} (1-\beta)} \\
     =& \max_{r \geq \rho} \;  \Rev \left( \Trad\left(S, \bm{r}_S^{T  ,r } \right) \right) \\
    \leq &  \Rev \left( \Trad\left(S, r_S^* \mathbf{1}_S \right) \right). 
\end{align*}
\hfill \Halmos
}

\section{Assortment Optimization}\label{sec:joint asst support}

In this section, we analyze the assortment optimization problem under the Opaque-MNL model, which involves determining the revenue-maximizing subset of $\cal N$ and its optimal opaque support offered to customers. For instance, in the Amazon example illustrated in Figure~\ref{fig:examples}\subref{fig:amazon}, the seller may choose to offer traditional products in black, blue, and red along with the ``Color Will Vary'' option to potentially increase revenue. 

\cor{
In the assortment optimization problem, we treat traditional prices as fixed parameters, while allowing the seller to choose the offered assortment, the opaque support, and the opaque price. This is motivated by the operational realities of online travel agencies and third-party marketplaces. Figure \ref{fig:hotel} provides a concrete example of an opaque hotel offering on Hotwire. Typically, these platforms lack the power to change the prices of branded products (such as specific hotels or flights), which are often set by service providers or locked by rate parity agreements.  In contrast, the platform has full control over the opaque product because it is a special offer created by the platform itself.
This distinction allows the agency to freely determine the opaque price and decide which products to include in the assortment and opaque support. Therefore, even though traditional prices are fixed, the platform can effectively segment the market and maximize total revenue by jointly optimizing the opaque offer and its structure.
\begin{figure}[h!]
\FIGURE{\centering
 \includegraphics[height=0.3\textwidth]{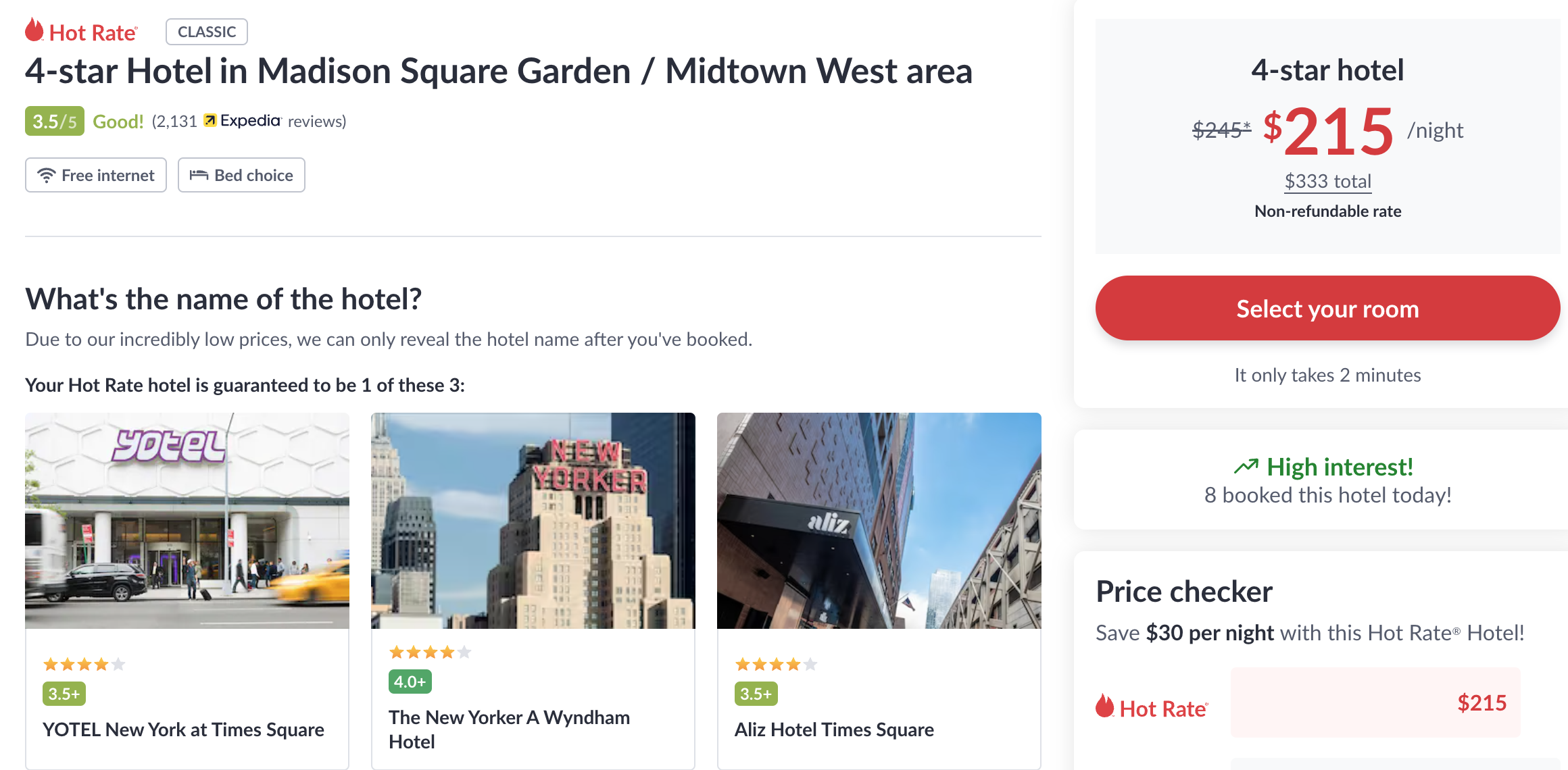}     \label{fig:hotel}
}{An example of an opaque product -- ``Hot Rate Hotel" on Hotwire.}
{A search of hotels in New York on Hotwire reveals an opaque hotel with ``Hot Rate'' at \$215 among three options, including YOTEL New York at Times Square (priced at \$246), The New Yorker A Wyndham (\$299), and Aliz Hotel Times Square (\$375). } 
\end{figure}
}

With the price vector $\bm r_{\cal N}$ of traditional products fixed, for any assortment $S \subseteq \mathcal{N} $, \cor{opaque support $ T \subseteq S $, and $k \in \{1,2,\ldots,n\} $, we define the selling model in which the opaque price is endogenously optimized as $\Opq^*_{(k)}(T, S, \bm{r}_S)$, where the optimal opaque price is denoted by $\rho^*(T)$ when the assortment $S$ is clear from context. Recall that when $k > |T| $, we assume that customers value the opaque product as the maximum of valuations of traditional products, i.e., $V_q^T =  \max_{j \in T} V_j  $.  We are interested in the following optimization problem:
\begin{align}\label{asst-opt: opaque support}
     \max_{S \subseteq \mathcal{N}} \max_{ T \subseteq S } \max_{\rho \geq 0 }   \rev{\Opq_{(k)}( T, S, \bm{r}_S,\rho )  }.  \tag{\sf Opq-Assortment-Support}
\end{align}
}

The unique feature of \eqref{asst-opt: opaque support} is that the opaque price is optimized for given assortment and opaque support, and we highlight the necessity as follows: For fixed $k \in \{1,\ldots,n \} $, recall that the valuation of the opaque product depends on the offered opaque support, i.e., $V_q^T = V_{(k)}^T $. Therefore, for two distinct assortments and opaque supports, it is not reasonable to set the same opaque price while the valuation distributions of the opaque product differ. Additionally, the constraint $\rho \leq r^T_{(k)} $ may render a feasible opaque price for one assortment (opaque support) infeasible for another. Furthermore, Proposition \ref{prop:condition} and Figure \ref{fig:unnested-averse} (when $T = S$) demonstrate that introducing an opaque product with an arbitrary feasible opaque price may reduce revenue compared to the traditional MNL model, thereby incentivizing the seller to set an optimal opaque price. In light of this analysis, we assume that the opaque price is endogenously optimized in the assortment optimization problem. 

\cor{
Despite the complexity of \eqref{asst-opt: opaque support}, we surprisingly show that the optimal assortment is nested by revenue, matching the structure under the traditional MNL model. The intuition is that opaque support choice lets the seller control which products shape the opaque option, excluding products that would otherwise weaken the opaque product or generate excessive cannibalization. This ability to isolate the opaque product yields behavior closer to standard MNL assortment optimization. We formalize this in the following theorem.
\begin{theorem}[Nested-by-revenue]\label{Thm: revenue-ordered} Fix $k \in \{1,2,\ldots,n\} $. Assume w.l.o.g. that revenues are ordered as $r_1 \geq r_2 \geq \cdots \geq r_n$. Then, there exists $S^*$ an optimal assortment of \eqref{asst-opt: opaque support}  such that $S^* = \{1,\ldots,j\}$ for some $j \in \{1,2,\ldots,n\}  $.
\end{theorem}
}

\cor{
A key insight in proving Theorem \ref{Thm: revenue-ordered} is that the seller never needs a large opaque support: for any assortment $S$, there exists an optimal opaque support $T^*\subseteq S$ with $|T^*|\le k$. Intuitively, once $|T|>k$, additional supported traditional products beyond $k$  contribute little to the customer’s $k$-th order-statistic guarantee but create downside by restricting the opaque price or inducing additional cannibalization through price truncation ($\min\{\rho, r_i \}$). This support reduction is pivotal because when $|T^*|\le k$, the opaque revenue collapses to a single traditional MNL revenue with a modified price vector (Eq. \eqref{eq:rev-k>|S|}), after which classical MNL monotonicity arguments imply a nested-by-revenue structure for the optimal assortment.  }

\cor{
Although the optimal assortment is nested-by-revenue, it remains to efficiently compute the optimal opaque support and opaque price for a given assortment, i.e., solving $\max_{\rho\ge 0}  \max_{\substack{T\subseteq S, |T|\le k} } \Rev \left(\Trad  \left(S,\bm r_S^{T,   \rho\wedge r}\right)\right)$, where the reduction to $|T|\le k$ follows from Lemma~\ref{lem: T* <=k} below and Eq.~\eqref{eq:rev-k>|S|}. For any $\epsilon > 0$, we develop a fully polynomial-time approximation scheme (FPTAS) for this problem, i.e., we compute a $(1-\epsilon)$-approximate solution in time polynomial in both the input size and $1/\epsilon$. Specifically, for any $\rho>0$, we reformulate the support problem $\max_{\substack{T\subseteq S, |T|\le k} } \Rev \left(\Trad  \left(S,\bm r_S^{T,   \rho\wedge r}\right)\right)$ as an MNL assortment optimization under totally unimodular constraints, which can be solved in polynomial time \citep{sumida2021revenue}. We then perform a geometric grid search on $\rho$ and return the opaque price with the largest revenue. We defer the details of the reformulation and geometric grid search to Appendix \ref{app:algorithm}.  
}

\subsection{Proof of Theorem \ref{Thm: revenue-ordered}.}
\cor{
Let $F(S)$ denote the optimal revenue of assortment $S$, i.e.,
 \begin{align*}
     F(S):= \max_{ T \subseteq S } \max_{\rho \geq 0 }   \rev{\Opq_{(k )}(T , S , \bm{r}_S,\rho )  }.
 \end{align*}
 For any assortment $S$, we show in Lemma \ref{lem: T* <=k} below that there exists an optimal opaque support $T^*$ with $|T^*| \leq k$. Then, according to Eq. \eqref{eq:rev-k>|S|}, the opaque revenue can be simplified as a traditional MNL revenue with price vector $\bm r_S^{T^*,\rho \wedge r}$, facilitating the subsequent analysis. The proof is provided in Appendix \ref{app:asst-general}.
\begin{lemma}\label{lem: T* <=k}
  For every offered assortment $S$, there exists an optimal opaque support $T^* \subseteq S $ with $|T^*| \leq k$ attaining $F(S)$.
\end{lemma}   }

\cor{
We next establish two monotonicity properties of $F(\cdot)$: adding a high-revenue ($ \geq F(S) $) product to $S$ or removing a low-revenue ($ \leq F(S) $) product from $S$ cannot decrease $F(S)$, i.e.,
 \begin{align}
     F(S\setminus\{i\}) \geq &  F(S), \quad i \in S, \ r_{i} \leq F(S),  \label{eq: remove->high} \\
     F(S\cup\{j\}) \geq&  F(S), \quad j \notin S, r_{j} \geq F(S). \label{eq: add->high}
 \end{align}
 To prove these results, we apply a well-known property of traditional MNL revenue.
\begin{claim}\label{claim:mnl-add-remove} 
Fix any assortment $S$ and price vector $\bm r_S$. (i) If $i\in S$ and $ r_i \leq \Rev(\Trad(S,\bm r_S))$, then $ \Rev(\Trad(S\setminus\{i\},\bm r_{S\setminus\{i\}})) \ge \Rev(\Trad(S,\bm r_S)).$ (ii) If $j\notin S$ and $r_j\ge \Rev(\Trad(S,\bm r_S))$, then $\Rev(\Trad(S\cup\{j\},\bm r_{S\cup\{j\}}))\ge \Rev(\Trad(S,\bm r_S))$.
\end{claim} 
}

\cor{
We first prove \eqref{eq: remove->high}. According to Lemma \ref{lem: T* <=k}, we can pick an optimal pair $(T^*,\rho^*)$ attaining $F(S)$ with $|T^*| \le k$. Then, Eq. \eqref{eq:rev-k>|S|} shows that the opaque revenue is simplified as:
\begin{align*}
F(S) = \Rev\left(\Opq_{(k)}(T^* , S,\bm r_S,\rho^*)\right)
= \Rev\left(\Trad(S,\bm r_S^{T^*,\rho^* \wedge r})\right).
\end{align*}
 Based on Claim \ref{claim:mnl-add-remove}(i) and the fact that $ (\bm r_S^{T^*,\rho^*\wedge r})_{i} \le r_{i} \le F(S) $, we can remove product $i$ from $S$ in the traditional MNL model without decreasing revenue:
\begin{align*}
\Rev\left(\Trad(S\setminus\{i\},\bm r_{S\setminus\{i\}}^{T^* \setminus \{i\} ,\rho^*\wedge r})\right) \ge \Rev\left(\Trad(S,\bm r_S^{T^*,\rho^*\wedge r})\right) = F(S) .
\end{align*}
Since $T^*\setminus\{i\}\subseteq S\setminus\{i\}$ and $|T^*\setminus\{i\}|\le k$, $(T^*\setminus\{i\},\rho^*)$ is feasible for $F(S\setminus\{i\})$. We conclude that
\begin{align*}
F(S\setminus\{i\}) &\ge \Rev\left(\Opq_{(k)}(T^*\setminus\{i\},S\setminus\{i\},\bm r_{S\setminus\{i\}},\rho^*)\right) = \Rev\left(\Trad(S\setminus\{i\},\bm r_{S\setminus\{i\}}^{T^* \setminus \{i\} ,\rho^*\wedge r})\right) \geq F(S),
\end{align*}
where the equality follows from Eq.~\eqref{eq:rev-k>|S|}.
}
\cor{
Next, we prove Eq. \eqref{eq: add->high}. Similarly, since $r_{j} \geq F(S) $, Claim \ref{claim:mnl-add-remove}(ii) implies that adding product $j$ to $S$ does not decrease revenue: 
\begin{align*}
\Rev\left(\Trad(S\cup\{j\}, \bm r_{S\cup\{j\}}^{T^*  ,\rho^*\wedge r})\right) \geq \Rev\left(\Trad(S,\bm r_S^{T^*,\rho^*\wedge r})\right) = F(S) .
\end{align*}
Moreover, $(T^*, \rho^*)$ is feasible for $F(S \cup \{j\})$ with $|T^*| \leq k $. Hence we obtain the desired result:
\begin{align*}
F(S \cup \{j\}) \geq \Rev\left(\Opq_{(k)}(T^* , S\cup\{j\},\bm r_{S\cup\{j\}},\rho^*)\right) = \Rev\left(\Trad(S\cup\{j\},\bm r_{S\cup\{j\}}^{T^* ,\rho^*\wedge r})\right) \geq F(S).
\end{align*} 
}

\cor{
We now complete the proof of Theorem \ref{Thm: revenue-ordered}. Let $S^*$ be an optimal assortment to problem \eqref{asst-opt: opaque support}. We want to show there exists an optimal assortment $\tilde{S}$ that is revenue-ordered.  To see this, starting from $S^*$, repeatedly remove any product $ i \in S^* $ with $r_i < F(S^*)$. Each such removal does not decrease revenue by \eqref{eq: remove->high} and preserves optimality since $F(S^*)$ is the global maximum. Next, repeatedly add any product $j \notin S^*$ with $r_j \geq F(S^*)$. By \eqref{eq: add->high}, each addition also preserves optimality. After finitely many steps, we get an optimal assortment $\tilde{S}$ that is revenue-ordered. \hfill \Halmos 
}

\section{Optimization with Full Opaque Support}\label{sec:asst = support}
In this section, we consider the case where the opaque support is full, i.e., the opaque support is required to coincide with the offered assortment so that $T=S$. Under this restriction, we omit the support from the notation and write $\Opq_{(k)}(S,\bm r_S,\rho)$ for the corresponding Opaque-MNL model. 

\cor{
Before studying the associated optimization problems, we first provide sufficient conditions under which the Opaque-MNL model dominates (or is dominated by) the traditional MNL model, holding the assortment $S$ and the traditional prices $\bm r_S$ fixed. Specifically, we establish that when all prices (both traditional and opaque) do not exceed the optimal uniform price $r_S^*$, the revenue under the Opaque-MNL model is bounded above by that of the traditional MNL model. Conversely, when traditional prices exceed $r_S^*$, the introduction of an opaque product can strictly enhance revenue, provided the opaque price is chosen within an appropriate interval. 
\begin{proposition}\label{prop:condition}
    Fix $k \in \{1,2,\ldots, n \}$ and an assortment $S \subseteq \mathcal{N} $. Let $\underline{\rho}_i$ denote the unique root in the interval $(0, r_i)$ of the equation $r_i e^{\rho - r_i} - \rho + (r_i - \rho) \cdot \sum_{j \in S} e^{v_j -r_i} =0 $ when it exists, for all $i \in S$. The following holds  
\begin{itemize}
 \item[(a)] If $ r_i \leq r_S^* $, $\forall i \in S$ and $\rho \leq \min_{i \in S} r_i $, then $  \Rev \left( \Opq_{(k)}( S,\bm{r}_S,\rho) \right)  \leq \Rev \left( \Trad\left(S, \bm{r}_S \right) \right)  $.
 \item[(b)] If $ r_i > r_S^* $, $\forall i \in S$ and $\rho \in (\max_{i \in S } \underline{\rho}_i, \min_{i \in S} r_i)  $, then  $  \Rev \left( \Opq_{(k)}( S, \bm{r}_S,\rho) \right)  > \Rev \left( \Trad\left(S, \bm{r}_S \right) \right)  $.
\end{itemize}
\end{proposition}  
}

\cor{The proof of Proposition \ref{prop:condition} is presented in Appendix \ref{append-same-r}. 
The sufficient conditions in Proposition \ref{prop:condition} require the opaque price to be no greater than the minimum traditional price in the opaque support. We emphasize that this requirement is mild and consistent with industry practice. In our motivating examples of notebooks and car rentals, the opaque product is consistently priced below the cheapest traditional products (e.g., \$10.99 versus \$11.24 and \$71 versus \$74). 
}

\subsection{Price Optimization}
Fix an assortment $S \subseteq \cal N$ and $ k \in \{1,2,\ldots,n \} $. The price optimization problem is
\begin{align}\label{opt: opq-price-2}
    \max_{{\bm r_S} \geq \bm{0},  \rho \geq 0  }   \Rev \left( \Opq_{(k)}(S,\bm r_S,\rho) \right)  .
\end{align}
Applying Theorem \ref{thm-pricing} with $T=S$, we immediately obtain the following corollary, which shows that uniform pricing remains optimal.
\begin{corollary}
  Uniform pricing is optimal for the price optimization problem \eqref{opt: opq-price-2}. 
For any assortment $S \subseteq \cal N$, $k \in \{ 1, \ldots, n \} $, ${\bm r_S} \geq \bm{0}$, and $\rho \geq 0   $,   
\begin{align*}
  \Rev \left( \Opq_{(k)}(S,\bm r_S,\rho) \right)  \leq  \Rev \left( \Opq_{(k)}(S,r_S^*\ones_S,r_S^*) \right) 
  = \rev{\tradast{S,r_S^*\ones_S}}.
\end{align*}
\end{corollary}
Thus, even in the special case where the opaque support must equal the offered assortment, allowing the seller to introduce an opaque product does not improve upon the optimal uniform-pricing solution of the traditional MNL model.

\newcommand{\set}{\tilde{S}}

\subsection{Assortment Optimization} \label{sec:assort}
 We now study the assortment optimization problem under the Opaque-MNL model when the opaque support is required to coincide with the offered assortment. Recall that, given fixed price vector $\bm r_{\cal N}$ of traditional products, for any assortment $S \subseteq \mathcal{N} $ and $k \in \{1,2,\ldots,n\} $, we use $\Opq^*_{(k)}(S,\bm r_S)$ to denote the selling model in which the opaque price is endogenously optimized, and we let $\rho^*(S)$ denote the corresponding optimal opaque price. When $k > |S|$, we adopt the convention that customers value the opaque product as the maximum valuation among the traditional products in $S$, that is, $V_q^S = \max_{j \in S} V_j$. We are interested in the following optimization problem:
\begin{align}\label{opt:assort}
       \max_{S \subseteq \mathcal{N}} \rev{\Opq^*_{(k)}(S, \bm{r}_S)} =  \max_{S \subseteq \mathcal{N}} \max_{\rho \geq 0 }   \rev{\Opq_{(k )}(S, \bm{r}_S,\rho )  } . \tag{\sf Opq-Assortment}
\end{align} 

Compared with \eqref{asst-opt: opaque support}, problem \eqref{opt:assort} loses the flexibility to choose the opaque support separately, and as a result the optimal assortment no longer admits a similarly clean structure in general. In particular, the opaque revenue does not reduce to a traditional MNL revenue and continues to involve the inclusion-exclusion summation. The above assortment optimization problem is therefore challenging because of the endogenously optimized opaque price, the exponential terms in the objective function, and its combinatorial nature, which we will explore further in Section \ref{sec: challenge}. Nevertheless, Section \ref{sec:assort-uniform} demonstrates that problem \eqref{opt:assort} is tractable when all traditional products are uniformly priced. In this case, Theorem \ref{thm:nested-val} shows that the optimal assortment has an elegant structure: the assortment is nested-by-valuation, meaning we choose all products with intrinsic valuation $v_i$ above a threshold. This scenario, where traditional products are uniformly priced, is also closely related to the sales maximization problem studied in \cite{berbeglia2020assortment} and \cite{housni2023placement}, which can be viewed as a special case of revenue maximization. For general prices, we propose a natural heuristic in Section \ref{sec:general-assort} that performs well in numerical experiments (Section \ref{sec:numerics}).

\subsubsection{Challenges of Assortment Problem \eqref{opt:assort}}\label{sec: challenge}
The main distinction between problem \eqref{opt:assort} and traditional assortment problems is the endogenously optimized opaque price and the structure of the opaque product's utility. Given these features, it is hard to compare the revenues of an assortment as the optimal opaque price is not explicitly given, thereby making the assortment optimization problem more challenging. Additionally, it turns out that the fundamental properties of the traditional MNL model, such as substitutability,  the inclusion of the highest revenue product in the optimal assortment, and the nested-by-revenue structure of the optimal assortment, do not extend to opaque selling.

{\bf Violation of substitutability.} 
Substitutability is a standard assumption in the literature on assortment optimization, which asserts that adding another product to an assortment does not increase the probability of selling other products in the assortment. Mathematically, this is expressed as $\pi \rpr{i \mid S} \geq \pi \rpr{i \mid S \cup \{j\}} $, $\forall i \in S, j \notin S, S \subseteq \cal N$. A very general choice model called the Random Utility model satisfies this assumption (including MNL, Nested Logit, Markov Chain, Mixture of MNLs, and rank-based choice models). However, this property does not hold under opaque selling. The following example illustrates this violation.
\begin{example}
     Suppose $k=1$ and there are 3 products with $v_1 = 4, v_2 = 5, v_3 =3 $ and $r_1 = 5.5, r_2 = 4, r_3 = 6$. For assortments $S_1 = \{1,2\}$ and $S_2 =\{1,2,3\}$, the optimal opaque prices are computed as $\rho^* (S_1) = 3.826 $ and $ \rho^*(S_2) = 3.729 $.  Using the expressions in Lemma \ref{lem:choice-prob}, the probabilities of choosing product 2 under assortments $S_1$ and $S_2$ are approximately $ \pi \rpr{2 \mid \textsc{Opq}_{(1)}^*\rpr{S_1, \bm r_{S_1}} }  \approx 0.554, 
    \pi \rpr{2 \mid \textsc{Opq}_{(1)}^*\rpr{S_2, \bm r_{S_2}} }  \approx 0.656. $ \Halmos


The rationale is that when product 3 is added, the valuation of the opaque product decreases under the risk-averse assumption. As a result, the opaque product becomes a weaker substitute for product 2, which increases the probability that product 2 is chosen in the larger assortment.
\end{example}
The violation of substitutability implies that all existing (approximation) algorithms designed for RUM models are ineffective in the context of opaque selling model, necessitating the development of new techniques.

\begin{figure}[h!]
\FIGURE{\centering
  \includegraphics[width=0.7\textwidth]{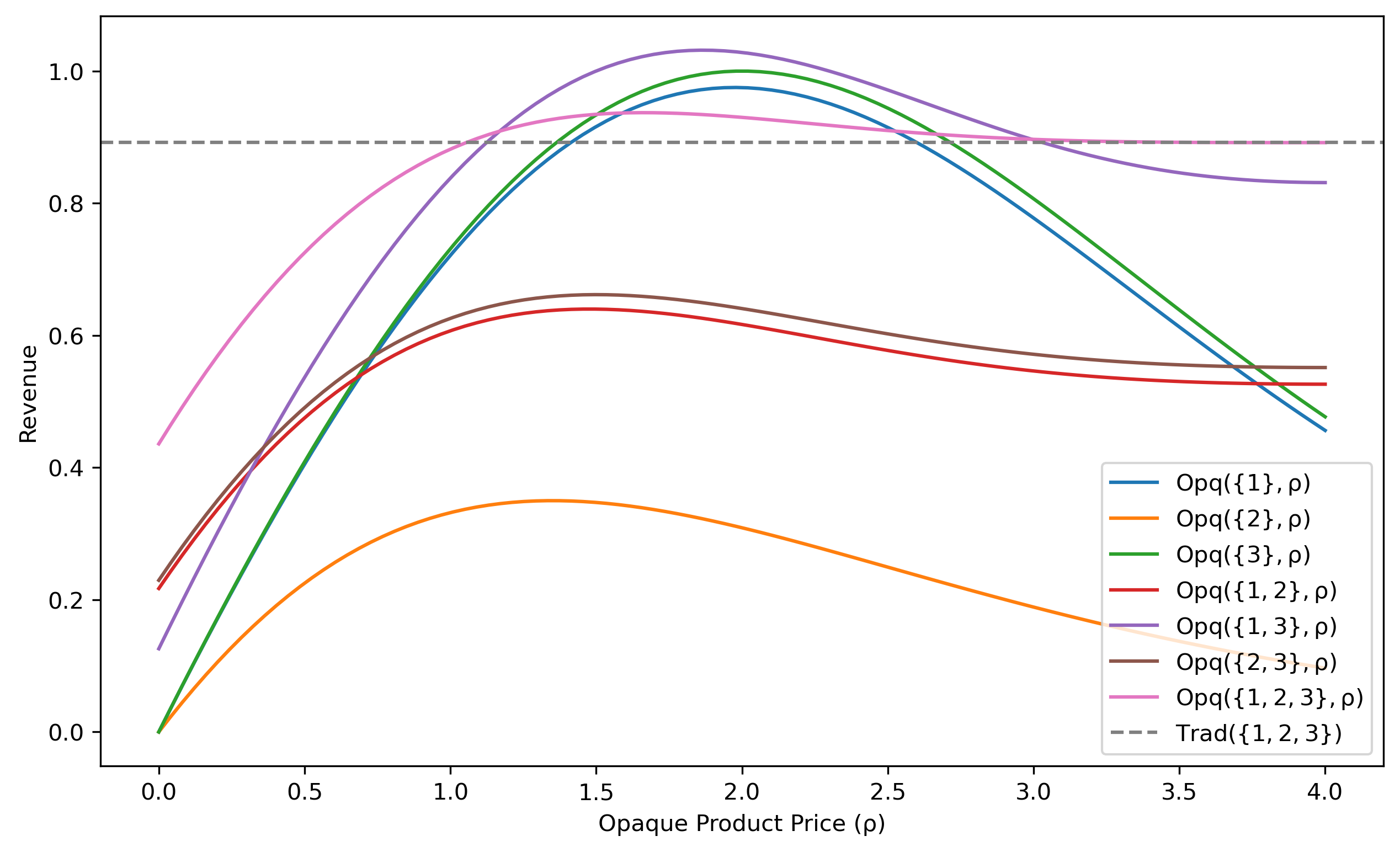} \label{fig:unnested-averse}
}{An example where the optimal assortment under the Opaque-MNL model is not nested-by-revenue. }
{Suppose $k=1$ and there are three products, with $\bm r = (4.02, 4.01, 4)$ and $\bm v = (1.95, 0.3, 2)$. The plots are generated for each assortment by varying the price of the opaque product. The optimal assortment is to offer products $\cpr{1,3}$ along with the opaque product, which is not nested-by-revenue since product $2$ is excluded from the assortment.} 
\end{figure}

{\bf Nested-by-revenue assortments are sub-optimal.}
It is well known that the optimal assortment under the traditional MNL model is \textit{nested-by-revenue}, i.e., offering all products whose price is above a certain threshold (\citealt{talluri2004revenue}). As a result, only $n$ candidate assortments need to be evaluated to identify the optimal assortment. However, this nested-by-revenue property does not hold under the Opaque-MNL model, as indicated by the earlier fact that the optimal assortment may exclude the highest revenue product. 
In Figure \ref{fig:unnested-averse}, we give an example with $k=1$ and 3 products, where it is optimal to offer $\cpr{1,3}$ along with the opaque product and exclude product $2$ from the assortment, even though it has a higher revenue than product $3$. The reason the nested-by-revenue property fails in this context is that, while adding a product with high revenue but low valuation to the assortment does no harm under the traditional MNL model, it can be greatly detrimental for opaque selling, since the valuation of the opaque product is impacted by such a product.

\textbf{The highest revenue product may not be in the optimal assortment.} 
Not only can the optimal assortment fail to be nested-by-revenue, but it may also exclude the product with the highest revenue. This contrasts with standard Random Utility models, where the highest-revenue product is always included in the optimal assortment. Under opaque selling, this property may fail because adding a high-revenue but low-valuation product reduces the valuation of the opaque product, thereby affecting the overall revenue.
\cor{
\begin{example}
     Consider a scenario with $k=2$ and three products, where $v_1 = 1.7, v_2 = 1.8, v_3 =2.4 $ and $r_1 = 6.75, r_2 = 5, r_3 = 4.25$. The optimal assortment is determined to be $S^* = \{2,3\}$, yielding $\rev{\textsc{Opq}_{(2)}^*\rpr{S^*, \bm r_{S^*}}} = 1.459 $ and $\rho^*(S^*) = 2.459  $. While for assortment $S = \{1,2,3\}$, we have $\rev{\textsc{Opq}_{(2)}^*\rpr{S, \bm r_{S}}} = 1.215 $ and $\rho^*(S) = 2.065  $. \Halmos
\end{example}
}

{\bf Evaluating revenue is challenging.} Finally, it is important to note that, although the purchase probabilities admit closed-form expressions, evaluating the expected revenue for even a single assortment takes $O(n2^n)$ operations.  While exact revenue computations are feasible when $n$ is small, we deal with this computational challenge in Lemma \ref{lem:rho-saa} by using Monte Carlo simulations, which provide a more practical alternative for estimating revenues as $n$ becomes larger.

\subsection{Uniformly Priced Traditional Products}\label{sec:assort-uniform}
In this section, we consider an important case where all traditional products are uniformly priced, i.e., $r_i = r$ for $i \in \cal N$. This case is particularly relevant as it represents a natural setting to sell an opaque product when traditional products differ by only one attribute (e.g., notebooks in Figure~\ref{fig:examples}\subref{fig:amazon} differ only in color). 
Additionally, this case is related to the sales maximization problem. However, it differs because the opaque price does not necessarily align with the prices of traditional products. Theorem \ref{thm:nested-val} demonstrates that the optimal assortment in this setting is nested-by-valuation. The proof is in Section \ref{sec:ass_proof}.
\begin{theorem}[Nested-by-valuation]\label{thm:nested-val}
Fix $k \in \{1,2,\ldots,n\}  $ and assume that $r_i = r$ for all $i \in \mathcal{N}$. Let $S^*$ be the optimal assortment of \eqref{opt:assort}. Assume w.l.o.g. that intrinsic valuations are ordered, i.e., $v_1 \geq v_2 \geq \cdots \geq v_n$. Then, there exists \cor{ $j \in \{k, k+1,\ldots,n\}  $ } such that $S^* = \{1,\ldots,j\}$.
\end{theorem}


Recall that in the traditional MNL model, it is optimal to offer all products when prices are equal. However, in the Opaque-MNL model, offering all products may be suboptimal, as it can reduce the valuation of the opaque product and, consequently, its purchase probability.  In fact, the optimal assortment is nested-by-valuation. To demonstrate this structure, we prove that between two assortments of equal cardinality differing by only one product, the assortment with the higher valuation yields greater revenue. Assume w.l.o.g. that the valuations are ordered, i.e., $v_1 \geq v_2 \geq \cdots \geq v_n$. Leveraging this nested-by-valuation structure, we can identify the optimal assortment by exploring only $n$ assortments: $\{1\}, \{1,2\}, \ldots, \{1,2,\ldots,n\}$. For each assortment, we can apply sample-average approximation (Monte Carlo simulation) to jointly approximate the optimal opaque price and the corresponding revenue. The proof of Lemma \ref{lem:rho-saa} is in Appendix \ref{append-same-r}.
\cor{
\begin{lemma}\label{lem:rho-saa}
Let $r_{\max}:=\max_{i\in \mathcal{N}} r_i$ and let parameters $(\epsilon, \delta) \in (0,1)^2 $ be given. For any assortment $S\subseteq \mathcal{N}$, draw $Q$ i.i.d. samples of $ \{V_i\}_{i \in S \cup \{0\} }  $, and let $\hat{\rho}$ maximize the sample-average revenue estimator $\hat{\Rev}^Q \left( \Opq_{(k)} \left(S, \bm r_{S}, \hat{\rho} \right) \right) $. If $ Q = O \left( \frac{r_{\max}^2}{\epsilon^2} \left[ \log \left(\frac{r_{\max}}{\epsilon}\right) +  \log  \left(\frac{1}{\delta}\right)\right] \right) $, then with 
probability at least $1-\delta$,  $  \hat{\Rev}^Q \left( \Opq_{(k)} \left(S, \bm r_{S}, \hat{\rho} \right) \right) \geq \Rev \left( \Opq^*_{(k)} \left(S, \bm r_{S} \right) \right) - \epsilon$. Moreover, both $\hat{\rho}$ and its associated sample-average revenue can be computed in total time $O ( Q |S| + Q\log Q ) $.
\end{lemma} 
}
\cor{
The key here is that the sample-average revenue estimator $\hat{\Rev}^Q \left( \Opq_{(k)} \left(S, \bm r_{S}, \rho \right) \right)$ is piecewise affine in $\rho$ over $[0,r^S_{(k)} ]$ with at most $Q+2$ breakpoints. Thus, one can compute the $\hat{\rho}$ exactly by evaluating $  \hat{\Rev}^Q\left( \Opq_{(k)} \left(S, \bm r_{S}, \rho \right) \right) $ at these breakpoints. Finally, when we restrict attention to nested-by-valuation candidate assortments $S \subseteq \{ \{1\}, \{1,2\}, \ldots, \{1,2,\ldots,n\} \} $, we can generate $Q$ i.i.d. samples of $\{ V_i\}_{i=0}^n $ only once and then apply these samples to evaluate each candidate’s sample-average optimal opaque price and its revenue estimate. Those revenue estimates then allow us to identify an assortment whose revenue is close to the best one. }

\subsection{Proof of Theorem~\ref{thm:nested-val}}  \label{sec:ass_proof} 
We now proceed to prove Theorem \ref{thm:nested-val}. For $k \geq |S| $, since $\rho \leq r$, \eqref{eq:rev-k>|S|} shows that the opaque revenue only contains one term:
\begin{align*}
    \rev{\Opq^*_{(k)}(S, r \ones_S )} = \max_{ \rho \leq r }  \Rev \left( \Trad\left(S, \rho \ones_S \right) \right).
\end{align*}
We also observe that for fixed $\rho$, assortment $S \cup \{j\} $ always yields a higher revenue than $ S $, i.e., 
\begin{align*}
  &  \Rev \left( \Trad\left(S \cup \{j\} , \rho \ones_{S \cup \{j\} } \right) \right) - \Rev \left( \Trad\left(S, \rho \ones_S \right) \right) \\
    = &\frac{ \rho \sum_{i \in S \cup \{j\} }e^{v_i - \rho}  }{ 1 + \sum_{i \in S \cup \{j\} }e^{v_i - \rho} } - \frac{ \rho \sum_{i \in S  }e^{v_i - \rho}  }{ 1 + \sum_{i \in S  }e^{v_i - \rho} } =\frac{\rho e^{v_j - \rho}}{(1 + \sum_{i \in S \cup \{j\} }e^{v_i - \rho})(1 + \sum_{i \in S  }e^{v_i - \rho} )} \geq 0.
\end{align*}
Moreover, let $S_m = S \cup \{m \}  $ and $S_l = S \cup \{l\}$ where  $v_m> v_l$ for $m,l \notin S \subseteq \mathcal{N}$. The assortment with higher valuation yields higher revenue:
\begin{align*}
  &  \Rev \left( \Trad\left(S_m, \rho \ones_{S_m } \right) \right) - \Rev \left( \Trad\left(S_l, \rho \ones_{S_l} \right) \right) \\
    = &\frac{ \rho \sum_{i \in S_m }e^{v_i - \rho}  }{ 1 + \sum_{i \in S_m }e^{v_i - \rho} } - \frac{ \rho \sum_{i \in S_l  }e^{v_i - \rho}  }{ 1 + \sum_{i \in S_l  }e^{v_i - \rho} } =\frac{\rho (e^{v_m- \rho} - e^{v_l -\rho} )}{(1 + \sum_{i \in S_m  }e^{v_i - \rho})(1 + \sum_{i \in S_l  }e^{v_i - \rho} )} \geq 0.
\end{align*}
Thus, when $|S| \leq k$, since products are ordered by valuations, offering assortment $ \{1,2,\ldots,k\} $ maximizes the expected revenue.  

For $ S \subseteq \mathcal{N} $ with $|S| > k$, our proof relies on a key observation: the revenue under the opaque selling model is a convex combination of the revenues from the traditional MNL model with uniform prices $r$ and $\rho$. Recall that $q$ is the opaque product.

\begin{lemma} \label{lem:same_r}
If $r_i = r$, $i \in \mathcal{N}$, then for any assortment $S$, every opaque price $\rho$, and any $k\in\{1,\ldots,|S|\}$, it holds that
\begin{align}\label{eq:opq-tra-r-rho}
\Rev \left(\Opq_{(k)}(S,r \ones_S,\rho)\right)
= \theta_k \cdot \Rev(\Trad(S, \rho \ones_S)) + (1-\theta_k) \cdot \Rev(\Trad(S,r \ones_S)),
\end{align}
where $ \theta_k := \frac{\pi  \left(q \middle| \Opq_{(k)}(S,r\ones_S,\rho)\right)}{\sum_{i\in S}\pi \left(i \middle| \Trad(S,\rho \ones_S)\right)} \in[0,1]$
\end{lemma}


The proof of Lemma~\ref{lem:same_r} is provided in Appendix~\ref{append-same-r}. 
Lemma \ref{lem:same_r}  allows us to either determine the optimal opaque price (in case 1 and case 2) or narrow down the range in which the optimal opaque price lies (in case 3), as shown later. 
To prove Theorem \ref{thm:nested-val}, we establish a sufficient condition: between two assortments differing by only one product, the assortment with the higher valuation yields greater revenue. Specifically, let $S_m = S \cup \{m \}  $ and $S_l = S \cup \{l\}$ where  $v_m> v_l$ for $m,l \notin S \subseteq \mathcal{N}$. Then, for any $ r > 0$, it is always beneficial to offer assortment $S_m$, i.e.,
\begin{align}\label{eq: S_k>S_l}
    \Rev \left(\Opq_{(k)}^*(S_m,r\ones_{S_m})\right)
     \geq    \Rev \left(\Opq_{(k)}^*(S_l,r\ones_{S_l})\right).
\end{align}


We now prove inequality \eqref{eq: S_k>S_l}. Let $r^*_{S_m}$ and $r^*_{S_l}$ be the unique optimal price of $\rev{\trad{S_m}{\cdot  }}$ and $\rev{\trad{S_l}{\cdot  }}$, respectively. Recall that $r_{S_m}^*$ can be solved by $1 + \sum_{i \in S_j}e^{v_i - r_{S_j}^*} = r_{S_j}^* $ for $j = m,l$. Thus, it is clear that $r^*_{S_m} > r^*_{S_l}$. Moreover, $\rev{\trad{S_j}{\cdot}}$ is increasing in $(0,r_{S_j}^*)$ and decreasing in $(r_{S_j}^*, \infty) $  for $j=m,l$. We divide the proof into three cases based on the value of $r$. \\

\underline{\textit{Case 1: $ r \leq r^*_{S_l}$}.} In this case, we know that $\rev{\trad{S_m}{r\ones_{S_m}}} \geq \rev{\trad{S_m}{\rho\ones_{S_m}}}  $ and $\rev{\trad{S_l}{r\ones_{S_l}}} \geq \rev{\trad{S_l}{\rho\ones_{S_l}}} $ since $\rho \leq r$. Then, Eq. \eqref{eq:opq-tra-r-rho} implies that for any $\rho \leq r$, 
\begin{align*}
  \Rev \left(\Opq_{(k)}(S_m,r\ones_{S_m}, \rho)\right) 
\leq \rev{\trad{S_m}{r\ones_{S_m}}}, \   \Rev \left(\Opq_{(k)}(S_l,r\ones_{S_l}, \rho)\right) \leq \rev{\trad{S_l}{r\ones_{S_l}}},
\end{align*}
and the inequalities are attained when $\rho  = r$. Hence, it holds that 
\begin{align*}
  \Rev \left(\Opq_{(k)}^*(S_m,r\ones_{S_m})\right) = \rev{\trad{S_m}{r\ones_{S_m}}} \geq \rev{\trad{S_l}{r\ones_{S_l}}} =   \Rev \left(\Opq_{(k)}^*(S_l,r\ones_{S_l})\right),
\end{align*}
where the inequality follows from $v_m > v_l$. \\

\underline{\textit{Case 2: $ r^*_{S_l} <  r \leq r^*_{S_m}$}.} Since $r \leq r^*_{S_m} $, the same argument shows that $     \Rev \left(\Opq_{(k)}^*(S_m,r\ones_{S_m})\right)  = \rev{\trad{S_m}{r\ones_{S_m}}}.$ While for assortment $S_l$, we have that
\begin{align*}
    \Rev \left(\Opq_{(k)}^*(S_l,r\ones_{S_l})\right)  \leq& \max \{\rev{\trad{S_l}{r\ones_{S_l}}}, \rev{\trad{S_l}{\rho\ones_{S_l}}}  \} \\
    \leq & \rev{\trad{S_l}{r_{S_l}^*\ones_{S_l}}} \\
    \leq & \rev{\trad{S_m}{r_{S_l}^*\ones_{S_m}}} \\
    \leq& \rev{\trad{S_m}{r\ones_{S_m}}} =   \Rev \left(\Opq_{(k)}^*(S_m,r\ones_{S_m})\right).
\end{align*}
The first inequality follows from the convex combination \eqref{eq:opq-tra-r-rho}. The second inequality holds since $r_{S_l}^*$ is the optimal price. The third inequality follows from $v_m > v_l$. The last inequality holds since $r_{S_l}^* < r \leq r_{S_m}^* $ and $\rev{\trad{S_m}{\cdot}}$ is increasing in $(0,r_{S_m}^*)$. \\

\underline{\textit{Case 3: $ r>r^*_{S_m}$}.} Because $r$ is larger than the optimal uniform price, the continuity and unimodality of the revenue under the traditional MNL model indicates that there exist $\rho_l < r^*_{S_l}$ and $\rho_m < r^*_{S_m}$ such that $ \rev{\trad{S_l}{r\ones_{S_l}}} = \rev{\trad{S_l}{\rho_l\ones_{S_l}}}$ and $  \rev{\trad{S_m}{r\ones_{S_m}}} = \rev{\trad{S_m}{\rho_m\ones_{S_m}}}$. Moreover, for all $ \rho_j<\rho < r$, it holds that $\rev{\trad{S_j}{r\ones_{S_j}}} < \rev{\trad{S_j}{\rho\ones_{S_j}}} $, $j = l,m$. Combining this observation and Eq. \eqref{eq:opq-tra-r-rho}, we conclude that for $j=m,l$:
\begin{align*}
  \Rev \left(\Opq_{(k)}^*(S_j,r\ones_{S_j})\right)  = \max_{\rho \leq r}   \Rev \left(\Opq_{(k)}(S_j,r\ones_{S_j}, \rho )\right) = \max_{ \rho_j \leq \rho \leq r}   \Rev \left(\Opq_{(k)}(S_j,r\ones_{S_j} , \rho)\right).
\end{align*}
To prove that $   \Rev \left(\Opq_{(k)}^*(S_m,r\ones_{S_m})\right) \geq   \Rev \left(\Opq_{(k)}^*(S_l,r\ones_{S_l})\right) $, it is sufficient to show that for all $ \rho_l = \min\{ \rho_l, \rho_m\}  \leq \rho \leq r$, assortment $S_m$ yields higher revenue, i.e.,
\begin{align}\label{eq:asst_k>asst_l}
    \Rev \left(\Opq_{(k)}(S_m,r\ones_{S_m}, \rho)\right) \geq   \Rev \left(\Opq_{(k)}(S_l,r\ones_{S_l},\rho)\right) .
\end{align}
Recall that $    \Rev \left(\Opq_{(k)}(S,r\ones_{S}, \rho)\right)=\rho \pi \rpr{q \mid \Opq_{(k)}(S,r\ones_{S}, \rho)} + r \sum_{i \in S} \pi \rpr{i \mid \Opq_{(k)}(S,r\ones_{S}, \rho)}$ by definition. Fix $r$ and $\rho$. For $j=m,l$, let $x_j$, $y_j$ denote $ \pi \rpr{q \mid \Opq_{(k)}(S_j,r\ones_{S_j}, \rho)}$ and $ \sum_{i \in S_j} \pi \rpr{i \mid \Opq_{(k)}(S_j,r\ones_{S_j}, \rho)}$, respectively. It is clear that the customer's purchase probability under assortment $S_m$ is greater than that under assortment $S_l$, i.e., $x_m + y_m > x_l + y_l$
since $S_m = S \cup \{m\}, S_l = S \cup \{l\}$ with $v_m > v_l$. 
We next prove inequality \eqref{eq:asst_k>asst_l} by considering three separate subcases: 

\textit{Subcase 1: $ x_m \leq x_l$.} For this subcase, we have
    \begin{align*}
\Rev \left(\Opq_{(k)}(S_m,r\ones_{S_m}, \rho)\right)  
 = \rho x_m + r y_m =& \rho x_l + r y_l + \rho (x_m - x_l) + r (y_m -y_l) \\
        \stackrel{(a)}\geq & \rho x_l + r y_l + \rho (x_m - x_l) + r (x_l - x_m) \\
        =& \rho x_l + r y_l + (x_l - x_m)(r-\rho) \\ 
        \stackrel{(b)}\geq& \rho x_l + r y_l  =  \Rev \left(\Opq_{(k)}(S_l,r\ones_{S_l}, \rho)\right),
    \end{align*}
    where (a) follows from $x_m + y_m > x_l + y_l $ and (b) holds since $x_m \leq x_l $ and $\rho \leq r$.
    
\textit{Subcase 2:} $ x_m >x_l $ and $  y_m \geq y_l  $. This subcase is trivial since it directly implies that $  \Rev \left(\Opq_{(k)}(S_m,r\ones_{S_m}, \rho)\right)  \geq \Rev \left(\Opq_{(k)}(S_l,r\ones_{S_l}, \rho)\right). $
    
\textit{Subcase 3:} $ x_m >x_l $ and $  y_m < y_l  $. For this subcase, we claim that 
\begin{align}\label{eq:convex-coeff}
    \frac{x_m}{\sum_{i\in S_m}  \ctradequal{i}{S_m}{\rho} } > \frac{x_l}{\sum_{i\in S_l}  \ctradequal{i}{S_l}{\rho} }.
\end{align}
First, we note that Eq. \eqref{eq:opq-tra-r-rho} and the definition of opaque revenue imply that 
   \begin{align*} 
        & \Rev \left(\Opq_{(k)}(S,r\ones_{S}, \rho)\right) = \rho \pi \rpr{q \mid \Opq_{(k)}(S,r\ones_{S}, \rho)} + r \sum_{i \in S} \pi \rpr{i \mid \Opq_{(k)}(S,r\ones_{S}, \rho)}  \\
         = & \frac{\pi \rpr{q \mid \Opq_{(k)}(S,r\ones_{S}, \rho) } }{\sum_{i\in S} \ctradequal{i}{S}{\rho\ones_S} }  \rev{\trad{S}{\rho\ones_{S}}} + \left[1-\frac{\pi \rpr{q \mid \Opq_{(k)}(S,r\ones_{S}, \rho)} }{\sum_{i\in S} \ctradequal{i}{S}{\rho\ones_S} } \right] \rev{\trad{S}{r\ones_S}} \\ 
        =& \rho  \pi \rpr{q \mid \Opq_{(k)}(S,r\ones_{S}, \rho) }  + r \left[1-\frac{\pi \rpr{q \mid \Opq_{(k)}(S,r\ones_{S}, \rho) } }{\sum_{i\in S} \ctradequal{i}{S}{\rho\ones_S} } \right] \sum_{i\in S}  \ctradequal{i}{S}{r\ones_S}.
\end{align*} 
Comparing the above two expressions of opaque revenue, we obtain that 
\begin{equation}\label{eq: prob-identity}
    \sum_{i \in S} \pi \rpr{i \mid \Opq_{(k)}(S,r\ones_{S}, \rho)} = \left[1-\frac{\pi \rpr{q \mid \Opq_{(k)}(S,r\ones_{S}, \rho)}}{\sum_{i\in S} \ctradequal{i}{S}{\rho\ones_S} } \right] \sum_{i\in S}  \ctradequal{i}{S}{r\ones_S}.
\end{equation}
Based on the above equality, we have that
\begin{align*}
     \frac{x_m}{\sum_{i\in S_m}  \ctradequal{i}{S_m}{\rho} } =& 1 - \frac{y_m}{\sum_{i\in S_m}  \ctradequal{i}{S_m}{r\ones_{S_m}}} \\
     \geq & 1 - \frac{y_l}{\sum_{i\in S_m}  \ctradequal{i}{S_m}{r\ones_{S_m}}} \\
     \geq & 1 - \frac{y_l}{\sum_{i\in S_l}  \ctradequal{i}{S_l}{r\ones_{S_l}} } =  \frac{x_l}{\sum_{i\in S_l}  \ctradequal{i}{S_l}{\rho} },
\end{align*}
where the first and last equalities follow from Eq. \eqref{eq: prob-identity}, the first inequality holds by $y_m < y_l$, the second inequality holds since $v_m > v_l$ and thus the purchase probability of $S_m$ is larger than that of $S_l$. Therefore, we have proved inequality \eqref{eq:convex-coeff}. Then, for all $ \rho_l = \min\{ \rho_l, \rho_m\}  \leq \rho \leq r$, it holds that
\begin{align*}
    &\rev{ \Opq_{(k)}(S_m,r\ones_{S_m}, \rho) } \\
        = & \frac{x_m \cdot \rev{\trad{S_m}{\rho\ones_{S_m}}} }{\sum_{i\in S_m}  \ctradequal{i}{S_m}{\rho\ones_{S_m}}   } + \left[1-\frac{x_m}{\sum_{i\in S_m}  \ctradequal{i}{S_m}{\rho\ones_{S_m}} }\right] \rev{\trad{S_m}{r\ones_{S_m}}} \\
        \geq & \frac{x_m \cdot \rev{\trad{S_l}{\rho\ones_{S_l}}}}{\sum_{i\in S_m}  \ctradequal{i}{S_m}{\rho\ones_{S_m}} }  + \left[1-\frac{x_m}{\sum_{i\in S_m}  \ctradequal{i}{S_m}{\rho\ones_{S_m}} }\right] \rev{\trad{S_l}{r\ones_{S_l}}} \\
        \geq& \frac{x_l \cdot \rev{\trad{S_l}{\rho\ones_{S_l}}}}{\sum_{i\in S_l}  \ctradequal{i}{S_l}{\rho\ones_{S_l}} }  + \left[1-\frac{x_l}{\sum_{i\in S_l}  \ctradequal{i}{S_l}{\rho\ones_{S_l}} }\right] \rev{\trad{S_l}{r\ones_{S_l}}} \\
        =& \rev{\Opq_{(k)}(S_l,r\ones_{S_l}, \rho)},
\end{align*}
where the first inequality follows from $v_m > v_l$ and the second inequality follows from \eqref{eq:convex-coeff} and $\rev{\trad{S_l}{\rho}} \geq \rev{\trad{S_l}{r}}$ when $ \rho_l \leq \rho \leq r $. \hfill \Halmos

\subsection{Generally Priced Products}\label{sec:general-assort}
In this section, we study the assortment optimization problem \eqref{opt:assort} with general prices.  Unfortunately, the convex combination \eqref{eq:opq-tra-r-rho} does not hold when prices are not equal. As a result, the proof techniques applicable to uniformly priced products cannot be directly extended to this more general case.

Motivated by the nested-by-valuation structure of the Opaque-MNL model under uniform prices and the nested-by-revenue structure of the traditional MNL model, we propose a natural Nested-by-Revenue-and-Valuation (NRV) heuristic. \cor{Throughout, we assume that products are labeled according to a fixed original index order $1,2,\ldots,n$. For each $i,j\in\{1,\ldots,n\}$, let $S_{r(i)}$ denote the set of the top $i$ products in decreasing revenue order, and let $S_{v(j)}$ denote the set of the top $j$ products in decreasing valuation order, where ties in both rankings are broken lexicographically according to the original product index. }
The NRV heuristic then considers candidate assortments of the form $S_{r(i)}\cap S_{v(j)}$. This heuristic has the advantage of being efficient computationally, since it only evaluates $O(n^2)$ candidate assortments. It is important to note that while assortment $\cpr{1, 3}$ in Figure \ref{fig:unnested-averse} is not nested-by-revenue, it is indeed nested by revenue and valuation. In fact, in our numerical experiments in Section \ref{sec:numerics}, we see that this heuristic is optimal in a wide majority of instances.

One subtle aspect of the NRV heuristic is tie-breaking. In the traditional MNL model, tie-breaks are not necessary, since if $r_i = r_j$ for two products $i$ and $j$, then $j$ is in the optimal assortment if $i$ is in the optimal assortment. Therefore, it is sufficient to select all products whose revenues exceed the threshold. Under Opaque-MNL, however, this simplification no longer applies. Even when products are identical, adding an additional product may either increase or decrease revenue. The NRV heuristic accounts for this feature by fixing a lexicographic tie-breaking rule and evaluating intersections of prefixes in the revenue and valuation rankings, rather than simply selecting all products above given thresholds.

    
\begin{example}[Adding identical products increases revenue]
Consider a setting with $k=1$ and two products where $\bm r = (1, 1)$ and $\bm v = (1, 1)$. For this setting, $\rev{\Opq^*_{(1)}( \{1,2\} )} = 0.67$ and $\rev{\Opq^*_{(1 )}( \{1\} )}  = 0.5$. 
Thus, it is optimal to include both identical products. \Halmos
\end{example}
\begin{example}[Adding identical products decreases revenue] \label{ex:infty-idnt-prod} 
Consider a setting with $k=1$ and two products where $\bm r = (\infty, \infty)$ and $\bm v = (1, 1)$. For this setting, $\rev{\Opq^*_{(1)}( \{1,2\} )}  = 0.34$ and $\rev{\Opq^*_{(1)}( \{1\} )}  = 0.57$, i.e., the optimal assortment contains only one of the identical products. This observation can be extended to $n$ products as well, since when the prices are arbitrarily large, products can only be sold via the opaque product, and the valuation of the opaque product decreases when more products are added to the assortment. \Halmos
\end{example}

The two examples above illustrate why tie-breaking must be handled carefully under Opaque-MNL. We summarize the NRV heuristic below.
\vspace*{0.3in}
\begin{center}
    \fbox{
    \begin{minipage}{0.95\linewidth}
        \vspace*{0.1in}
        \underline{\textbf{Nested by Revenue and Valuation (NRV)}}\\
        \textbf{Input:} $k \in \{1,\ldots,n\}  $, price vector $\bm{r}_{\mathcal{N}}$ and mean valuations $\bm{v} = \cpr{v_i \mid i \in \prodset}$ \\
        \begin{itemize}
            \item  \textbf{Step 1:} Sort products in decreasing revenue order and decreasing valuation order, breaking ties lexicographically by the original product index.
            \item \textbf{Step 2:} Construct candidate assortments $\mathcal{S} = \cpr{S_{v_{(j)}, r_{(i)}} = S_{r(i)} \cap S_{v(j)} \mid i,j \in  \{1,2,\ldots,n\} }$ 
            \item \textbf{Step 3:} Return revenue maximizing assortment $\hat{S}$ from $\mathcal{S}$,
            $$
             \hat{S} = \argmax_{S \in \mathcal{S}} \rev{\Opq^*_{(k)}(S, \bm{r}_S)}.
            $$
        \end{itemize}
    \end{minipage}}
\end{center}
\vspace*{0.3in}

Note that despite performing well in practice, the NRV heuristic is also provably sub-optimal.
\cor{
\begin{example}[NRV heuristic is sub-optimal]
Consider the instance with $k=1$ and three products, where $\bm r=(25,5.5,3) $ and $ \bm v=(1.85,1.8,2.1).$ For this example, the optimal nested-by-revenue-and-valuation assortment is $\{1,2,3\}$ with expected revenue $1.051$. However, the optimal opaque assortment is $\{2,3\}$ with expected revenue $1.073$. In this example, product 1 adds almost no revenue because its price is extremely high, so its direct purchase probability is negligible. At the same time, for $k=1$, adding product 1 does not increase the opaque product’s valuation, so NRV is hurt by being forced to include a product that brings essentially no benefit. \Halmos
\end{example}
}

Next, we demonstrate that NRV is optimal when $n=2$ and achieves a 1/2 approximation ratio when $n \geq 3$. Let $S^*$ denote the optimal assortment under the Opaque-MNL model. Formally, 

\begin{theorem}\label{prop: NRV}
Fix $k \in \{1,2,\ldots,n\} $ and let $\mathcal{S}$ denote the set of candidate assortments considered by the NRV heuristic. 
\begin{itemize}
    \item[(a)] For $n=2$, NRV is optimal, i.e., $S^* \in  \argmax_{S \in \mathcal{S}}   \rev{\Opq_{(k)}^*( S,  \bm r_S)} $ .
    \item[(b)] For $n\geq 3$, NRV has a 1/2 approximation ratio:
    \begin{align*}
    \max_{S \in \mathcal{S}}  \rev{\Opq^*_{(k)}( S,  \bm r_S)}\geq \frac{1}{2}   \rev{\Opq^*_{(k)}(S^*,  \bm r_{S^*})}.
    \end{align*}
\end{itemize}  
\end{theorem}

\cor{
\textit{Proof of Theorem \ref{prop: NRV}.  } (a) Let $\{\hat{i} \}$ denote the best Opaque-MNL singleton assortment, consisting of a traditional product $\hat{i}$ together with an opaque product whose opaque price $\rho$ is optimized, i.e.,
\begin{align*}
     \hat{i} =  \argmax_{ i \in \cal N } \Rev \rpr{ \Opq_{(k)}^* \left( \{i\}, r_i  \right) } = \argmax_{i \in \cal N} \max_{\rho \leq r_i } \frac{\rho e^{v_i - \rho}}{1 + e^{v_i - \rho}} = \argmax_{i \in \cal N} \max_{\rho \leq r_i } \rev{\trad{\cpr{i}}{\rho}}.
\end{align*}
We claim that $\{ \hat{i} \} \in \cal S $. Then,  since both $\mathcal{N}=\{1,2\}$ and the best singleton $\{\hat{i}\} $ are included in $\cal S$, it is clear that NRV is optimal. To this end, we prove it by contradiction. Suppose there exists $i'$ such that $v_{i'} > v_{\hat i}$ and $r_{i'} \geq r_{\hat i}$. It follows that
\begin{align*}
      \Rev \rpr{\Opq_{(k)}^* \rpr{ \{i'\}, r_{i'}  }  }  = \max_{\rho \leq r_{i'}} \frac{\ropq e^{v_{i'} - \ropq}}{1 +  e^{v_{i'} - \ropq}} \geq  \max_{\rho \leq r_{\hat i}} \frac{\ropq e^{v_{i'} - \ropq}}{1 +  e^{v_{i'} - \ropq}} > \max_{\rho \leq r_{\hat i}} \frac{\ropq e^{v_{\hat i} - \ropq}}{1 +  e^{v_{\hat i} - \ropq}}  =   \Rev \rpr{\Opq_{(k)}^* \rpr{ \{\hat{i}\}, r_{\hat{i}}  }  },
    \end{align*}
    where the first inequality follows from the fact $r_{i'} \geq r_{\hat i}$ and the second inequality follows from monotonicity and the fact that $v_{i'} > v_{\hat i}$. This contradicts the claim that $\cpr{\hat i}$ is the optimal singleton assortment. Therefore, we conclude that $\hat{i}$ is the only product with a valuation greater than $v_{\hat{i}}$ and a revenue of at least $r_{\hat{i}}$, meaning that $\cpr{\hat i} \in \mathcal{S}$.  
}

\cor{(b) Next, we prove the 1/2 approximation ratio stated for $n\geq 3$. Let $\hat S $ denote the optimal assortment under the traditional MNL model. Formally,
\begin{align*}
    \hat{S} = \argmax_{ S \subseteq \cal N }  \rev{\trad{S}{\bm r_S  }}.
\end{align*}
For simplicity, we denote $\rho^* :=  \rho^*(S^*).$ Next, recall that the opaque revenue is given by
\begin{align*}
    \Rev \rpr{ \Opq_{(k)}^* \left( S^*, \bm{r}_{S^*} \right) } &= \rho^* \pi \rpr{q \mid \Opq_{(k)}^* \left( S^*, \bm r_{S^*}  \right) } + \sum_{i \in S^*} r_i \pi \rpr{i \mid \Opq_{(k)}^* \left( S^*, \bm r_{S^*}  \right) }.
\end{align*}
We proceed by bounding both terms above using the revenues from candidate assortments in NRV. It is clear that $\hat{S} \in \cal S $ according to the nested-by-revenue structure of the traditional MNL model. We first bound the second term above as follows:
 \begin{align}\label{eq:nested-half-ub-2}
        \sum_{i \in S^*} r_i \pi \rpr{i \mid \Opq_{(k)}^* \left( S^*, \bm r_{S^*}  \right) } \leq&   \sum_{i \in S^*} r_i \pi \rpr{i \mid \Trad \left( S^*, \bm r_{S^*}  \right) }= \rev{\trad{S^*}{\bm r_{S^*}  }}\\
        \nonumber & \leq \rev{\trad{\hat{S}}{\bm r_{\hat{S}}  }} \leq   \Rev \rpr{ \Opq^*_{(k)} \left( \hat{S}, \bm{r}_{\hat{S}} \right) },
    \end{align}
where the first inequality follows from the fact that introducing the opaque product can only decrease the choice probabilities of traditional products, and the second inequality follows from the definition of $\hat S$.
}

\cor{ We next bound the first term (revenue from opaque product) using an NRV assortment.
With a slight abuse of notation, define $\Low = \{ i \in S^*: r_i \geq \rho^* \}$. We also define the nested-by-revenue set induced by $\rho^*$ as  $\mathcal{T} := \{ i\in \mathcal N : r_i \ge \rho^* \}$. Note that $ \Low = S^* \cap \mathcal{T} $. Then, for $j=1,\ldots,n$ , we consider assortments $S_j$ in NRV:
\begin{align*}
 S_j := \mathcal{T}\cap S_{v(j)}, \quad j = 1,\ldots,n ,
\end{align*}
where $S_{v(j)} $ denotes the assortment consisting of the top-$j$ products  by valuation in $\mathcal N$. For any $S_j $, $j=1,\ldots,n$, we provide a useful lower bound on its opaque revenue that will be used later.
\begin{align}
    \Rev \left(\Opq^*_{(k)}(S_j,\bm r_{S_j})\right)   \geq& \Rev \left(\Opq_{(k )}(S_j,\bm r_{S_j},\rho^*)\right) \notag \\
    =& \rho^* \pi \rpr{q \mid \Opq_{(k )} \left( S_j, \bm r_{S_j}, \rho^*  \right) } 
    + \sum_{i \in S_j} r_i \pi \rpr{i \mid \Opq_{(k)} \left( S_j, \bm r_{S_j}, \rho^*  \right) } \notag \\
    \stackrel{(\star)}\geq& \rho^* \left[ \pi \rpr{q \mid \Opq_{(k )} \left( S_j, \bm r_{S_j}, \rho^*  \right) }  + \sum_{i \in S_j} \pi \rpr{i \mid \Opq_{(k )} \left( S_j, \bm r_{S_j}, \rho^*  \right) } \right] \notag \\
    =& \rho^* \left(1 - \pi \rpr{0 \mid \Opq_{(k )} \left( S_j, \bm r_{S_j}, \rho^*  \right) }  \right) \notag \\
    =& \rho^* \left( 1 -  \mathbb{P} \left[ V_0  > V^{S_j}_{(k)} - \rho^* \;\&\; V_0 > V_i -r_i, \forall i \in S_j  \right]  \right) \notag \\
    \geq& \rho^* \left( 1 -  \mathbb{P} \left[ V_0  > V^{S_j}_{(k)} - \rho^*   \right]  \right) \notag \\ 
    =& \rho^*   \mathbb{P} \left[ V^{S_j}_{(k )} - \rho^* > V_0  \right],  \label{ineq: revenue_floor}
\end{align}
where the crucial inequality $(\star)$ follows from $r_i \geq \rho^* $, for all $i \in S_j \subseteq \mathcal{T}$. We next divide the analysis into two cases. 
\textbf{Case (i):} $|\Low| \geq k$. In this case, since $ \Low \subseteq S^* $ and $ |\Low| \geq k $, we observe that $V_{(k)}^{S^*} \leq V_{(k)}^{\Low} $ almost surely because removing elements cannot decrease the $k$-th smallest value. Then, we have that
\begin{align*}
       \pi \rpr{q \mid \Opq_{(k)} \left( S^*, \bm r_{S^*}, \rho^*  \right) } \leq \pi \rpr{q \mid \Opq_{(k)} \left( \Low, \bm r_{\Low}, \rho^*  \right) } \leq&  \mathbb{P} \left[ V_{(k)}^{\Low} - \rho^* > V_0  \right].
\end{align*}
Let $\tilde{j} \in\{1,\ldots,n\} $ be an index such that $S_{\tilde{j}} = \mathcal{T} \cap S_{v(\tilde{j})} $ contains the $|\Low|$ products in $\mathcal{T}$ with the largest valuations $v_i$. Such a $\tilde{j}$ exists because $ \Low = S^* \cap \mathcal{T} \subseteq \mathcal{N} \cap \mathcal{T} = S_n  $. We claim that 
\begin{align*}
     \mathbb{P} \left[ V_{(k)}^{\Low} - \rho^* > V_0  \right]  \leq \mathbb{P} \left[V^{S_{\tilde{j}}}_{(k)}-\rho^* > V_0\right].
\end{align*}
This is because both $\Low$ and $S_{\tilde{j}}$ are subsets of $\mathcal{T}$ with the same cardinality, but $S_{\tilde{j}}$ uses the $|\Low|$ largest mean valuations $v_i$. Since $V_i = v_i + \epsilon_i $ with i.i.d. random variables $\epsilon_i$, replacing products in $\Low$ by those in $S_{\tilde{j}}$ shifts the corresponding valuations upward. Thus, $V_{(k)}^{\Low} $ is stochastically dominated by $V^{S_{\tilde{j}}}_{(k)}$, which yields the inequality above.
Combining the above inequalities and applying \eqref{ineq: revenue_floor} with $j=\tilde{j}$ shows
\begin{align*}
\rho^* \pi \left(q\mid \Opq_{(k)}(S^*,\bm r_{S^*},\rho^*)\right)
&\le \rho^* \mathbb{P}\left[V^{S_{\tilde{j}}}_{(k)}-\rho^*>V_0\right] \leq \Rev \left(\Opq^*_{(k)}(S_{\tilde{j}},\bm r_{S_{\tilde{j}}})\right).
\end{align*}
\textbf{Case (ii):} $|\Low| < k$. If $L=\emptyset$, the opaque product will never be chosen. Therefore, $\pi \left(q\mid \Opq_{(k)}(S^*,\bm r_{S^*}, \rho^*) \right)=0$, and the desired bound holds trivially. We next assume $L\neq\emptyset$ and claim that
\begin{align*}
    \pi \rpr{q \mid \Opq^*_{(k)} \left( S^*, \bm r_{S^*}  \right) }    
    =& \mathbb{P} \left[ V_{(k)}^{S^*} - \rho^* > V_0 \;\&\; V_{(k)}^{S^*} - \rho^* > V_i - r_i, \forall i \in S^*  \right] \\
   =& \sum_{j \in S^* } \mathbb{P} \left[ V_{(k)}^{S^*} - \rho^* > V_0 \;\&\; V_{(k)}^{S^*} - \rho^* > V_i - r_i, \forall i \in S^*,  V_{(k)}^{S^*} = V_j \right] \\
   =& \sum_{j \in \Low } \mathbb{P} \left[ V_{(k)}^{S^*} - \rho^* > V_0 \;\&\; V_{(k)}^{S^*} - \rho^* > V_i - r_i, \forall i \in S^*,  V_{(k)}^{S^*} = V_j \right] \\
  \leq & \sum_{j \in \Low } \mathbb{P} \left[ V_{(k)}^{\Low}  - \rho^* > V_0 \;\&\; V_{(k)}^{\Low}  - \rho^* > V_i - r_i, \forall i \in S^*,  V_{(k)}^{S^*} = V_j \right] \\
 \leq & \mathbb{P} \left[ V_{(k)}^{\Low}  - \rho^* > V_0 \;\&\; V_{(k)}^{\Low}  - \rho^* > V_i - r_i, \forall i \in S^* \right] \\ 
  \leq& \mathbb{P} \left[ V_{(k)}^{\Low} - \rho^* > V_0  \right] .
\end{align*}
The first equality follows from the definition of opaque probability. The second equality holds by the law of total probability. The third equality follows from the fact that for $j \notin \Low $, $ V_{(k)}^{S^*} - \rho^* =  V_j - \rho^* < V_j - r_j $ and thus the corresponding probability is 0. The first inequality is true since for any $j \in \Low $, on the event $ \{V_{(k)}^{S^*} = V_j \} $, we have that $V^{S^*}_{(k)} = V_j \leq \max_{i\in \Low}V_i =V^{\Low}_{(k)}  $ since $|\Low| <k $, and hence replacing $V_{(k)}^{S^*}$ by $V^{\Low}_{(k)}$ can only increase the probability. The second inequality holds as $\{ V^{S^*}_{(k)} = V_j  \}$ are pairwise disjoint almost surely. Finally, dropping constraints increases probability, which yields the third inequality.
Similarly, recall that $S_{\tilde{j}}$ ($|S_{\tilde{j}}| = |\Low|$) consists of the $|\Low|$ largest $v_i$ among $\mathcal{T}$. The same argument in case (i) shows that
\begin{align*}
 \mathbb{P} \left[ V_{(k)}^{\Low} - \rho^* > V_0  \right] \leq  \mathbb{P} \left[V^{S_{\tilde{j}}}_{(k)}-\rho^*>V_0\right].
\end{align*}
Combining the above inequalities and applying \eqref{ineq: revenue_floor} with $j=\tilde{j}$ yields
\begin{align*}
\rho^* \pi \left(q\mid \Opq_{(k)}(S^*,\bm r_{S^*},\rho^*)\right)
&\leq \rho^* \mathbb{P} \left[V^{S_{\tilde{j}}}_{(k)}-\rho^*>V_0\right] \leq \Rev \left(\Opq_{(k)}^*(S_{\tilde{j}},\bm r_{S_{\tilde{j}}})\right).
\end{align*}
Therefore, we conclude that
\begin{align*}
     \Rev \rpr{ \Opq_{(k)}^* \left( S^*, \bm{r}_{S^*} \right) }      \leq& \Rev \left(\Opq_{(k)}^*(S_{\tilde{j}},\bm r_{S_{\tilde{j}}})\right) +  \Rev \rpr{ \Opq^*_{(k)} \left( \hat{S}, \bm{r}_{\hat{S}} \right) } \\
     \leq &  2  \max_{S \in \mathcal{S}}  \rev{\Opq^*_{(k)}( S,  \bm r_S)}
\end{align*}
This completes the proof of part (b) in Theorem \ref{prop: NRV}.  \hfill \Halmos }

\subsection{Numerics} \label{sec:numerics}

We numerically evaluate the NRV heuristic proposed in Section~\ref{sec:general-assort}. For each candidate assortment, the optimized opaque price is computed using the sample-average approximation procedure in Lemma~\ref{lem:rho-saa}, which allows us to evaluate candidate assortments in a computationally tractable and theoretically justified manner. We generate 2000 base instances with 9 products. For each $k\in{1,3,7}$ and each $n=2,\ldots,9$, we construct the corresponding $n$-product instance by retaining the first $n$ products from each base instance. Detailed instance construction and summary tables are reported in Appendix~\ref{appendix:tables}.

The numerical results lead to three main observations. First, the NRV heuristic is nearly optimal across all tested instances. For $k=1$, NRV is optimal when $n=2$, as predicted by Theorem~\ref{prop: NRV}(a), and remains highly accurate for $n\ge 3$: the largest optimality gap across all baseline instances is only $2.498\%$, and the average optimality gap is below $0.013\%$. Second, opaque selling is more likely to improve revenue when traditional prices are relatively high, consistent with Proposition~\ref{prop:condition}. In particular, when we increase the mean of the price distribution, the number of instances in which the opaque product is offered increases relative to the baseline; when we reduce the dispersion of the price distribution, this number decreases substantially. Third, as $k$ increases, opaque selling becomes beneficial in many more instances. This is because customers value the opaque product according to a higher order statistic, and the feasible range for the opaque price expands from $\rho\le r^S_{(1)}$ to $\rho\le r^S_{(k)}$. Across all robustness checks, NRV continues to perform well: the maximum optimality gap remains below $3.8\%$, while the average optimality gap remains below $0.02\%$.

\section{Conclusions}
In this work, we develop a customer choice model for customers facing an opaque product by integrating the opaque product into the traditional MNL choice model. For \cor{$k$-risk-averse customers}, the introduced Opaque-MNL model provides closed-form choice probabilities, and we examine the associated price and assortment optimization problems.
 Remarkably, for the pricing problem, we show that uniform pricing is optimal, implying that adding an opaque product does not enhance revenue when prices are optimized rather than fixed exogenously. It would be interesting to theoretically verify whether this phenomenon extends to other contexts, such as risk-neutral customers, a setting with costs, or entirely different valuation paradigms. 
\cor{ For the joint assortment and opaque support problem, we prove that the optimal assortment is nested-by-revenue.} For the assortment problem with support equal to the offered assortment, we show that the revenue-maximizing assortment is nested by valuation when traditional products share the same price. However, for general prices, the gap between theory and experiments suggests that there is potential for improving the theoretical guarantees of the proposed NRV heuristic. It would also be useful to provide new heuristics with stronger performance guarantees or even characterize the optimal assortment policy for $k$-risk-averse customers.

\FloatBarrier
\bibliographystyle{informs2014} 
\bibliography{refs} 

\newpage
\begin{APPENDICES}

\cor{
\section{Computing the Optimal Opaque Support and Price}\label{app:algorithm}
In this section, for any offered assortment $S$ and $\epsilon>0$, we develop a fully polynomial-time approximation scheme (FPTAS) that computes an opaque support and opaque price whose resulting revenue is within a factor of $(1-\epsilon)$ of $F(S)$:
\begin{align*}
    F(S) = \max_{\substack{T\subseteq S, |T|\le k} } \max_{\rho\ge 0} \Rev \left(\Trad  \left(S,\bm r_S^{T,   \rho\wedge r}\right)\right) = \max_{\rho\ge 0}  \max_{\substack{T\subseteq S, |T|\le k} } \Rev \left(\Trad  \left(S,\bm r_S^{T,   \rho\wedge r}\right)\right) ,
\end{align*}
where the reduction to $|T|\le k$ follows from Lemma~\ref{lem: T* <=k} and Eq.~\eqref{eq:rev-k>|S|}. For notational convenience, we define
\begin{align*}
   G_S(\rho) := \max_{T \subseteq S: |T|\le k} \Rev \left(\Trad \left(S,\bm r_S^{T,\rho\wedge r}\right)\right).
\end{align*}
We first show that for any $\rho$, the corresponding optimal opaque support can be computed in polynomial time. 
}
\cor{
\begin{lemma}\label{lem: fixed-rho-LP}
 For any $\rho \ge 0$, $G_S(\rho)$ can be computed in polynomial time by solving a linear program.
\end{lemma}
}
\cor{
\underline{\textit{Proof of Lemma \ref{lem: fixed-rho-LP}. }} Fix $\rho \ge 0$ and recall that $\Low = \{i \in S: r_i \geq \rho\} $. We reformulate the support optimization problem as an MNL assortment problem under totally unimodular (TU) constraints:
\begin{align}\label{TU LP}
\begin{aligned}
\max_{\bm x,\bm x'} \quad &
\frac{ \sum_{i\in S \setminus \Low } r_i e^{v_i-r_i} x_i + \sum_{i\in \Low} r_i e^{v_i-r_i}x_i +\sum_{i\in \Low} \rho e^{v_i-\rho}x_i'}{1 + \sum_{i\in S \setminus \Low } e^{v_i-r_i}x_i +\sum_{i\in \Low} e^{v_i-r_i}x_i +\sum_{i\in \Low} e^{v_i-\rho}x_i'}\\
\text{s.t.}\quad
& x_i = 1, && \forall i\in S \setminus \Low, \\
& x_i + x_i' = 1, && \forall i\in \Low, \\
& \sum_{i\in \Low } x_i' \le k, \\
& x_i \in \{0,1\}, && \forall i\in S \setminus \Low , \\
& x_i, x_i' \in \{0,1\}, && \forall i\in \Low .
\end{aligned}
\end{align}
To see this, note that any product $i\in S \setminus \Low$ is unaffected by truncation, since $\min\{\rho,r_i\}=r_i$, so it is always included. For each product $i\in \Low$, we create two copies: an original copy with revenue $r_i$ and MNL weight $e^{v_i-r_i}$, and a modified copy with revenue $\rho$ and MNL weight $e^{v_i-\rho}$. The constraint $x_i+x_i'=1$ enforces that exactly one of these two copies is selected, while $\sum_{i\in \Low} x_i'\le k$ ensures that at most $k$ products are switched to their modified copies. Hence, optimizing the opaque support is equivalent to solving the above constrained MNL assortment problem. It remains to show that the constraint matrix is TU. We partition the rows into two parts $\{x_i=1,\ i\in S\setminus \Low\}\cup \{x_i+x_i'=1,\ i\in \Low\}$ and $\{\sum_{i\in \Low}x_i'\le k\}$.  For a column corresponding to $x_i$, it has exactly one nonzero entry. For a column corresponding to $x_i'$ with $i\in \Low$, it has exactly two nonzero entries that have the same sign and lie in different parts of the row partition. By the standard row-partition characterization of total unimodularity, the resulting matrix is TU. Therefore, its LP relaxation is integral. Since MNL assortment optimization under TU constraints can be solved in polynomial time as a linear program \citep{sumida2021revenue}, it follows that $G_S(\rho)$ can be computed in polynomial time. \hfill \Halmos
}
\cor{
Let $\rho^*$ denote the optimal opaque price, i.e., $\rho^*\in \arg\max_{\rho\ge 0} G_S(\rho)$. We next show that $G_S(\rho^*)$ can be approximated efficiently via a geometric grid search. In particular, we evaluate $G_S(\rho)$ at all points in the grid $\{ \frac{\max_{i \in S} r_i }{(1+\epsilon)^m} : m =0,1,\ldots, \left\lceil  \log_{1+\epsilon} \max_{i \in S} r_i \right\rceil  \} \cup \{1\} $ and return the price that attains the largest value. The next lemma shows that this yields the desired approximation guarantee. Since the grid contains $O(\frac{1}{\epsilon} \log \max_{i\in S}r_i )  $ prices, we obtain a FPTAS whose running time is $O(\frac{1}{\epsilon} \log \max_{i\in S}r_i )  $ times the time required to solve the LP relaxation of \eqref{TU LP}.
}
\cor{
\begin{lemma}\label{lem: grid-search-rho}
Fix $\epsilon>0$. We can find an opaque price $\rho$ in polynomial time such that $G_S(\rho) \geq (1- \epsilon) G_S(\rho^*)$.
\end{lemma}
}
\cor{
\underline{\textit{Proof of Lemma \ref{lem: grid-search-rho}. }} We first identify the search interval. Note that for any $\rho \geq \max_{i \in S}r_i$, we have $\bm r_S^{T,\rho\wedge r}=\bm r_S$ for every $T\subseteq S$, and hence $G_S(\rho)=G_S(\max_{i\in S} r_i)$. Thus, an optimal opaque price $\rho^*$ is always upper bounded by $\max_{i\in S} r_i$. To derive the lower bound, we observe that $G_S(\rho)$ is non-decreasing on $[0,1]$.
\begin{claim}\label{claim：G_S increase}
 $G_S(\rho)$ is non-decreasing on $[0,1]$.
\end{claim}
}

\cor{
When $\max_{i \in S}r_i  \leq 1$, Claim \ref{claim：G_S increase} implies that $\rho^* = \max_{i \in S}r_i $ and one can choose opaque price $\rho = 1$ and return $G_S(1).$ We next consider that $\max_{i \in S}r_i  > 1$ and it suffices to search over $[1,\max_{i \in S}r_i]$. We will use the following claim to complete the proof.
\begin{claim}\label{claim: rho1<rho2}
    For any $0 \leq \rho_1 \leq \rho_2$, it holds that $ G_S(\rho_1) \geq \frac{\rho_1}{\rho_2} G_S(\rho_2).$
\end{claim}
}
\cor{
Fix $\epsilon > 0$ and consider a geometric grid on $[1, \max_{i \in S}r_i]$ with multiplier $(1 + \epsilon)$, i.e., $\{ \frac{\max_{i \in S} r_i }{(1+\epsilon)^m} : m =0,1,\ldots, \left\lceil  \log_{1+\epsilon} \max_{i \in S} r_i \right\rceil  \} \cup \{1\} $. The grid contains $O(\frac{1}{\epsilon} \log \max_{i\in S}r_i )  $ prices. For each grid point, we solve the support optimization problem exactly using the TU formulation as shown in Lemma \ref{lem: fixed-rho-LP}. Recall that $\rho^*$ denotes the optimal opaque price. If $\rho^* \leq 1$, we have $G_S(1) \geq G_S(\rho^*) $. Hence, $\rho=1$ is optimal among the candidate prices considered. If $\rho^* >1$, then by construction of the geometric grid, there exists a grid point $\rho$ such that $\rho \leq \rho^* \leq  (1+\epsilon) \rho $. Claim \ref{claim: rho1<rho2} indicates that $G_S(\rho) \geq \frac{\rho}{\rho^*} G_S(\rho^*) \geq \frac{1}{1+ \epsilon} G_S(\rho^*) \geq (1-\epsilon) \cdot  G_S(\rho^*)   $. Thus, in all cases, one of the candidate prices achieves at least a $ (1-\epsilon) $ fraction of the optimal value. \hfill \Halmos
}
\cor{
\subsection{Proofs of Claims \ref{claim：G_S increase} and \ref{claim: rho1<rho2}}
\underline{\textit{Proof of Claim \ref{claim：G_S increase}. }} Fix any feasible support $T \subseteq S$ with $|T|\le k$. On any interval of $\rho$ over which the set $  \{ i \in T : r_i \geq \rho  \}  $ remains unchanged, $\Rev \left(\Trad  \left(S,\bm r_S^{T,\rho\wedge r}\right)\right)$ is differentiable and its derivative is
\begin{align*}
    \frac{d}{d \rho } \Rev \left(\Trad \left(S,\bm r_S^{T,\rho\wedge r}\right)\right) = \frac{  \sum_{i \in \Low} e^{v_i -\rho}  \left[\sum_{i \in S \setminus \Low}r_i e^{v_i- r_i}  + (1-\rho) \cdot (1 + \sum_{i \in S \setminus \Low} e^{v_i- r_i} ) + \sum_{i \in \Low} e^{v_i -\rho}  \right] }{ \left(1 + \sum_{i \in S \setminus \Low} e^{v_i- r_i}   +  \sum_{i \in \Low} e^{v_i -\rho} \right)^2  } \geq 0,
\end{align*}
where the inequality follows from $\rho \le 1$. Since the revenue function is continuous in $\rho$, it follows that $\Rev \left(\Trad \left(S,\bm r_S^{T,\rho\wedge r}\right)\right)$ is non-decreasing in $\rho$ on $[0,1]$. Taking the maximum over all feasible supports $T$ shows that $G_S(\rho)$ is also non-decreasing on $[0,1]$. \hfill \Halmos  \\
}

\cor{
\underline{ \textit{Proof of Claim \ref{claim: rho1<rho2}. }} Fix $0 \leq \rho_1 \leq \rho_2$. Let $T_2$ attain $G_S(\rho_2)$, i.e., $ G_S(\rho_2)  =  \Rev \left(\Trad  \left(S,\bm r_S^{T_2,\rho_2 \wedge r}\right)\right)  $. Since $\min\{r_i,\rho_2\} = r_i $ when $\rho_2 \geq r_i   $, we may without loss of generality assume that $r_i > \rho_2 $ for all $i \in T_2 $. 
}
\cor{
If $T_2 = \emptyset$, then $ G_S(\rho_2)  =  \Rev \left(\Trad  \left(S,\bm r_S \right)\right)  $. Since $\emptyset $ is also feasible for $G_S(\rho_1)$, we have that 
\begin{align*}
    G_S(\rho_1) \geq \Rev \left(\Trad  \left(S,\bm r_S \right)\right)  = G_S(\rho_2) \geq \frac{\rho_1}{\rho_2} G_S(\rho_2).
\end{align*}
}
\cor{
We next assume $T_2 \neq \emptyset$. We first show that $G_S(\rho_2)\le \rho_2. $ Pick any $i \in T_2$, by optimality of $T_2$,
\begin{align*}
 \Rev \left(\Trad \left(S,r_S^{T_2 \setminus\{i\}, \rho_2 \wedge r}\right)\right) \leq \Rev \left(\Trad \left(S,r_S^{T_2,\rho_2 \wedge r}\right)\right) =G_S(\rho_2).
\end{align*}
Using the expressions of MNL revenue for the LHS and $G_S(\rho_2)$, this inequality is equivalent to
\begin{align*}
 \left(G_S(\rho_2 )-\rho_2\right)e^{v_i-\rho_2} +\left(r_i-G_S(\rho_2)\right)e^{v_i-r_i}\leq 0.
\end{align*}
Suppose by contradiction that $G_S(\rho_2)>\rho_2$. Since $e^{v_i-r_i}>0$, dividing by $e^{v_i-r_i}$ yields
\begin{align*}
    \left(G_S(\rho_2)-\rho_2 \right)e^{r_i-\rho_2} + r_i-G_S(\rho_2) \leq 0.
\end{align*}
However, the LHS equals $0$ at $r_i=\rho_2$ and has derivative in terms of $r_i$, $ (G_S(\rho_2)-\rho_2 )e^{r_i-\rho_2} + 1>0. $ Therefore, the LHS is positive since $r_i>\rho_2$ ($i \in T_2$). That's a contradiction. Thus, we obtain that $G_S(\rho_2)\le \rho_2$.
}
\cor{
Now, we use the same support $T_2$ at $\rho_1$. Since $T_2$ is feasible for $G_S(\rho_1)$, it holds that
\begin{align*}
    G_S(\rho_1) \geq \Rev \left(\Trad \left(S,r_S^{T_2 , \rho_1 \wedge r}\right)\right).
\end{align*}
Thus, it suffices to show $ \Rev \left(\Trad \left(S,r_S^{T_2 , \rho_1 \wedge r}\right)\right) \geq \frac{\rho_1}{\rho_2}  G_S(\rho_2) $. Since $r_i > \rho_2 \geq \rho_1$ for all $i \in T_2 $, we note that $ r_S^{T_2 , \rho_1 \wedge r} = r_S^{T_2 , \rho_1 }   $ and $ r_S^{T_2 , \rho_2 \wedge r} = r_S^{T_2 , \rho_2 }$. We observe that
\begin{align*}
    & \Rev \left(\Trad  \left(S,\bm r_S^{T_2, \rho_1 }\right)\right)  - \frac{\rho_1}{\rho_2} G_S(\rho_2) \\
    =& \frac{ \sum_{i \in S \setminus T_2} r_i e^{v_i -r_i} + \sum_{i \in T_2  }\rho_1 e^{v_i - \rho_1} - \frac{\rho_1}{\rho_2} G_S(\rho_2)  (1  + \sum_{i \in S \setminus T_2} e^{v_i -r_i}  + \sum_{i \in T_2  }e^{v_i - \rho_1})   }{1 + \sum_{i \in S \setminus T_2} e^{v_i -r_i}  + \sum_{i \in T_2  }e^{v_i - \rho_1} } \\
    =& \frac{  (1- \frac{\rho_1}{\rho_2}) \cdot \sum_{i \in S \setminus T_2} r_i e^{v_i -r_i} + \frac{\rho_1}{\rho_2} (\rho_2 - G_S(\rho_2) )\sum_{i \in T_2} (e^{v_i -\rho_1} -e^{v_i -\rho_2} ) }{1 + \sum_{i \in S \setminus T_2} e^{v_i -r_i}  + \sum_{i \in T_2  }e^{v_i - \rho_1} }  \geq 0.
\end{align*}
The second equality follows from the fact that $ \sum_{i \in S \setminus T_2} r_i e^{v_i -r_i} - G_S(\rho_2)\cdot(1 + \sum_{i \in S \setminus T_2} e^{v_i -r_i}) = (G_S(\rho_2) -\rho_2 ) \sum_{i \in T_2} e^{v_i - \rho_2}  $. The inequality holds by $ G_S(\rho_2) \leq \rho_2  $ and $e^{v_i - \rho_1} \geq e^{v_i - \rho_2} $. Therefore, we conclude that $  G_S(\rho_1) \geq \frac{\rho_1}{\rho_2} G_S(\rho_2)  $. \hfill \Halmos
}

\section{Numerical Details for the NRV Heuristic}\label{appendix:tables}
In this appendix, we provide the detailed numerical results supporting the summary in Section~\ref{sec:numerics}. For the instance bed, we considered instances with different values of $n$ ranging from 2 to 9. For each product $i$, $r_i$ and $v_i$ are sampled independently from distributions $Lognormal(0.5, 1.5)$ and $Lognormal(0, 0.3)$ respectively. Specifically, we generate 2000 instances of $\{r_1,r_2,\ldots,r_9\}$ and $\{v_1,v_2,\ldots,v_9\}$. For each $n$, we only use the first $n$ prices and valuations, i.e., $\{r_1,r_2,\ldots,r_n\}$ and $\{v_1,v_2,\ldots,v_n\}$. We choose the lognormal distribution for its heavy tail, and choose settings where the price of the product will tend to exceed the valuation to create instances where adding the opaque product is incentivized. 

 We first consider risk-averse customers ($k=1$). In Table \ref{tab:summary}, we summarize the properties of the optimal assortments under Opaque-MNL for our instances, which we denote by $S^*$, as well as the performance of the NRV heuristic on our instance bed. We consider the opaque product to be offered in the optimal assortment only if the optimal opaque price is less than the minimum traditional product prices, i.e., $\rho < r^{S^*}_{(1)}  = \min_{i \in S^*} r_i$.

\begin{table}[htbp]
    \newcolumntype{P}{>{\centering\arraybackslash}c}
    \centering
    \begin{tabular}{c|c|c|c|c|c}
        $n$ & \# instances w/ opaque & \# sub-optimal & max. opt gap (\%) & avg opt gap (\%) & avg size $|S^*|$ \\
        \hline
        2 & 1278 & 0  & 0.0000 & 0.0000 & 1.32 \\
        3 & 1248 & 15 & 2.4978 & 0.0066 & 1.74 \\
        4 & 1110 & 24 & 1.5620 & 0.0068 & 2.25 \\
        5 & 1018 & 41 & 2.0777 & 0.0086 & 2.74 \\
        6 & 918  & 54 & 2.0777 & 0.0094 & 3.28 \\
        7 & 826  & 59 & 1.9671 & 0.0124 & 3.85 \\
        8 & 728  & 64 & 1.9493 & 0.0098 & 4.43 \\
        9 & 672  & 79 & 2.2990 & 0.0083 & 4.96
    \end{tabular}
    \caption{Summary for the optimal assortments under Opaque-MNL when $k=1$.}
    \label{tab:summary}
\end{table}
In Table \ref{tab:summary}, for $n =2$, it is expected for NRV to be optimal as shown in Theorem \ref{prop: NRV}(a).  
For $n \geq 3$, we find that NRV is near-optimal. The largest optimality gap across all instances for NRV is 2.498\%, which is much smaller than the 1/2 approximation ratio. In addition, there is a decreasing trend in the percentage of instances in which the opaque product is offered as $n$ increases. This can be explained as follows. Offering opaque products is generally beneficial when the prices of traditional products are high. To model this, we choose $v\sim Lognormal(0, 0.3)$ and $r \sim Lognormal(0.5,1.5)$, creating instances where prices typically exceed customer valuations. However, numerical experiments show that the revenue from offering opaque products with higher traditional prices is lower than that from the regular MNL model with moderate traditional prices. As the number of products increases, moderate prices appear more, and the optimal assortment will include these moderate-priced products (also narrows the feasible region for the opaque price), where adding an opaque product does not increase revenue. 

Furthermore, although the NRV heuristic becomes less frequently optimal for larger $n$ since $2^n$ (the number of all possible assortments) grows faster than $n^2$ (the number of assortments considered by NRV), the maximum and average optimality gaps do not necessarily increase. This is because as the assortment size increases, adding more products often contributes less to the total revenue due to diminishing returns, which limits the maximum possible revenue difference between the NRV heuristic and the optimal assortment.

\cor{We next vary the price distribution to $Lognormal(1,1.5)$ and $Lognormal(0.5,1)$ to further validate the conditions in Proposition \ref{prop:condition} regarding when offering an opaque product leads to higher revenue. The results are reported in Tables \ref{tab:summary_mean_plus_half} and \ref{tab:summary_var_one}. Table \ref{tab:summary_mean_plus_half} corresponds to the case $Lognormal(1,1.5)$, which generates relatively higher prices on average. Compared with Table \ref{tab:summary}, we observe that the number of instances with opaque products increases for all n from 2 to 9. This is consistent with Proposition \ref{prop:condition}, which suggests that opaque products are more likely to improve revenue when traditional prices are relatively high. By contrast, Table \ref{tab:summary_var_one} reports the results for $Lognormal(0.5,1)$, which generates less dispersed prices (e.g., the 0.99-quantile under $Lognormal(0.5,1.5)$ is 54.03, whereas under $Lognormal(0.5,1)$ it is only 16.88). Thus, it is expected that the number of instances in which opaque selling is beneficial is substantially smaller, decreasing from 1068 when $n=2$ to only 300 when $n=9$, which agrees with Proposition \ref{prop:condition}: when traditional prices are lower, the opaque product is less likely to create additional value. In both cases, the NRV heuristic remains highly effective as the maximum and average optimality gap remains very small.
}

\begin{table}[htbp]
    \newcolumntype{P}{>{\centering\arraybackslash}c}
    \centering
    \cor{
    \begin{tabular}{c|c|c|c|c|c}
        $n$ & \# instances w/ opaque & \# sub-optimal & max. opt gap (\%) & avg opt gap (\%) & avg size $|S^*|$ \\
        \hline
        2 & 1513 & 0   & 0.0000 & 0.0000 & 1.30 \\
        3 & 1470 & 17  & 1.7197 & 0.0058 & 1.69 \\
        4 & 1336 & 44  & 1.7407 & 0.0097 & 2.22 \\
        5 & 1235 & 96  & 2.8752 & 0.0184 & 2.79 \\
        6 & 1116 & 123 & 2.8752 & 0.0191 & 3.43 \\
        7 & 1005 & 133 & 3.5219 & 0.0171 & 4.07 \\
        8 & 917  & 148 & 1.8808 & 0.0134 & 4.71 \\
        9 & 835  & 139 & 1.6723 & 0.0096 & 5.36
    \end{tabular}
    \caption{Summary for the optimal assortments under Opaque-MNL with prices generated from $Lognormal(1, 1.5)$.}
    \label{tab:summary_mean_plus_half}
    }
\end{table}

\begin{table}[htbp]
    \newcolumntype{P}{>{\centering\arraybackslash}c}
    \centering
    \cor{
    \begin{tabular}{c|c|c|c|c|c}
        $n$ & \# instances w/ opaque & \# sub-optimal & max. opt gap (\%) & avg opt gap (\%) & avg size $|S^*|$ \\
        \hline
        2 & 1068 & 0  & 0.0000 & 0.0000 & 1.53 \\
        3 & 883  & 8  & 2.1268 & 0.0040 & 2.15 \\
        4 & 728  & 13 & 2.5867 & 0.0032 & 2.80 \\
        5 & 594  & 18 & 1.2539 & 0.0038 & 3.47 \\
        6 & 506  & 20 & 1.2539 & 0.0025 & 4.10 \\
        7 & 392  & 13 & 0.5850 & 0.0008 & 4.72 \\
        8 & 321  & 15 & 0.5850 & 0.0005 & 5.28 \\
        9 & 300  & 12 & 0.1427 & 0.0003 & 5.80
    \end{tabular}
    \caption{Summary for the optimal assortments under Opaque-MNL with prices generated from $Lognormal(0.5, 1)$.}
    \label{tab:summary_var_one}
    }
\end{table}

\cor{
Finally, we examine the performance of the NRV heuristic for larger order statistics using the same instances in Table~\ref{tab:summary}, namely $k=3$ and $k=7$ (see Tables~\ref{tab:summary_k3} and \ref{tab:summary_k7}). Compared with the benchmark case $k=1$, the number of instances in which offering an opaque product is beneficial becomes substantially larger. This is consistent with the interpretation of larger $k$ as representing less conservative customers. When $k$ increases, the opaque product is no longer valued according to the worst product in the assortment, but instead according to a higher order statistic, so its value becomes significantly higher. At the same time, the feasible price region for the opaque product also becomes larger, since the seller can choose $\rho$ up to $r_{(k)}^S$ instead of $r_{(1)}^S = \min_{i \in S } r_i$. As a result, the opaque product becomes attractive in many more instances. A similar effect explains why the number of instances with the opaque product also increases as $n$ grows for fixed $k$. With more products available, the seller has greater flexibility to construct an assortment whose $k$-th smallest valuation remains high while still maintaining a sufficiently large $k$-th lowest price. Moreover, the NRV heuristic consistently performs well and is optimal in the vast majority of instances.
}

\begin{table}[htbp]
    \newcolumntype{P}{>{\centering\arraybackslash}c}
    \centering
    \cor{
    \begin{tabular}{c|c|c|c|c|c}
        $n$ & \# instances w/ opaque & \# sub-optimal & max. opt gap (\%) & avg opt gap (\%) & avg size $|S^*|$ \\
        \hline
        2 & 1497 & 0  & 0.0000 & 0.0000 & 1.61 \\
        3 & 1712 & 0  & 0.0000 & 0.0000 & 2.20 \\
        4 & 1836 & 7  & 2.7053 & 0.0041 & 2.60 \\
        5 & 1890 & 9  & 2.7053 & 0.0035 & 2.84 \\
        6 & 1924 & 18 & 3.4438 & 0.0079 & 3.00 \\
        7 & 1932 & 33 & 2.7053 & 0.0087 & 3.16 \\
        8 & 1916 & 51 & 2.4450 & 0.0094 & 3.34 \\
        9 & 1894 & 68 & 3.7969 & 0.0113 & 3.54
    \end{tabular}
    \caption{Summary for the optimal assortments under opaque selling when $k=3$.}
    \label{tab:summary_k3}
    }
\end{table}

\begin{table}[htbp]
    \newcolumntype{P}{>{\centering\arraybackslash}c}
    \centering
    \cor{
    \begin{tabular}{c|c|c|c|c|c}
        $n$ & \# instances w/ opaque & \# sub-optimal & max. opt gap (\%) & avg opt gap (\%) & avg size $|S^*|$ \\
        \hline
        2 & 1497 & 0 & 0.0000 & 0.0000 & 1.61 \\
        3 & 1712 & 0 & 0.0000 & 0.0000 & 2.20 \\
        4 & 1841 & 0 & 0.0000 & 0.0000 & 2.74 \\
        5 & 1909 & 0 & 0.0000 & 0.0000 & 3.25 \\
        6 & 1949 & 0 & 0.0000 & 0.0000 & 3.74 \\
        7 & 1971 & 0 & 0.0000 & 0.0000 & 4.22 \\
        8 & 1981 & 0 & 0.0000 & 0.0000 & 4.69 \\
        9 & 1984 & 3 & 1.2650 & 0.0008 & 5.13
    \end{tabular}
    \caption{Summary for the optimal assortments under opaque selling when $k=7$.}
    \label{tab:summary_k7}
    }
\end{table}

\section{Missing Proofs in Section~\ref{sec:price-support} } \label{app:pricing-general} 

\cor{
\underline{ \textit{Proof of Lemma \ref{lemma: k - k-j}. }} Recall from Lemma \ref{lem:choice-prob} that the opaque revenue is given by
\begin{align*}
    \Rev \left( \Opq_{(k)}(T, S, \bm r_S,\rho) \right) =& \sum_{\substack{I\subseteq T\\ |I|\ge k}}  (-1)^{|I|-k}\binom{|I|-1}{k-1} \Rev \left( \Trad\left(S, \bm{r}_S^{I,\rho \wedge r} \right) \right),  \\
    =& \sum_{\substack{I\subseteq T\\ |I|\ge k}}  (-1)^{|I|-k}\binom{|I|-1}{k-1} \Rev \left( \Trad\left(S, \bm{r}_S^{I \cap \Low ,\rho } \right) \right) ,
\end{align*}
where the second equality holds by the definition of $\Low = \{ i \in T, r_i \geq \rho \}  $. Because $ r^T_{(j)} < \rho \leq r^T_{(j+1)} $, it holds that $ | T  \setminus \Low| =j $. 
}
\cor{
Consider any fixed $J \subseteq \Low$ with $|J|=a$. The terms in the summation that share the same price vector $\bm{r}_S^{J,\rho}$ correspond to sets $I$ of the form $I=J\cup I'$ with $|J\cup I'|\ge k$, where $I' \subseteq T \setminus \Low$. Letting $t := |I'|$, the condition $|I| \ge k$ implies $t \in \{\max(0,k-a), \ldots, j \}$. The combined coefficient in front of $\Rev \left( \Trad\left(S, \bm{r}_S^{J  ,\rho } \right) \right)$  is
\begin{align*}
\mathcal C(a):=\sum_{t=\max(0,k-a)}^{j}
(-1)^{(a+t)-k}\binom{a+t-1}{k-1}\binom{j}{t}.
\end{align*}
Using the convention $\binom{n}{m}=0$ for $m<0$ or $m>n$, we extend the sum to $t = 0,\ldots,j$:
\begin{align*}
\mathcal C(a) &=(-1)^{a-k}\sum_{t=0}^{j}(-1)^t\binom{j}{t}\binom{a+t-1}{k-1} = (-1)^{a-(k-j)}\binom{a-1}{ k-j-1 },
\end{align*}
where the equality follows from the alternating Vandermonde identity: $\sum_{t=0}^{j}(-1)^t\binom{j}{t}\binom{m+t}{r}
=(-1)^j\binom{m}{ r-j }$ for all integers $m,r,j\ge 0$. Therefore, substituting this coefficient back into the revenue expression yields the desired result:
\begin{align*}
    \Rev \left( \Opq_{(k)}(T, S, \bm{r}_S,\rho) \right)  =& \sum_{\substack{J \subseteq \Low \\ |J|\ge k-j}}  (-1)^{|J|-(k-j)}\binom{|J|-1}{k -j -1} \Rev \left( \Trad\left(S, \bm{r}_S^{J,\rho } \right) \right) .  
\end{align*}
This completes the proof. \hfill \Halmos  \\
}
\cor{
{\underline{\textit{Proof of Claim~\ref{claim: Binomial}}.}  }   Fix $i\in T$. Let $a_j = w_j (1-\beta_j) \in (0,1) $ for $j \in T$. Since $1 > \sum_{j \in T} a_j $, we have 
\begin{align*}
\frac{1}{1-\sum_{j\in T\setminus I}a_j}
=\int_{0}^{\infty} e^{-\left(1-\sum_{j\in T\setminus I}a_j\right)t} dt
=\int_{0}^{\infty} e^{-t} e^{t\sum_{j\in T\setminus I}a_j} dt.
\end{align*}
Therefore, the target sum equals
\begin{align*}
    \sum_{\substack{I\subseteq T \setminus \{i\}  \\ |I|\ge k}}\frac{  (-1)^{|I|-k}\binom{|I|-1}{k-1} }{1 - \sum_{j \in T \setminus I} w_j(1-\beta_j)}  =& \int_{0}^{\infty} e^{-t}
\sum_{\substack{I\subseteq T \setminus \{i\} \\ |I|\ge k}}
(-1)^{|I|-k}\binom{|I|-1}{k-1} e^{t \sum_{j\in T\setminus I} a_j} dt \\
=&\int_{0}^{\infty} e^{-t\left(1-\sum_{j\in T} a_j\right)}
\sum_{\substack{I\subseteq T \setminus \{i\} \\ |I|\ge k}}
(-1)^{|I|-k}\binom{|I|-1}{k-1}\prod_{j\in I} e^{-t a_j} dt
\end{align*}
Let $N_t=\sum_{j\in T\setminus \{i\} }X_j$ where $X_j$ are independent Bernoulli with $\mathbb P(X_j=1)=e^{-ta_j} \in [0,1] $.  Let $\{A_j\}_{j\in S}$ be events. The indicator that at least \(k\) of the \(A_j\)'s occur can be written as an inclusion–exclusion sum:
\begin{align}\label{eq:identity}
    \mathbf{1} \left\{\sum_{j\in S} \mathbf{1}_{A_j} \geq k \right\}
= \sum_{\substack{I\subseteq S\\ |I|\ge k}} (-1)^{|I|-k}\binom{|I|-1}{ k-1 } \mathbf{1} \left\{\bigcap_{j\in I} A_j\right\}.
\end{align}
}
\cor{
Then, for every $k\ge 1$, applying \eqref{eq:identity} gives
\begin{align*}
    \sum_{\substack{I\subseteq T \setminus \{i\} \\ |I|\ge k}}
(-1)^{|I|-k}\binom{|I|-1}{k-1}\prod_{j\in I} e^{-t a_j} = \mathbb P(N_t \ge k).
\end{align*}
Because $ 0 \leq \mathbb P(N_t \ge k) \leq  1$, it holds that
\begin{align*}
  0 \leq  \sum_{\substack{I\subseteq T \setminus \{i\}  \\ |I|\ge k}}\frac{  (-1)^{|I|-k}\binom{|I|-1}{k-1} }{1 - \sum_{j \in T \setminus I} a_j}  = \int_{0}^{\infty} e^{-t\left(1-\sum_{j\in T}a_j\right)} \cdot \mathbb P(N_t\ge k) dt
\le \int_{0}^{\infty} e^{-t\left(1-\sum_{j\in T}a_j\right)} dt =\frac{1}{1-\sum_{j\in T}a_j}.
\end{align*} 
This completes the proof. \hfill \Halmos \\
}

\section{Missing Proofs in Section~\ref{sec:joint asst support} } \label{app:asst-general} 

\cor{
{\underline{\textit{Proof of Lemma~\ref{lem: T* <=k}}.}  } Fix an assortment $S$. If $|S| \leq k$, then it is clear that the optimal opaque support $|T^*| \leq k $. In the rest of the proof, assume that $|S | > k$. For any $T \subseteq S $ with $|T| >k $, let $\rho^*$ be the corresponding optimal opaque price, i.e., $\rho^* = \argmax_{\rho \geq 0 }   \rev{\Opq_{(k)}( T , S, \bm{r}_S,\rho )  } $. Note that $ \rho^* \in [0, r_{(k)}^T] $. Recall that $\Low := \{i \in T: r_i \geq \rho^* \}$. We claim the following identity: 
\begin{align*}
    & \sum_{ i \in \Low } \frac{e^{v_i - \rho^*} - e^{v_i -r_i}}{1 + \sum_{j \in S } e^{v_j - r_j}}  \left[    \Rev \left( \Opq_{(k)}( T , S, \bm r_S,\rho^*) \right) -  \Rev \left( \Opq_{(k)}( T \setminus \{i\}  , S, \bm r_S,\rho^*) \right) \right] \\
    =& \Rev \left( \Trad\left(S, \bm r_S \right) \right) -    \Rev \left( \Opq_{(k)}( T , S, \bm r_S,\rho^*) \right).
\end{align*}
We now use the identity to show that one can drop an item from $T$ without decreasing the revenue. Consider first that $r_i = \rho^*$ for all $i \in \Low$. Then the LHS is 0, forcing $\Rev \left( \Trad\left(S, \bm r_S \right) \right) =   \Rev \left( \Opq_{(k)}( T , S, \bm r_S,\rho^*) \right)$. In this situation, we can equivalently take $T' = \emptyset$ without decreasing the achieved revenue. Otherwise, there exists at least one product in $\Low$ with $r_i > \rho^*$. Since $\Rev \left( \Trad\left(S, \bm r_S \right) \right) \leq    \Rev \left( \Opq_{(k)}( T , S, \bm r_S,\rho^*) \right)$,  there exists some $ j \in \Low$ such that 
\begin{align*}
    \Rev \left( \Opq_{(k)}( T \setminus \{ j \} , S, \bm r_S,\rho^*) \right) \geq  \Rev \left( \Opq_{(k)}( T , S, \bm r_S,\rho^*) \right).
\end{align*}
Therefore, we conclude that offering opaque support $T \setminus \{j\} $ with its corresponding optimal opaque price does not decrease the revenue:
\begin{align*}
 \max_{\rho \geq 0 }   \rev{\Opq_{(k)}( T \setminus \{j\} , S, \bm{r}_S,\rho )  } \geq &  \Rev \left( \Opq_{(k)}( T \setminus \{j\} , S, \bm r_S,\rho^*) \right) \\
 \geq&   \max_{\rho \geq 0 }\rev{\Opq_{(k)}( T , S, \bm{r}_S,\rho )  },
\end{align*}
Now, starting from any support $T \subseteq S$ with $|T | >k$, we can remove one product and weakly improve the best achievable revenue. Iterating this argument at most $|T| -k$ times yields a support of size at most $k$ that attains revenue at least as large as the original $T$. Consequently, for the fixed assortment $S$, there exists an optimal support $T^* \subseteq S$ with $|T^*| \leq k $ attaining $F(S)$.
}
\cor{
It remains to prove the claimed identity. Since $\rho^* \in [0, r_{(k)}^T ] $, we apply \eqref{eq:rev-k-order}
and observe that:
\begin{align*}
    & \sum_{ i \in \Low } \frac{e^{v_i - \rho^*} - e^{v_i -r_i}}{1 + \sum_{j \in S } e^{v_j - r_j}}  \left[    \Rev \left( \Opq_{(k)}( T , S, \bm r_S,\rho^*) \right) -  \Rev \left( \Opq_{(k)}( T \setminus \{i\}  , S, \bm r_S,\rho^*) \right) \right] \\
    =& \sum_{ i \in \Low } \frac{e^{v_i - \rho^*} - e^{v_i -r_i}}{1 + \sum_{j \in S } e^{v_j - r_j}} \sum_{\substack{I\subseteq T, i \in I\\ |I|\ge k}}   (-1)^{|I|-k}\binom{|I|-1}{k-1} \Rev\left(\Trad \left(S,\bm r_S^{I,\rho^*\wedge r}\right)\right) \\
  =& \sum_{\substack{I\subseteq T\\ |I|\ge k}} (-1)^{|I|-k}\binom{|I|-1}{k-1}   \Rev\left(\Trad \left(S,\bm r_S^{I,\rho^*\wedge r}\right)\right)  \sum_{i\in I\cap \Low} \frac{e^{v_i-\rho^*}-e^{v_i-r_i}}{1 + \sum_{j \in S } e^{v_j - r_j}} \\
 =& \sum_{\substack{I\subseteq T\\ |I|\ge k}} (-1)^{|I|-k}\binom{|I|-1}{k-1}   \Rev\left(\Trad \left(S,\bm r_S^{I,\rho^*\wedge r}\right)\right) \left[\frac{1 + \sum_{j \in I \cap \Low  } e^{v_j - \rho^*} + \sum_{j \in S \setminus (I \cap \Low  ) } e^{v_j -r_j} }{ 1 + \sum_{j \in S} e^{v_j - r_j}  } - 1 \right] \\
   =& \sum_{\substack{I\subseteq T\\ |I|\ge k}}    (-1)^{|I|-k}\binom{|I|-1}{k-1} \frac{\sum_{j\in I \cap \Low } \rho^* e^{v_j - \rho^* } + \sum_{j\in S \setminus ( I \cap \Low) } r_j e^{v_j - r_j}   }{ 1 + \sum_{j \in S}e^{v_j - r_j}  }  - \Rev \left( \Opq_{(k)}( T , S, \bm r_S,\rho^*) \right) \\
   =& \frac{\sum_{j \in S}r_j e^{v_j -r_j}  }{ 1 + \sum_{j \in S}e^{v_j -r_j} }  -\Rev \left( \Opq_{(k)}( T , S, \bm r_S,\rho^*) \right)  \\
 =& \Rev \left( \Trad\left(S, \bm r_S \right) \right)   -\Rev \left( \Opq_{(k)}( T , S, \bm r_S,\rho^*) \right) .
\end{align*}
Here, the first equality follows from Eq. \eqref{eq:rev-k-order} and the fact that when $|T|>k$, the difference retains exactly those terms indexed by subsets $I\subseteq T$ with $|I|\ge k$ that contain $i$. The second equality holds by exchanging the order of summation. The third equality holds because $  1+\sum_{j\in S}e^{v_j-r_j}+\sum_{i\in I\cap \Low}\bigl(e^{v_i-\rho^*}-e^{v_i-r_i}\bigr)  =1+\sum_{j\in I\cap \Low}e^{v_j-\rho^*}+\sum_{j\in S\setminus (I\cap \Low)}e^{v_j-r_j}$. The fourth equality follows by  $\Rev(\Trad(S,\bm r_S^{I,\rho^*\wedge r}))= \frac{\sum_{j\in I \cap \Low } \rho^* e^{v_j - \rho^* } + \sum_{j\in S \setminus ( I \cap \Low) } r_j e^{v_j - r_j}  }{1 + \sum_{j \in I \cap \Low  } e^{v_j - \rho^*} + \sum_{j \in S \setminus (I \cap \Low  ) } e^{v_j -r_j}}$. and Eq. \eqref{eq:rev-k-order}. The fifth equality holds by the fact that the weighted sum of numerators collapses to $\sum_{j\in S} r_j e^{v_j-r_j}$. To this end, we observe that
\begin{align*}
& \sum_{\substack{I\subseteq T\\ |I|\ge k}}    (-1)^{|I|-k}\binom{|I|-1}{k-1} \frac{\sum_{j\in I \cap \Low } \rho^* e^{v_j - \rho^* } + \sum_{j\in S \setminus ( I \cap \Low) } r_j e^{v_j - r_j}   }{ 1 + \sum_{j \in S}e^{v_j - r_j}  } \\
=&  \sum_{\substack{I\subseteq T\\ |I|\ge k}}    (-1)^{|I|-k}\binom{|I|-1}{k-1} 
\frac{ \sum_{j \in S } r_j e^{v_j - r_j} + \sum_{j \in I \cap \Low  } (\rho^* e^{v_j - \rho^*} - r_j e^{v_j -r_j})  }{ 1 + \sum_{j \in S}e^{v_j - r_j}  } \\
=& \frac{\sum_{j \in S} r_j e^{v_j-r_j}}{1 + \sum_{j \in S}e^{v_j - r_j} }  + \frac{\sum_{ j \in \Low}\left(\rho^* e^{v_j-\rho^*}-r_j e^{v_j-r_j}\right) \sum_{\substack{I\subseteq T, j \in I \\ |I|\ge k}}(-1)^{|I|-k}\binom{|I|-1}{k-1} }{1 + \sum_{j \in S}e^{v_j - r_j}} = \frac{\sum_{j \in S} r_j e^{v_j-r_j}}{1 + \sum_{j \in S}e^{v_j - r_j} },
\end{align*}
where the second equality follows from the identity $\sum_{\substack{I\subseteq T, |I|\ge k}}(-1)^{|I|-k}\binom{|I|-1}{k-1} = 1 $  and exchanging the order of summation, and the final equality is true as $ \sum_{\substack{I\subseteq T, j \in I, |I|\ge k}}(-1)^{|I|-k}\binom{|I|-1}{k-1}=0$ for any $j \in T$ (when $|T| > k$). Finally, the last equality follows from the definition of MNL revenue. \hfill \Halmos \\
}
\cor{
{\underline{\textit{Proof of Claim~\ref{claim:mnl-add-remove}.}  }}
For any $i\in S$, we have
\begin{align*}
\Rev(\Trad(S\setminus\{i\},\bm r_{S\setminus\{i\}})) - \Rev(\Trad(S,\bm r_S)) =& \frac{\sum_{j\in S\setminus\{i\}} r_j e^{v_j-r_j}}{1+\sum_{j\in S\setminus\{i\}} e^{v_j-r_j}} - \frac{\sum_{j\in S} r_j e^{v_j-r_j}}{1+\sum_{j\in S} e^{v_j-r_j}} \\
=& \frac{e^{v_i-r_i}\left(\Rev(\Trad(S,\bm r_S)) - r_i\right) }{1+\sum_{j\in S\setminus\{i\}} e^{v_j-r_j} },
\end{align*}
which is nonnegative whenever $r_i \le \Rev(\Trad(S,\bm r_S))$. By the same calculation, if $j \notin S$ and $ r_j \geq \Rev(\Trad(S,\bm r_S))$, then adding $j$ weakly increases revenue: $ \Rev(\Trad(S\cup\{j\},\bm r_{S\cup\{j\}})) \ge \Rev(\Trad(S,\bm r_S)).$ \hfill \Halmos 
}

 \section{Missing proofs in Section~\ref{sec:asst = support} }\label{append-same-r}
\cor{
{\underline{\textit{Proof of Proposition~\ref{prop:condition}}.}  } By setting $ I = \emptyset $ and $T=S$ in the decomposition formula \eqref{eq:decomposition}, we have that
\begin{align*}
    \Rev \left( \Trad\left(S, \bm{r}_S \right) \right) =  \Rev \left( \Trad\left(S, \bm{r}_S^{S  ,\rho } \right) \right) + \sum_{i \in S }\frac{w_i f(\beta_i, \rho, \bm{r}_{S \setminus S}) }{1 - \sum_{j \in S } w_j(1-\beta_j) }.
\end{align*}
 To determine the relationship between the opaque and traditional revenues, we analyze the sign $f(\beta_i, \rho, \bm{r}_{S \setminus S})$. Recall that $\beta_i = e^{\rho - r_i}$. Substituting the explicit forms into the definition of $f$ gives
  \begin{align*}
        f(\beta_i, \rho, \bm{r}_{S \setminus S}) :=& -\ln \beta_i \cdot \beta_i + \left[\rho -\Rev \left( \Trad\left(S, \bm{r}_S^{S  ,\rho } \right) \right)   \right] \cdot (\beta_i -1). \\
        =& (r_i - \rho) e^{\rho - r_i} + \left( \rho - \frac{\rho \sum_{j \in S} e^{v_j - \rho}  }{1 + \sum_{j \in S} e^{v_j - \rho} } \right) (e^{\rho - r_i } -1) \\
        =& \frac{r_i e^{\rho - r_i} - \rho + (r_i - \rho) \cdot \sum_{j \in S} e^{v_j -r_i}  }{1 + \sum_{j \in S} e^{v_j - \rho}} := \frac{h(\rho)}{1 + \sum_{j \in S} e^{v_j - \rho}  }.
    \end{align*}
    Since the denominator is strictly positive, the sign of $f$ depends entirely on $h(\rho)$. 
   Taking the derivative of $h(\rho)$ gives $h'(\rho)= r_i e^{\rho - r_i}  -1 - \sum_{j \in S }e^{v_j - r_i}$.  Recall that the optimal uniform price $r_S^*$ can be solved by $ 1 + \sum_{j \in S} e^{v_j - r_S^*} = r_S^* .$ This allows us to classify the behavior of $h(\rho)$ based on the relationship between $r_i$ and $r_S^*$.
}
\cor{
   \textbf{Case (a):} $ r_i \leq r_S^* $, for all $ i \in S$ and $\rho \leq \min_{i \in S} r_i $.  In this case, $r_i \leq 1+ \sum_{j \in S}e^{v_j -r_i }$. Since $\rho \leq r_i$, we have $e^{\rho - r_i} \le 1$. Therefore,
\begin{align*}
    h'(\rho) \leq h'(r_i) = r_i -  1 - \sum_{j \in S}e^{v_j -r_i } \leq 0.
\end{align*}
Thus, $h(\rho)$ is non-increasing. Given that $h(r_i) = 0$, it follows that $h(\rho) \geq 0$ for all $\rho \leq r_i$. Consequently, $f(\beta_i, \rho, \bm{r}_{S \setminus S}) \geq 0$ for all $i \in S$. For $k \leq |S| -1$, applying Claim \ref{claim: Binomial}, the coefficient of the $f$ term in the opaque revenue \eqref{eq:opq-rev-simplifies} is bounded by the traditional revenue coefficient:
\begin{align*}
     \Rev \left( \Opq_{(k)}( S,\bm{r}_S,\rho) \right) =&  \Rev \left( \Trad\left(S, \bm{r}_S^{S  ,\rho } \right) \right)  + \sum_{i \in S}w_i  f(\beta_i, \rho, \bm{r}_{S \setminus S})  \sum_{\substack{I\subseteq S \setminus \{i\}  \\ |I|\ge k}}\frac{  (-1)^{|I|-k}\binom{|I|-1}{k-1} }{1 - \sum_{j \in S \setminus I} w_j(1-\beta_j)}  \\
    \leq & \Rev \left( \Trad\left(S, \bm{r}_S^{S  ,\rho } \right) \right)  + \sum_{i \in S} \frac{w_i  f(\beta_i, \rho, \bm{r}_{S \setminus S})}{1 - \sum_{j \in S} w_j (1 - \beta_j) }   =   \Rev \left( \Trad\left(S, \bm{r}_S \right) \right) .
    \end{align*}
For $k \geq |S| $, applying Eq. \eqref{eq:rev-k>|S|} and the fact that the coefficient is non-negative shows that
\begin{align*}
     \Rev \left( \Opq_{(k)}( S,\bm{r}_S,\rho) \right) =  \Rev \left( \Trad\left(S, \bm{r}_S^{S  ,\rho } \right) \right)  \leq&  \Rev \left( \Trad\left(S, \bm{r}_S^{S  ,\rho } \right) \right)  + \sum_{i \in S} \frac{w_i  f(\beta_i, \rho, \bm{r}_{S \setminus S})}{1 - \sum_{j \in S} w_j (1 - \beta_j) }  \\
     =&  \Rev \left( \Trad\left(S, \bm{r}_S \right) \right) .
    \end{align*}
}
\cor{
 \textbf{Case (b):}  $ r_i > r_S^* $, for all $ i \in S$ and $\rho \in ( \max_{j \in S} \underline{\rho}_j , \min_{i \in S} r_i)  $ with $h(\underline{\rho}_i) = 0 $ and $\underline{\rho}_i \in (0,r_i)$. In this case, $r_i > 1+ \sum_{j \in S}e^{v_j -r_i }$ . We observe that $h(r_i) = 0$ and its derivative at $r_i$ is positive: 
 \begin{align*}
    h'(r_i) = r_i - \left( 1 + \sum_{j \in S}e^{v_j -r_i } \right) > 0. 
 \end{align*}
Furthermore, computing the second derivative yields $h''(\rho) = r_i e^{\rho - r_i} > 0$ for all $\rho$, implying that $h(\rho)$ is strictly convex. Since $h(\rho)$ is convex, $h(r_i) = 0$, and $h'(r_i) > 0$, it follows that $h(\rho)$ must be strictly negative in an interval immediately to the left of $r_i$. Additionally, since $h(0) > 0$, the strict convexity implies there exists a unique root $\underline{\rho}_i \in (0, r_i)$ such that $h(\underline{\rho}_i) = 0$. Consequently, $h(\rho) < 0$ for all $\rho \in (\underline{\rho}_i, r_i)$, which implies $f(\beta_i, \rho, \bm{r}_{S \setminus S}) < 0$ in this interval. Then, for any $\rho \in (\max_{j \in S} \underline{\rho}_j, \min_{i \in S} r_i )$, we have $f(\beta_i, \rho ,\bm{r}_{S \setminus S}) < 0$ for all $i \in S$.  
}
\cor{
 For $k \leq |S| -1$, since the coefficient of the $f$ term in \eqref{eq:opq-rev-simplifies} is non-negative and smaller than the coefficient in the traditional revenue (Claim \ref{claim: Binomial}), multiplying by the negative factor $f$ reverses the inequality:
     \begin{align*}
         \Rev \left( \Opq_{(k)}( S,\bm{r}_S,\rho) \right) =&  \Rev \left( \Trad\left(S, \bm{r}_S^{S  ,\rho } \right) \right)  + \sum_{i \in S}w_i  f(\beta_i, \rho, \bm{r}_{S \setminus S})  \sum_{\substack{I\subseteq S \setminus \{i\}  \\ |I|\ge k}}\frac{  (-1)^{|I|-k}\binom{|I|-1}{k-1} }{1 - \sum_{j \in S \setminus I} w_j(1-\beta_j)}  \\
        > & \Rev \left( \Trad\left(S, \bm{r}_S^{S  ,\rho } \right) \right)  + \sum_{i \in S} \frac{w_i  f(\beta_i, \rho, \bm{r}_{S \setminus S})}{1 - \sum_{j \in S} w_j (1 - \beta_j) }   =   \Rev \left( \Trad\left(S, \bm{r}_S \right) \right) .
    \end{align*}
For $k \geq |S| $, applying Eq. \eqref{eq:rev-k>|S|} and the fact that the coefficient is non-negative shows that
\begin{align*}
     \Rev \left( \Opq_{(k)}( S,\bm{r}_S,\rho) \right) =  \Rev \left( \Trad\left(S, \bm{r}_S^{S  ,\rho } \right) \right)  > &  \Rev \left( \Trad\left(S, \bm{r}_S^{S  ,\rho } \right) \right)  + \sum_{i \in S} \frac{w_i  f(\beta_i, \rho, \bm{r}_{S \setminus S})}{1 - \sum_{j \in S} w_j (1 - \beta_j) }  \\
     =&  \Rev \left( \Trad\left(S, \bm{r}_S \right) \right) .
    \end{align*}
This completes the proof.  \hfill \Halmos \\
}

\cor{
{\underline{\textit{Proof of Lemma~\ref{lem:rho-saa}.}}  } 
When $|S| \leq k$, since the opaque revenue collapses to one traditional MNL revenue, we don not need Monte Carlo simulations. We assume $ |S| > k $ for the rest of the proof. Consider $Q $ i.i.d. valuation samples $\{V_i^t \}_{i \in S \cup \{0\} }$. We use superscript $t$ to represent quantities of the $t$-th sample: $V_0^t$, $U_i^t = V_i^t  - r_i$ ($i \in S$), $V_{(k)}^t:= k$-th order statistic of $\{V_i^t : i \in S \}$, and $M^t:= \max \{V_0^t, \max_{i \in S} U_i^t \} $. We also define $\tau^t := \min \{ r^S_{(k)}, \max \{0, V_{(k)}^t - M^t \}  \}  $ (i.e., $V_{(k)}^t - M^t$ projected onto $[0,r^S_{(k)} ]$ ) and $r^t = r_{\argmax_{i \in S} U_i^t  } \cdot \ones \{ \max_{i \in S} U_i^t\ge V_0^t\}\in[0,r_{\max}]$ as the price of the traditional product with the highest utility. For a given $\rho$, the opaque utility is $ V^t_{(k)} - \rho $. Therefore, the opaque product is chosen if and only if $ V^t_{(k)} - \rho \geq M^t  \Leftrightarrow \rho \leq V^t_{(k)} - M^t  $. Because we restrict $\rho$ to $[0,r^S_{(k)}]$, this condition is equivalent to $\rho\le\tau^t$. If $\rho\le\tau^t$, the opaque product is chosen and the realized revenue equals $\rho$. 
When the opaque product ties with the best non-opaque alternative, we break ties in favor of the opaque product.
If instead $\rho>\tau^t$, then the opaque product is not chosen and the realized revenue is $r^t$. Therefore, for every $\rho\in[0,r^S_{(k)}]$, the realized revenue from the  $t$-th sample, $Z^t(\rho)$, is given by
\begin{align*}
   Z^t(\rho):=\rho  \ones \{\rho\le\tau^t\}+r^t \ones \{\rho>\tau^t\},\qquad \rho\in[0, r^S_{(k)}].
\end{align*}
We denote  $\hat{\Rev}^Q \left( \Opq_{(k)} \left(S, \bm r_{S}, \rho \right) \right)$ as the sample-average revenue estimator
\begin{align*}
    \hat{\Rev}^Q\left( \Opq_{(k)} \left(S, \bm r_{S}, \rho \right) \right) : = \frac{1}{Q } \sum_{t=1}^Q Z^t(\rho) .
\end{align*}
We first note that $  \hat{\Rev}^Q\left( \Opq_{(k)} \left(S, \bm r_{S}, \rho \right) \right)$ is affine in $\rho$ with breakpoints contained in $ \{0, r^S_{(k)}\} \cup \{\tau^1,\dots,\tau^Q\}$. In fact, let $\tau^{(1)}\le \tau^{(2)}\le \cdots \le \tau^{(Q)}$ be the sorted breakpoints, and set $\tau^{(0)}:=0$ and $\tau^{(Q+1)}:=r^S_{(k)}$. Fix any interval $(\tau^{(j)},\tau^{(j+1)})$. For any $\rho$ in this interval, the indicator $\ones \{\rho\le\tau^t\}$ is constant in $\rho$ for each $t$, hence each $Z^t(\rho)$ is either the affine function $\rho$ (if $\tau^t\ge \rho$ throughout the interval) or the constant $r^t$ (if $\tau^t< \rho$ throughout the interval). Therefore the average $   \hat{\Rev}^Q\left( \Opq_{(k)} \left(S, \bm r_{S}, \rho \right) \right)$ is affine in $\rho$ on $(\tau^{(j)},\tau^{(j+1)})$. 
}
\cor{
An affine function on an interval achieves its maximum over the closed interval $[\tau^{(j)},\tau^{(j+1)}]$ at one of the endpoints. Because the entire feasible set $[0,r^S_{(k)}]$ is the union of these closed intervals, there exists a maximizer, $\hat{\rho} \in \argmax_{\rho \in [0, r^S_{(k)}] } \hat{\Rev}^Q\left( \Opq_{(k)} \left(S, \bm r_{S}, \rho \right) \right)  $, that lies in $ \{0, r^S_{(k)}\} \cup \{\tau^1,\dots,\tau^Q\}$, which has cardinality at most $Q+2$. Thus, one can compute $\hat\rho$ by evaluating $  \hat{\Rev}^Q\left( \Opq_{(k)} \left(S, \bm r_{S}, \rho \right) \right) $ at at most $Q+2$ points. The time complexity is as follows: It takes $O(Q|S|)$ to compute the statistics ($M^t, \tau^t, r^t$) for all $Q$ samples, and $O(Q \log Q)$ to sort the $Q$ breakpoints to find the optimal $\hat{\rho}$. Then, using prefix sum, the evaluation can be done in $O (Q)$.
}

\cor{
We next use pseudo-dimension to establish a uniform high-probability convergence bound and the number of samples ($Q$) required to get a $\epsilon$ additive loss estimator. Define the function class
\begin{align*}
    \mathcal{F}_S := \left\{ f_\rho:(\tau,r)\mapsto \rho \ones\{\rho\le \tau\}+ r\ones\{\rho>\tau\}: \rho\in[0,r^S_{(k)}] \right\}
\end{align*}
and use $\text{Pdim}(\mathcal F_S)$ to denote its pseudo-dimension.  Recall that $ \Rev(\Opq_{(k)}(S,\bm r_S,\rho)) = \mathbb{E} [Z^1(\rho)] = \mathbb{E} \left[f_\rho(\tau^1, r^1)\right]   $ and $ \hat{\Rev}^Q\left( \Opq_{(k)} \left(S, \bm r_{S}, \rho \right) \right) = \frac{1}{Q}\sum_{t=1}^Q f_{\rho} (\tau^t, r^t)  $. Note that the functions in $\mathcal F_S$ take values in $[0, r_{\max} ]$.
Applying Theorem 11.8 in \cite{mohri2018foundations} gives
\begin{align*}
    \sup_{f_{\rho} \in\mathcal F_S} \left| \frac{1}{Q}\sum_{t=1}^Q f_{\rho}(\tau^t, r^t)- \mathbb{E} \left[f_{\rho}(\tau^1, r^1) \right] \right| \leq r_{\max} \sqrt{ \frac{2 \text{Pdim}(\mathcal F_S) \log \frac{e Q}{\text{Pdim}(\mathcal F_S)}  }{Q} } + r_{\max} \sqrt{\frac{\log(1/\delta)}{2Q}}  
\end{align*}
with probability at least $1-\delta$. Applying this with $\text{Pdim}(\mathcal F_S) \leq 2$ yields
\begin{align*}
      \sup_{\rho\in \left[0,r^S_{(k)} \right] } \left|   \hat{\Rev}^Q\left( \Opq_{(k)} \left(S, \bm r_{S}, \rho \right) \right)  - \Rev(\Opq_{(k)}(S,\bm r_S,\rho))   \right|  \leq r_{\max} \sqrt{ \frac{4 \log \frac{e Q}{2}  }{Q} } + r_{\max} \sqrt{\frac{\log(1/\delta)}{2Q}}
\end{align*}
with probability at least $1-\delta$.  Thus, if $ r_{\max} \sqrt{ \frac{4 \log \frac{e Q}{2}  }{Q} } \leq \frac{\epsilon}{2} \Longleftrightarrow Q \geq \frac{16 (r_{\max})^2}{\epsilon^2} \log \frac{eQ}{2}  $ and $ r_{\max} \sqrt{\frac{\log(1/\delta)}{2Q}} \leq \frac{\epsilon}{2} \Longleftrightarrow Q\geq \frac{2 (r_{\max})^2}{\epsilon^2} \log \frac{1}{\delta}  $, then the right-hand side is at most $\epsilon$, which gives the desired uniform bound. In particular, $Q = O \left( \frac{r_{\max}^2}{\epsilon^2}
\left[ \log \frac{r_{\max}}{\epsilon} +  \log  \frac{1}{\delta}\right] \right) $. Let $\rho^*$ denote the optimal opaque price for assortment $S$. Based on the uniform convergence bound and the optimality of $\hat{\rho}$, we conclude that
\begin{align*}
    \hat{\Rev}^Q\left( \Opq_{(k)} \left(S, \bm r_{S}, \hat{\rho} \right) \right) \geq   \hat{\Rev}^Q\left( \Opq_{(k)} \left(S, \bm r_{S}, \rho^* \right) \right)  \geq \Rev(\Opq_{(k)}(S,\bm r_S,\rho^*))  - \epsilon.
\end{align*}
}
\cor{
Finally, we finish the proof by showing that $\text{Pdim}(\mathcal F_S)\le 2$. For any fixed  $((\tau,r),s)$, the superlevel set $\{\rho : f_{\rho} (\tau, r) \geq s \}$ is an interval in $\rho$. Hence the induced threshold labelings correspond to interval classifiers on $\mathbb{R}$, whose VC-dimension is 2 (See \citealp{mohri2018foundations} Ch. 3). Thus,  $\text{Pdim}(\mathcal F_S)\le 2$. \hfill \Halmos \\
}

{\underline{\textit{Proof of Lemma~\ref{lem:same_r}.}}  } 
Fix  $k\in\{1,\ldots,|S|\}$. Since $r_i = r$, $i \in \mathcal{N}$ and $\rho\le r$, we have $\Low = \{i \in S, r_i \geq \rho  \} =S$.    By Lemma~\ref{lem:choice-prob}, the opaque revenue admits
\begin{align*}
\Rev \left(\Opq_{(k)}(S,r\ones_S,\rho)\right)
= \sum_{\substack{I\subseteq S\\ |I|\ge k}} (-1)^{|I|-k}\binom{|I|-1}{k-1}  \Rev \left(\Trad(S,\bm{r}_S^{I,\rho})\right),
\end{align*}
and the opaque purchase probability satisfies
\begin{align}\label{eq: convex-opaque}
\pi \left(q \middle| \Opq_{(k)}(S,r\ones_S,\rho)\right) = \sum_{\substack{I\subseteq S\\ |I|\ge k}}
(-1)^{|I|-k}\binom{|I|-1}{k-1} \sum_{i\in I}\pi \left(i \middle| \Trad(S,\bm{r}_S^{I,\rho})\right).
\end{align}

We first show that each mixed-price MNL revenue is a convex combination of $\Rev(\Trad(S, \rho \ones_S)) $ and $\Rev(\Trad(S,r \ones_S))$. Fix any $I\subseteq S$. Direct computation yields that
\begin{align*}
&\Rev(\Trad(S, \bm{r}_S^{I,\rho})) \\
=& \frac{ \rho \sum_{i \in I} e^{v_i - \rho} + r \sum_{i \in S \setminus I} e^{v_i - r} }{1 + \sum_{i \in I} e^{v_i - \rho} +  \sum_{i \in S \setminus I} e^{v_i - r}} \\
=& \frac{ \rho \sum_{i \in I} e^{v_i - \rho}  }{1 + \sum_{i \in I} e^{v_i - \rho} +  \sum_{i \in S \setminus I} e^{v_i - r}} + \frac{  r \sum_{i \in S \setminus I} e^{v_i - r} }{1 + \sum_{i \in I} e^{v_i - \rho} +  \sum_{i \in S \setminus I} e^{v_i - r}} \\ 
= & \frac{\frac{ \sum_{i \in I} e^{v_i - \rho}  }{1 + \sum_{i \in I} e^{v_i - \rho} + \sum_{i \in S \setminus I}e^{v_i -r} } }{ \frac{\sum_{i\in S} e^{v_i-\rho}}{1+\sum_{i\in S} e^{v_i-\rho}} } \frac{ \rho\sum_{i\in S} e^{v_i-\rho}}{1+\sum_{i\in S} e^{v_i-\rho}} + ( 1- \frac{\frac{ \sum_{i \in I} e^{v_i - \rho}  }{1 + \sum_{i \in I} e^{v_i - \rho} + \sum_{i \in S \setminus I}e^{v_i -r}  }}{ \frac{\sum_{i\in S} e^{v_i-\rho}}{1+\sum_{i\in S} e^{v_i-\rho}} } ) \frac{ r\sum_{i\in S} e^{v_i-r}}{1+\sum_{i\in S} e^{v_i-r}}   \\ 
= &  \frac{\sum_{i\in I}\pi( i\mid \Trad(S,\bm{r}_S^{I,\rho}))}{\sum_{i\in S}\pi( i \mid \Trad(S, \rho \ones_S ))} \Rev(\Trad(S, \rho \ones_S))
  + \rpr{ 1- \frac{\sum_{i\in I}\pi( i\mid \Trad(S,\bm{r}_S^{I,\rho}))}{\sum_{i\in S}\pi( i \mid \Trad(S, \rho \ones_S ))}  } \Rev(\Trad(S,r  \ones_S)) .
\end{align*}
According to \eqref{eq: convex-opaque}, it holds that
\begin{align*}
     \sum_{\substack{I\subseteq S\\ |I|\ge k}} (-1)^{|I|-k}\binom{|I|-1}{k-1} \frac{\sum_{i\in I}\pi( i\mid \Trad(S,\bm{r}_S^{I,\rho}))}{\sum_{i\in S}\pi( i \mid \Trad(S, \rho \ones_S ))} = \frac{\pi \left(q \middle| \Opq_{(k)}(S,r\mathbf 1_S,\rho)\right)}{\sum_{i\in S}\pi( i \mid \Trad(S, \rho \ones_S ))} := \theta_k.  
\end{align*}
Thus, plugging the convex-combination representation of each $\Rev(\Trad(S,\bm{r}_S^{I,\rho}))$ into the opaque revenue shows that
\begin{align*}
\Rev \left(\Opq_{(k)}(S,r\ones_S,\rho)\right)  =\theta_k \cdot \Rev(\Trad(S, \rho \ones_S)) + (1-\theta_k) \cdot \Rev(\Trad(S,r \ones_S)).
\end{align*}

Finally, the convex combination can be proved by showing that $\pi \left(q  \middle|  \Opq_{(k)}(S,r\mathbf 1_S,\rho)\right) \leq \sum_{i\in S} \ctradequal{i}{S}{\rho\ones_S} $, which by definition is equivalent to 
\begin{align*}
    \PP\spr{  V_{(k)}- \rho  >  V_i - r , \forall i \in S ,  V_{(k)} - \rho > U_0} \leq  \PP\spr{ \exists j \in S,  V_j - \rho  >  U_0 } .
\end{align*}
Because the event of the LHS is included in the event of the RHS, the inequality holds.  \hfill \Halmos  \\

    \end{APPENDICES}

\end{document}